\documentclass[11pt,a4paper]{article}
\usepackage{amsmath, amsfonts, amssymb,amsthm}
\usepackage{a4wide}
\usepackage{parskip}
\usepackage{enumitem}
\usepackage{xcolor}
\usepackage{mathrsfs}
\usepackage{cite}
\usepackage{hyperref}

\usepackage{latexsym}
\usepackage{comment}
\usepackage{tikz}

\usepackage[left=1in, right=1in,top=1in,bottom=1in]{geometry}

\allowdisplaybreaks

\newcommand{\be} {\begin{equation}}
	\newcommand{\ee} {\end{equation}}
\newcommand{\bea} {\begin{eqnarray}}
	\newcommand{\eea} {\end{eqnarray}}
\newcommand{\Bea} {\begin{eqnarray*}}
	\newcommand{\Eea} {\end{eqnarray*}}

\numberwithin{equation}{section}

\newtheorem{theorem}{Theorem}[section]
\newtheorem{definition}[theorem]{Definition}
\newtheorem{corollary}[theorem]{Corollary}
\newtheorem{lemma}[theorem]{Lemma}
\newtheorem{prop}[theorem]{Proposition}
\newtheorem{remark}[theorem]{Remark}
\newtheorem{ex}[theorem]{Example}

\begin{document}
	\date{}
	\title{The first Fu\v{c}\'{\i}k curve of the $p$-Laplacian on compact quantum graphs: slopes at the diagonal and nonresonance}
	\author{ {\bf Archana Prajapati$\,^{1,}$\footnote{e-mail: {\tt matharchana3@gmail.com}}}, {\bf Divya Goel$\,^{1,}$\footnote{e-mail: {\tt divya.mat@iitbhu.ac.in}}} \\ $^1\,$Department of Mathematical Sciences, Indian Institute of Technology (BHU),\\ Varanasi 221005, India.}

	\maketitle
	\begin{abstract}
		We study the Fu\v{c}\'{\i}k spectrum of the $p$-Laplacian, $1<p<\infty$, on a compact connected quantum graph, with Dirichlet conditions on a possibly empty set $V_D$ of vertices whose removal does not disconnect the graph and Kirchhoff conditions at the others. We construct the first nontrivial curve $\mathfrak C$ of the spectrum by a mountain pass argument and show that it is Lipschitz and strictly decreasing. As the parameter tends to infinity, $\mathfrak C$ tends to the first eigenvalue if $V_D\neq\emptyset$, and otherwise to the smallest first eigenvalue obtained by a Dirichlet condition at a single point. If the second eigenvalue is simple, we compute the two one-sided slopes of $\mathfrak C$ at the diagonal. The curve is differentiable there if and only if the positive and negative parts of the second eigenfunction have the same $L^p$ norm. On three-edge stars these norms differ for suitable edge lengths and every $p$, so the Lipschitz regularity is optimal. For $p=2$ the simplicity assumption can be dropped, and for a simple eigenvalue the corner is the crossing of the two smooth branches of the spectrum through the diagonal point. Finally, we prove a nonresonance result in which, for $V_D=\emptyset$, the ratio $f(x,s)/|s|^{p-2}s$ may be arbitrarily large at $+\infty$, provided it stays below the threshold at $-\infty$.


		\medskip

	\end{abstract}

	\section{Introduction}
	A quantum graph is a metric graph equipped with a differential operator--typically a Hamiltonian acting on its edges, together with suitable vertex conditions. The field has grown rapidly over three decades, driven by its mathematical structure and by applications in physics, chemistry, engineering and biology. Linus Pauling introduced the subject in quantum chemistry \cite{pauling1936diamagnetic}, while modern developments stem from the work of Kottos and Smilansky on quantum chaos \cite{kottos1997quantum} and from the approximation of thin manifolds by graphs \cite{exner2005convergence}. For introductions, see \cite{berkolaiko2013introduction,kuchment2003quantum}. Spectral estimates for the Laplacian on metric graphs, including bounds on the spectral gap and Rayleigh-type estimates, are obtained in \cite{kurasov2013spectral,kurasov2014rayleigh}.

	A metric graph identifies the edges of a combinatorial graph with intervals of $\mathbb{R}$. Such models arise naturally in quantum wires, thin waveguides and other graph-shaped structures. The Kirchhoff condition is the most common vertex condition: it requires continuity at each vertex and that the outgoing fluxes sum to zero \cite{berkolaiko2013introduction,kuchment2003quantum}; see \cite{kostrykin1999kirchhoff} for the classification of vertex conditions.

	For $1<p<\infty$, the $p$-Laplacian $\Delta_p u=(|u'|^{p-2}u')'$ is a nonlinear generalization of the Laplacian, to which it reduces for $p=2$. It appears in non-Newtonian fluid mechanics, nonlinear elasticity and models of torsional creep. Its spectral theory on intervals and on bounded domains of $\mathbb{R}^N$ is well developed \cite{garcia1987existence,le2006eigenvalue}. On metric graphs, the combinatorial structure and the vertex conditions pose additional difficulties; nonlinear eigenvalue problems for the $p$-Laplacian on graphs are studied in \cite{del2016first,del2016clustering}. On an interval, the eigenfunctions are expressed through the generalized trigonometric functions $\sin_p$, which are studied in detail in \cite{lang2011eigenvalues}.

	The Fu\v{c}\'{\i}k spectrum was introduced by Fu\v{c}\'{\i}k \cite{fuvcik1976boundary} and Dancer \cite{dancer1982boundary} in the 1970s to describe resonance for problems with asymmetric nonlinearities. For the $p$-Laplacian on a bounded domain or on a metric graph, subject to boundary or vertex conditions, the Fu\v{c}\'{\i}k spectrum $\Sigma$ is the set of pairs $(\lambda^+,\lambda^-)\in\mathbb{R}^2$ for which
	\begin{align*}
		-\Delta_p u=\lambda^+(u^+)^{p-1}-\lambda^-(u^-)^{p-1}
	\end{align*}
	has a nontrivial solution, where $u^{\pm}=\max\{\pm u,0\}$. For $\lambda^+=\lambda^-=\lambda$ this is the eigenvalue problem, so $\Sigma$ contains the points $(\lambda,\lambda)$ with $\lambda$ an eigenvalue. The structure of $\Sigma$ governs the solvability of problems with jumping nonlinearities, and it has been studied extensively on Euclidean domains. Cuesta, de Figueiredo and Gossez \cite{cuesta1999beginning} constructed the first nontrivial curve of $\Sigma$ for the $p$-Laplacian on a bounded domain by a mountain pass argument. The construction has been extended to other boundary conditions \cite{motreanu2012fuvcik,dagui2023fucik} and to nonlocal operators \cite{goyal2014fuvcik}, and resonance and nonresonance problems related to $\Sigma$ are studied in \cite{cuesta1999beginning,drabek1999resonance, cuesta1992variational}. For the Neumann problem, the first curve was studied by de Figueiredo and Gossez \cite{defigueiredo1994first} for $p=2$ and by Arias, Campos and Gossez \cite{arias2000antimaximum} for the $p$-Laplacian (see also \cite{arias2008asymmetric}). In one space dimension its horizontal asymptote lies strictly above the first eigenvalue, and for $p=2$ this property is equivalent to a uniform antimaximum principle \cite{defigueiredo1994first}. Near a simple eigenvalue of the Laplacian, the Fu\v{c}\'{\i}k spectrum consists locally of two curves crossing the diagonal \cite{gallouet1981resultats,ruf1981nonlinear}; the minimal and maximal local curves through the diagonal points are described in \cite{perera2015note}, and one-sided tangent lines of Fu\v{c}\'{\i}k curves for discrete systems are studied in \cite{holubova2025fucik}. For the $p$-Laplacian, curves of the spectrum through the diagonal points are constructed in \cite{dancer2001some,micheletti2001fucik,perera2004fuvcik}. The first curve for a $p$-Laplacian with mixed boundary conditions and weights is studied in \cite{leadi2007first}.

	This paper studies the Fu\v{c}\'{\i}k spectrum of the $p$-Laplacian on compact connected quantum graphs. We consider two settings of vertex conditions:
	\begin{enumerate}[label=(\roman*)]
		\item \emph{Dirichlet--Kirchhoff case}: Dirichlet conditions on a nonempty set $V_D\subseteq V(\Gamma)$ and Kirchhoff conditions at the remaining vertices;
		\item \emph{Neumann--Kirchhoff case}: $V_D=\emptyset$, that is, Kirchhoff conditions at all vertices.
	\end{enumerate}
	The two cases differ in their variational structure: in the second, $\lambda_{1,p}=0$ with constant eigenfunctions, and the first nontrivial curve approaches a positive threshold instead of the trivial lines. Sections~\ref{as:4}--\ref{as:6} are written for the first case. Section~\ref{as:8} shows that their proofs apply to the second case (Proposition~\ref{ap:100}) and treats the asymptotics there, while the results of Sections~\ref{as:10} and~\ref{as:nonres} cover both.

	For $1 < p < \infty$, the Dirichlet--Kirchhoff problem is
	\begin{equation}{\label{a1:1}}
		\begin{cases}
			-\bigl(|u'|^{p-2}u'\bigr)'
			= \lambda^+(u^+)^{p-1} - \lambda^-(u^-)^{p-1}
			&\text{on the edges of } \Gamma, \\[2mm]
			u(v) = 0 & \forall\, v \in V_D, \\[2mm]
			\displaystyle\sum_{e \in E_v(\Gamma)} |u_e'(v)|^{p-2}u_e'(v)=0,
			&v\in V(\Gamma)\setminus V_D.
		\end{cases}
	\end{equation}
	Its natural variational setting is
	\begin{align*}
		X(\Gamma,V_D) = \Bigl\{u \in W^{1,p}(\Gamma) : u(v) = 0
		\text{ for all } v \in V_D\Bigr\}.
	\end{align*}
	When $V_D = \emptyset$, we have the Neumann--Kirchhoff problem
	\begin{equation}{\label{a1:2}}
		\begin{cases}
			-\bigl(|u'|^{p-2}u'\bigr)'
			= \lambda^+(u^+)^{p-1} - \lambda^-(u^-)^{p-1}
			&\text{on the edges of } \Gamma, \\[2mm]
			\displaystyle\sum_{e \in E_v(\Gamma)} |u_e'(v)|^{p-2}u_e'(v)=0,
			&\forall\, v \in V(\Gamma),
		\end{cases}
	\end{equation}
	with variational setting $W^{1,p}(\Gamma)$. Throughout the paper we assume that $\Gamma\setminus V_D$ is connected; see \eqref{ah:100} and Remark~\ref{ar:3}.

	The Fu\v{c}\'{\i}k spectra are defined as:
	\begin{align*}
		\Sigma_{p,\Gamma} = \Bigl\{(\lambda^+, \lambda^-) \in \mathbb{R}^2 :
		\text{\eqref{a1:1} admits a nontrivial solution in } X(\Gamma,V_D)\Bigr\},
	\end{align*}
	\begin{align*}
		\mathcal{F}_{p,\Gamma} = \Bigl\{(\lambda^+, \lambda^-) \in \mathbb{R}^2 :
		\text{\eqref{a1:2} admits a nontrivial solution in } W^{1,p}(\Gamma)\Bigr\}.
	\end{align*}
	We write $u = u^+ - u^-$ with $u^\pm = \max\{\pm u, 0\}$ defined edgewise: for
	$u = (u_e)_{e \in E(\Gamma)}$,
	\begin{align*}
		u_e^+(x) = \max\{u_e(x), 0\},
		\quad u_e^-(x) = \max\{-u_e(x), 0\},
		\quad x \in (0,\ell_e),\ e \in E(\Gamma).
	\end{align*}
	When $\lambda^+ = \lambda^- = \lambda$, both problems reduce to the standard eigenvalue problem $-\bigl(|u'|^{p-2}u'\bigr)' = \lambda|u|^{p-2}u$ on the edges. Del Pezzo and Rossi \cite{del2016clustering, del2016first} studied these eigenvalue problems, establishing existence, variational characterization, and qualitative properties.

	For the Dirichlet--Kirchhoff case, the first eigenvalue satisfies
	\begin{equation}{\label{a1:3}}
		\lambda_{1,p}(\Gamma,V_D)
		= \inf\left\{\frac{\displaystyle\int_{\Gamma}|u'(x)|^p\,dx}
		{\displaystyle\int_{\Gamma}|u(x)|^p\,dx}
		: u \in X(\Gamma,V_D),\; u \neq 0\right\}.
	\end{equation}
	Under the hypothesis that $\Gamma\setminus V_D$ is connected, it is simple, with an eigenfunction that is positive on $\Gamma\setminus V_D$ (Lemma~\ref{al:100}). Thus $\Sigma_{p,\Gamma}$ contains the diagonal points $(\lambda,\lambda)$ with $\lambda$ an eigenvalue, the vertical line $\{\lambda_{1,p}(\Gamma,V_D)\} \times \mathbb{R}$, and the horizontal line $\mathbb{R} \times \{\lambda_{1,p}(\Gamma,V_D)\}$. The spectrum is symmetric about the diagonal.

	Under Neumann--Kirchhoff conditions, the first eigenvalue is $\lambda_{1,p}(\Gamma,\emptyset) = 0$
	with the nonzero constants as eigenfunctions. Thus $\mathcal{F}_{p,\Gamma}$ contains
	the diagonal points $(\lambda,\lambda)$ with $\lambda$ a Neumann eigenvalue, the lines $\{0\} \times \mathbb{R}$ and $\mathbb{R} \times \{0\}$,
	and is also symmetric about the diagonal.

	In Euclidean domains, the Fu\v{c}\'{\i}k spectrum has found application to problems with jumping nonlinearities \cite{fuvcik1976boundary, dancer1982boundary} and to suspension bridge dynamics \cite{lazer1990large}. By establishing the Fu\v{c}\'{\i}k spectrum on quantum graphs, we lay the foundation for extending these classical results to network-shaped domains.

To our knowledge, the present paper is the first study of the Fu\v{c}\'{\i}k spectrum on metric graphs. We state the main results.

For $b\ge0$ let $\alpha(b)$ be the mountain pass value defined in \eqref{a:28} below. A point of $\Sigma_{p,\Gamma}$ is called nontrivial if it does not lie on the lines $\{\lambda_{1,p}(\Gamma,V_D)\}\times\mathbb R$ and $\mathbb R\times\{\lambda_{1,p}(\Gamma,V_D)\}$. For $V_D=\emptyset$ we put
\begin{equation}\label{a12:3}
	\bar{\lambda} = \inf \left\{ \int_{\Gamma} |u'|^p\,dx :
	u \in W^{1,p}(\Gamma),\,  \|u\|_p=1,\,
	u \text{ vanishes at some point in } \Gamma \right\}=\min_{x\in\Gamma}\lambda_{1,p}(\Gamma,\{x\}),
\end{equation}
where $\lambda_{1,p}(\Gamma,\{x\})$ is the first eigenvalue with a Dirichlet condition at the single point $x$, regarded as a vertex; $\bar\lambda>0$, and the second equality holds by \eqref{n:lbarmin} (Section~\ref{as:8}). Finally,
\begin{equation}\label{n:Lambda}
\Lambda=\lambda_{1,p}(\Gamma,V_D)\quad\text{if }V_D\neq\emptyset,\qquad \Lambda=\bar\lambda\quad\text{if }V_D=\emptyset .
\end{equation}

\begin{theorem}\label{at:1}
Let $1<p<\infty$ and let $V_D\subseteq V(\Gamma)$, possibly empty, be such that $\Gamma\setminus V_D$ is connected. If $V_D=\emptyset$, read $\mathcal F_{p,\Gamma}$ for $\Sigma_{p,\Gamma}$.
\begin{enumerate}[label=(\roman*)]
\item For every $b\ge0$, the point $(b+\alpha(b),\alpha(b))$ is the first nontrivial point of $\Sigma_{p,\Gamma}$ on the line $(b,0)+t(1,1)$, $t\in\mathbb R$: it belongs to $\Sigma_{p,\Gamma}$, the points $(b+\omega,\omega)$ with $\lambda_{1,p}(\Gamma,V_D)<\omega<\alpha(b)$ do not, and for $\omega\le\lambda_{1,p}(\Gamma,V_D)$ the point $(b+\omega,\omega)$ is not a nontrivial point of $\Sigma_{p,\Gamma}$. In particular, $\alpha(0)=\lambda_{2,p}(\Gamma,V_D)$.
\item The function $\alpha$ is Lipschitz continuous and strictly decreasing on $[0,\infty)$, so that
\begin{align*}
\mathfrak C=\{(b+\alpha(b),\alpha(b)),\ (\alpha(b),b+\alpha(b)):b\ge0\}
\end{align*}
is a Lipschitz curve, symmetric about the diagonal, through $(\lambda_{2,p}(\Gamma,V_D),\lambda_{2,p}(\Gamma,V_D))$.
\item $\alpha(b)>\Lambda$ for every $b\ge0$, and $\alpha(b)\to\Lambda$ as $b\to\infty$. If $V_D=\emptyset$, every nontrivial point of $\mathcal F_{p,\Gamma}$ lies in $(\bar\lambda,\infty)\times(\bar\lambda,\infty)$, and $\bar\lambda$ is the largest number with this property.
\end{enumerate}
\end{theorem}

\begin{figure}[ht]
\begin{center}
\begin{tikzpicture}[scale=0.62]
	\draw[->] (-0.3,0) -- (7.2,0) node[below] {$\lambda^+$};
	\draw[->] (0,-0.3) -- (0,7.2) node[left] {$\lambda^-$};
	\draw[thick] (1,-0.2) -- (1,7) ;
	\draw[thick] (-0.2,1) -- (7,1) ;
	\node[below] at (1,-0.2) {\small$\lambda_{1,p}$};
	\node[left] at (-0.2,1) {\small$\lambda_{1,p}$};
	\draw[dotted] (0,0) -- (7,7);
	\draw[very thick,domain=3:7,samples=60] plot (\x,{1+4/(\x-1)});
	\draw[very thick,domain=3:7,samples=60] plot ({1+4/(\x-1)},\x);
	\fill (3,3) circle (0.08);
	\node[above right] at (3,3) {\small$(\lambda_{2,p},\lambda_{2,p})$};
	\node at (6.3,2.4) {$\mathfrak C$};
	\node at (3.5,-1.4) {$V_D\neq\emptyset$};
	\begin{scope}[xshift=9.5cm]
	\fill[gray!20] (0,0) rectangle (7,1.3);
	\draw[->] (-0.3,0) -- (7.2,0) node[below] {$\lambda^+$};
	\draw[->] (0,-0.3) -- (0,7.2) node[left] {$\lambda^-$};
	\draw[thick] (0,-0.2) -- (0,7);
	\draw[thick] (-0.2,0) -- (7,0);
	\draw[dashed] (1.3,0) -- (1.3,7);
	\draw[dashed] (0,1.3) -- (7,1.3);
	\node[below] at (1.3,-0.2) {\small$\bar\lambda$};
	\node[left] at (-0.2,1.3) {\small$\bar\lambda$};
	\draw[dotted] (0,0) -- (7,7);
	\draw[very thick,domain=3.3:7,samples=60] plot (\x,{1.3+4/(\x-1.3)});
	\draw[very thick,domain=3.3:7,samples=60] plot ({1.3+4/(\x-1.3)},\x);
	\fill (3.3,3.3) circle (0.08);
	\node[above right] at (3.3,3.3) {\small$(\lambda_{2,p},\lambda_{2,p})$};
	\node at (6.3,2.7) {$\mathfrak C$};
	\node at (3.5,-1.4) {$V_D=\emptyset$};
	\end{scope}
\end{tikzpicture}
\end{center}
\caption{Schematic picture of Theorem~\ref{at:1}. Left: for $V_D\neq\emptyset$ the first nontrivial curve $\mathfrak C$ is asymptotic to the trivial lines through $\lambda_{1,p}(\Gamma,V_D)$. Right: for $V_D=\emptyset$ the trivial lines pass through $0$, while $\mathfrak C$ is asymptotic to the lines through $\bar\lambda>0$ (dashed); the shaded half-strip is the region of Theorem~\ref{n:thmstrip}.}
\label{n:figasym}
\end{figure}
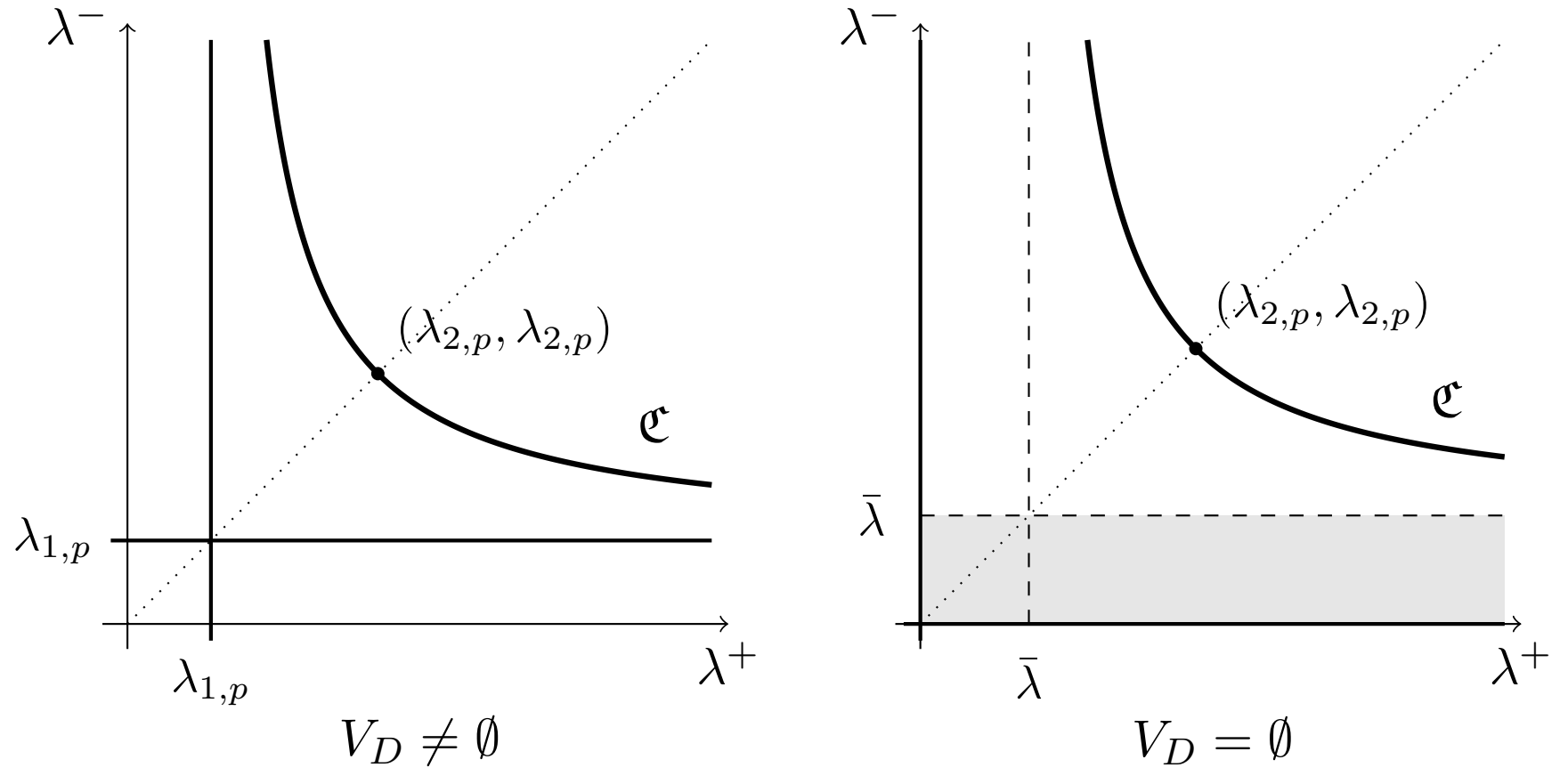

Part \textit{(i)} is proved in Subsection~\ref{as:5}, and its consequence $\alpha(0)=\lambda_{2,p}(\Gamma,V_D)$ is Corollary~\ref{ar:2}. Part \textit{(ii)} is Propositions~\ref{ap:5} and~\ref{ap:6}, and part \textit{(iii)} is Propositions~\ref{ap:7}, \ref{ap:8} and~\ref{at:6}. Thus the first curve is asymptotic to the trivial lines when $V_D\neq\emptyset$, as for the Dirichlet problem on a bounded domain \cite{cuesta1999beginning}, while for $V_D=\emptyset$ it stays at distance $\bar\lambda>0$ from the trivial lines $\{0\}\times\mathbb R$ and $\mathbb R\times\{0\}$; see Figure~\ref{n:figasym}. On an interval with Neumann conditions this threshold is the first eigenvalue of a mixed problem, and the zero of the minimizer lies at an endpoint \cite{defigueiredo1994first,arias2000antimaximum}. On a graph a minimizer in \eqref{a12:3} vanishes at exactly one point (Lemma~\ref{al:12}). Where this point lies depends on the topology of $\Gamma$, and it need not be a point of symmetry (Proposition~\ref{aex:100}).

The second result concerns the regularity of $\mathfrak C$ at the diagonal. By Theorem~\ref{at:1}, $\mathfrak C$ is the graph of a strictly decreasing involution $G$ of $(\Lambda,\infty)$ whose only fixed point is $\lambda_{2,p}=\lambda_{2,p}(\Gamma,V_D)$ (see \eqref{a10:1}); on the side $\lambda^+>\lambda^-$ the curve is $\lambda^-=G(\lambda^+)$, $\lambda^+\ge\lambda_{2,p}$.

\begin{theorem}\label{n:thmdiag}
Let $1<p<\infty$ and let $V_D$ be as in Theorem~\ref{at:1}.
\begin{enumerate}[label=(\roman*)]
\item Assume that every eigenfunction with eigenvalue $\lambda_{2,p}$ is a multiple of one $\phi_{2,p}$ with $\|\phi_{2,p}\|_p=1$, and put $r=\int_\Gamma(\phi_{2,p}^+)^p\,dx\big/\int_\Gamma(\phi_{2,p}^-)^p\,dx$. Then
\begin{align*}
G'(\lambda_{2,p}^+)=-\max\Bigl\{r,\frac1r\Bigr\},\qquad G'(\lambda_{2,p}^-)=-\min\Bigl\{r,\frac1r\Bigr\},
\end{align*}
and $\mathfrak C$ is differentiable at $(\lambda_{2,p},\lambda_{2,p})$ if and only if $r=1$, with slope $-1$.
\item Let $p=2$, let $E_2$ be the eigenspace of $\lambda_{2,2}$, and put $M_0=\max\{\int_\Gamma(u^+)^2\,dx:u\in E_2,\ \|u\|_2=1\}$. Then $\frac12\le M_0<1$,
\begin{align*}
G'(\lambda_{2,2}^+)=-\frac{M_0}{1-M_0},\qquad G'(\lambda_{2,2}^-)=-\frac{1-M_0}{M_0},
\end{align*}
and $\mathfrak C$ is differentiable at $(\lambda_{2,2},\lambda_{2,2})$ if and only if $\int_\Gamma(u^+)^2\,dx=\int_\Gamma(u^-)^2\,dx$ for every $u\in E_2$.
\end{enumerate}
\end{theorem}

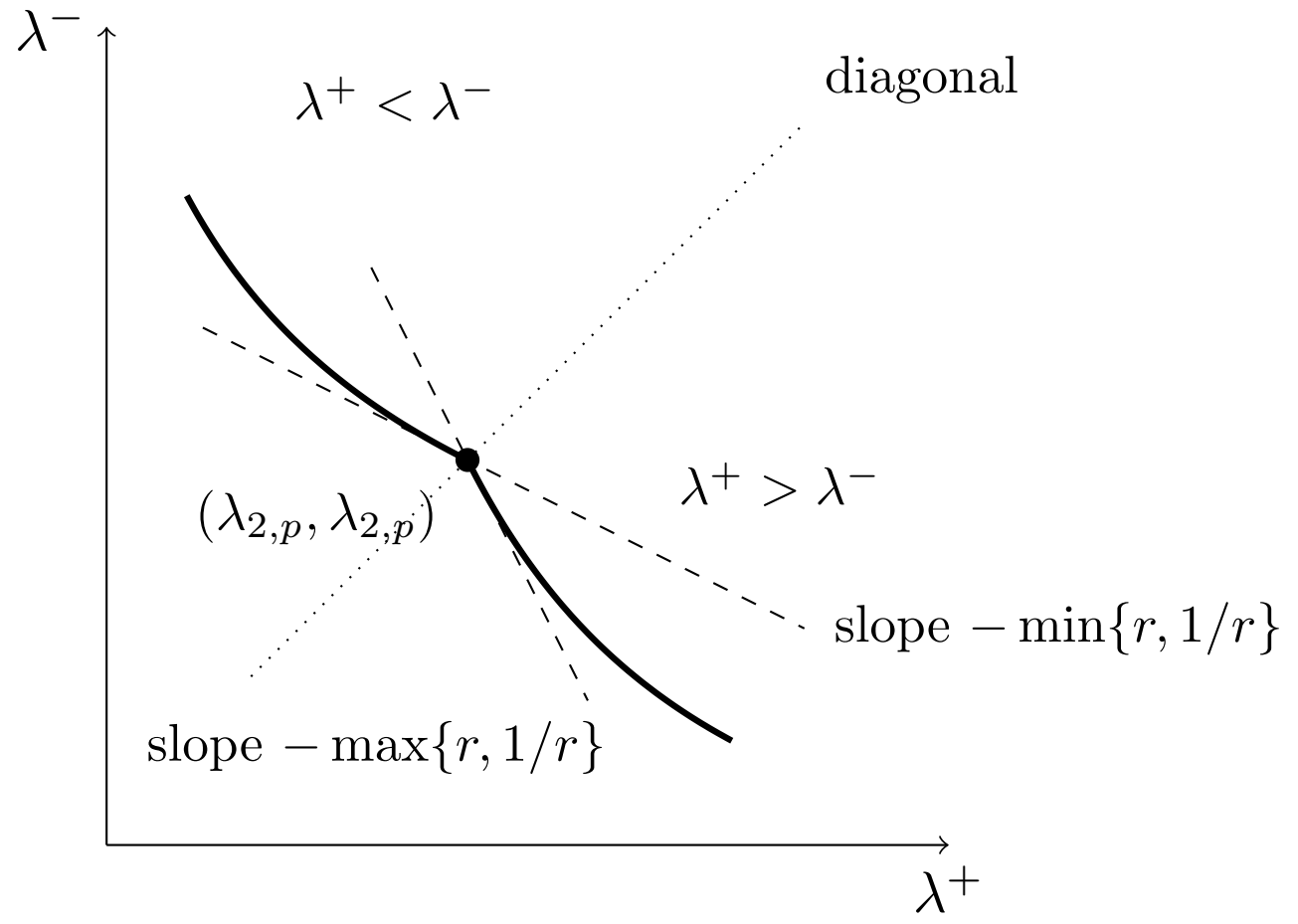
\begin{figure}[ht]
\begin{center}
\begin{tikzpicture}[scale=3.2]
	\draw[->] (0.25,0.2) -- (2.0,0.2) node[below] {$\lambda^+$};
	\draw[->] (0.25,0.2) -- (0.25,1.9) node[left] {$\lambda^-$};
	\draw[dotted] (0.55,0.55) -- (1.7,1.7) node[above right] {\small diagonal};
	\draw[very thick,domain=1:1.55,samples=60] plot (\x,{1/(\x*\x)});
	\draw[very thick,domain=1:1.55,samples=60] plot ({1/(\x*\x)},\x);
	\draw[dashed] (0.8,1.4) -- (1.25,0.5);
	\draw[dashed] (0.45,1.275) -- (1.7,0.65);
	\fill (1,1) circle (0.025);
	\node[below left] at (0.98,0.98) {\small$(\lambda_{2,p},\lambda_{2,p})$};
	\node[anchor=north east] at (1.33,0.5) {\small slope $-\max\{r,1/r\}$};
	\node[anchor=west] at (1.72,0.65) {\small slope $-\min\{r,1/r\}$};
	\node at (1.65,0.95) {\small$\lambda^+>\lambda^-$};
	\node at (0.85,1.75) {\small$\lambda^+<\lambda^-$};
\end{tikzpicture}
\end{center}
\caption{Schematic picture of Theorem~\ref{n:thmdiag} near $(\lambda_{2,p},\lambda_{2,p})$ when $r\neq1$, drawn with $\max\{r,1/r\}=2$, the value of Example~\ref{n:eqstar}. The curve $\mathfrak C$ (solid) is reflected into itself by the diagonal (dotted); its one-sided tangent lines (dashed) have slopes $-2$ and $-1/2$.}
\label{n:figkink}
\end{figure}

Part \textit{(i)} is Theorem~\ref{at10:1} and part \textit{(ii)} is Theorem~\ref{n:p2}. Both are obtained by computing the right derivative of $\alpha$ at $b=0$, and the steeper slope is always on the side $\lambda^+>\lambda^-$; see Figure~\ref{n:figkink}. The corner comes from the symmetry $u\mapsto-u$ of the problem at $b=0$: the critical points at the level $\lambda_{2,p}$ come in pairs $\pm\phi_{2,p}$, whose positive parts may have different weights. On an interval the reflection about the midpoint maps $\phi_{2,p}$ to $-\phi_{2,p}$, so $r=1$ and the classical curve is smooth (Corollary~\ref{ac10:1}). On a graph $r$ can be computed. On three-edge stars with Dirichlet ends, $r<1$ for suitable lengths and every $p$, and $r\neq1$ for almost all lengths when $p=2$ (Proposition~\ref{ap10:1}). Thus the Lipschitz regularity in Theorem~\ref{at:1}\textit{(ii)} cannot be improved in general. On the equilateral star with Dirichlet ends and $p=2$, $\lambda_{2,2}$ is double and Theorem~\ref{n:thmdiag}\textit{(ii)} gives the slopes $-2$ and $-1/2$ (Example~\ref{n:eqstar}), although the graph is symmetric. Away from the diagonal, Proposition~\ref{n:bpos} bounds the difference quotients of $\alpha$ at every $b\ge0$.

For $p=2$ and $\lambda_{2,2}$ simple, Theorem~\ref{n:thmdiag}\textit{(i)} also follows from the local description of the spectrum near a simple eigenvalue \cite{gallouet1981resultats,ruf1981nonlinear} combined with Theorem~\ref{at:1} (Remark~\ref{ar10:2}). Near $(\lambda_{2,2},\lambda_{2,2})$ the nontrivial part of the spectrum is then the union of two smooth curves with slopes $-r$ and $-1/r$, and $\mathfrak C$ is their lower envelope. Thus the corner of $\mathfrak C$ is the crossing of two smooth branches, not a singular point of one branch. This agrees with \cite{holubova2025fucik}, where a single curve leaving a diagonal point can lose smoothness only when the algebraic multiplicity of the eigenvalue exceeds the geometric one. For $p\neq2$ we do not know whether $\Sigma_{p,\Gamma}$ contains a second curve through $(\lambda_{2,p},\lambda_{2,p})$ with slope $-\min\{r,1/r\}$      on the side $\lambda^+>\lambda^-$ (end of Section~\ref{as:10}). The new content of Theorem~\ref{n:thmdiag} is the case $p\neq2$ and, through part \textit{(ii)}, the case of a multiple eigenvalue when $p=2$.

The third result is a nonresonance theorem for
\begin{equation}\label{n:introprob}
-\bigl(|u'|^{p-2}u'\bigr)'=f(x,u)\ \text{ on the edges of }\Gamma,\qquad \sum_{e\in E_v(\Gamma)}|u_e'(v)|^{p-2}u_e'(v)=0\quad\forall v\in V(\Gamma),
\end{equation}
where $f:\Gamma\times\mathbb R\to\mathbb R$ satisfies the $L^\infty$-Carath\'eodory conditions and $F(x,s)=\int_0^sf(x,t)\,dt$.

\begin{theorem}\label{n:thmstrip}
Let $V_D=\emptyset$, and let $\pi_\pm,\Pi_\pm,\theta_\pm,\Theta_\pm\in L^\infty(\Gamma)$ be such that, uniformly in $x\in\Gamma$,
\begin{align*}
\pi_\pm\le\liminf_{s\to\pm\infty}\frac{f(x,s)}{|s|^{p-2}s}\le\limsup_{s\to\pm\infty}\frac{f(x,s)}{|s|^{p-2}s}\le\Pi_\pm,\qquad
\theta_\pm\le\liminf_{s\to\pm\infty}\frac{pF(x,s)}{|s|^{p}}\le\limsup_{s\to\pm\infty}\frac{pF(x,s)}{|s|^{p}}\le\Theta_\pm,
\end{align*}
with $0\le\pi_+\le\Pi_+$, $0\le\pi_-\le\Pi_-\le\bar\lambda$, $0\le\theta_+\le\Theta_+$, $0\le\theta_-\le\Theta_-\le\bar\lambda$ a.e.\ in $\Gamma$, and $\theta_\pm>0$ on sets of positive measure. Then \eqref{n:introprob} has a solution in $W^{1,p}(\Gamma)$.
\end{theorem}

No upper bound is imposed on $\Pi_+$ and $\Theta_+$: the asymptotic quotient of $f$ may be as large as we like at $+\infty$, provided it stays below $\bar\lambda$ at $-\infty$. Theorem~\ref{n:thmstrip} is Corollary~\ref{ac:strip}. It follows from a nonresonance theorem between $(\lambda_{1,p}(\Gamma,V_D),\lambda_{1,p}(\Gamma,V_D))$ and $\mathfrak C$, valid for every $V_D$ (Theorem~\ref{at:5}), together with Theorem~\ref{at:1}\textit{(iii)}. On an interval the same half-strip follows from the one-dimensional Neumann results \cite{defigueiredo1994first,arias2000antimaximum}; on a graph its width is $\bar\lambda$.

	The proofs combine variational methods---Ekeland's principle and the Mountain Pass Theorem for a $C^1$ functional on a $C^1$ manifold \cite{cuesta2003minimax}---with the Palais--Smale compactness condition and the $(S_+)$ property of the $p$-Laplacian (Lemma~\ref{n:splus}). Central to the approach are the compact embeddings
	$W^{1,p}(\Gamma) \hookrightarrow L^p(\Gamma)$ and $W^{1,p}(\Gamma) \hookrightarrow C(\Gamma)$
	on compact quantum graphs (Theorems \ref{at:3}, \ref{at2:3}), which provide the
compactness required for variational arguments.

	The construction of $\mathfrak C$ and the proofs of its basic properties follow the scheme of Cuesta, de Figueiredo and Gossez \cite{cuesta1999beginning} for bounded domains and, for $V_D=\emptyset$, the one-dimensional Neumann results of \cite{defigueiredo1994first,arias2000antimaximum}. The graph enters at the following points. First, the positivity of the first eigenfunction, the simplicity of $\lambda_{1,p}$ and the fact that only $\lambda_{1,p}$ has eigenfunctions of constant sign (Lemma~\ref{al:100}) rest on the hypothesis that $\Gamma\setminus V_D$ is connected; this hypothesis replaces the connectedness of a domain, is weaker than the condition $V_D\subseteq T(\Gamma)$ of \cite{del2016first}, and cannot be dropped (Remark~\ref{ar:3}). Second, the enlargement of nodal domains in Lemma~\ref{al:9} needs a boundary point outside $V_D$ whose exterior has positive length. Third, the convergence $\alpha(b)\to\lambda_{1,p}(\Gamma,V_D)$ in Proposition~\ref{ap:7} uses an explicit test function at a Dirichlet vertex. Fourth, for $V_D=\emptyset$ the limit of $\alpha(b)$ is $\bar\lambda=\min_{x\in\Gamma}\lambda_{1,p}(\Gamma,\{x\})$, and the point where the minimum is attained depends on the topology of $\Gamma$ (Proposition~\ref{aex:100}). Fifth, on a graph the weights of the second eigenfunction can be computed, which makes Theorem~\ref{n:thmdiag} effective (Proposition~\ref{ap10:1}, Example~\ref{n:eqstar}). Proposition~\ref{ap:100} lists the places where $V_D\neq\emptyset$ is used. The remaining arguments of Sections~\ref{as:4}--\ref{as:6} and~\ref{as:nonres} use only the variational structure and the compact embeddings of Theorems~\ref{at:3} and~\ref{at2:3}; the routine ones are given in condensed form.

	Paper organization: Section \ref{as:2} recalls background on quantum graphs and
	embedding theorems. Section \ref{as:4} constructs the first nontrivial curve via the Mountain Pass
	Theorem, proves the isolatedness of the trivial lines and establishes
	Theorem \ref{at:1}, yielding a variational characterization of $\lambda_{2,p}(\Gamma, V_D)$.
	Section \ref{as:6} analyzes the first nontrivial curve: Lipschitz continuity,
	monotonicity, and asymptotics. Section \ref{as:8} addresses the Neumann--Kirchhoff setting and $\mathcal{F}_{p,\Gamma}$: it shows that the results of Sections \ref{as:4}--\ref{as:6} hold for $V_D=\emptyset$ and determines the asymptotics of $\mathfrak C$ there. Section \ref{as:10} studies the differentiability of the first nontrivial curve at the diagonal and treats a three-edge star, and Section \ref{as:nonres} proves the nonresonance theorem together with the half-strip corollary.

	\section{Preliminaries}{\label{as:2}}
	In this section, we present a brief overview of fundamental notions concerning quantum graphs; see \cite{berkolaiko2013introduction} and the references therein for further details.

	A graph $\Gamma$ consists of a finite set of vertices $V(\Gamma) = \{v_i\}$ and a set of edges $E(\Gamma) = \{e_i\}$ connecting the vertices, and is denoted by $(V, E)$. An edge and a vertex on that edge are called \emph{incident}. We write $v \in e$ when $e$ and $v$ are incident. The set of edges incident to $v$ is denoted by $E_{v}(\Gamma)$. A graph $\Gamma$ is \emph{connected} if any two vertices of $\Gamma$ are joined by a path, or equivalently if $\Gamma$ is connected as a topological space. A finite metric graph whose edges all have finite lengths is called \emph{compact}.

	\begin{definition}
		A graph $\Gamma$ is a metric graph if
		\begin{enumerate}
			\item each edge $e$ is assigned a positive length
			$\ell_e \in (0,\infty)$;
			\item each edge is assigned a coordinate $x_e \in I_e = [0,\ell_e]$, increasing in the direction of the edge.
		\end{enumerate}
	\end{definition}

	For two vertices $v_i$ and $v_j$, the distance $d(v_i, v_j)$ is the minimal length of a path connecting them. More generally, a compact metric graph $\Gamma$ is a metric measure space with the distance $d(x,y)$ between two points $x$ and $y$ of the graph (not necessarily vertices) defined as the length of the shortest path on $\Gamma$ connecting them, that is,
	\begin{align*}
		d(x,y) = \inf\left\{
		\int_0^1|\gamma'(t)|\,dt
		: \gamma : [0,1] \to \Gamma \text{ is Lipschitz},\;
		\gamma(0) = x,\; \gamma(1) = y
		\right\}.
	\end{align*}
	For $1 \le p \le \infty$, a function $u$ belongs to $L^p(\Gamma)$ if $u_e \in L^p(0,\ell_e)$ for all $e \in E(\Gamma)$, and
	\begin{align*}
		\|u\|_{L^p(\Gamma)}^p
		= \sum_{e \in E(\Gamma)}\|u_e\|_{L^p(0,\ell_e)}^p < \infty.
	\end{align*}
	With this definition, the direct sum decomposition
	$L^p(\Gamma) = \bigoplus_{e \in E(\Gamma)} L^p(e)$ holds. The Sobolev space $W^{1,p}(\Gamma)$ is the space of continuous functions $u$ on $\Gamma$ such that $u_e \in W^{1,p}(I_e)$ for all $e \in E(\Gamma)$, endowed with the norm
	\begin{align*}
		\|u\|_{W^{1,p}(\Gamma)}^p
		= \sum_{e \in E(\Gamma)}
		\Bigl(\|u_e\|_{L^p(0,\ell_e)}^p + \|u_e'\|_{L^p(0,\ell_e)}^p\Bigr)
		< \infty.
	\end{align*}
	The continuity condition in the definition of $W^{1,p}(\Gamma)$ implies that for each vertex $v \in V(\Gamma)$, the function $u$ takes the same value at $v$ on all edges $e \in E_v(\Gamma)$. The following results are proved in \cite{del2016first}.
	\begin{theorem}
		Let $\Gamma$ be a compact connected quantum graph and $1 \le p \le \infty$. Then the Sobolev space $W^{1,p}(\Gamma)$ is a Banach space. Moreover, $W^{1,p}(\Gamma)$ is reflexive for $1 < p < \infty$ and separable for $1 \le p < \infty$.
	\end{theorem}
	\begin{theorem}{\label{at:3}}
		Let $\Gamma$ be a compact connected quantum graph and $1 < p < \infty$. Then the embedding $W^{1,p}(\Gamma) \hookrightarrow L^{q}(\Gamma)$ is compact for all $1 \le q \le \infty$.
	\end{theorem}
	In particular, the embedding $W^{1,p}(\Gamma) \hookrightarrow L^{q}(\Gamma)$ is continuous, and there exists a constant $C > 0$ such that
	\begin{align*}
		\|u\|_{L^{q}(\Gamma)} \le C\|u\|_{W^{1,p}(\Gamma)},
		\quad \text{for all } u \in W^{1,p}(\Gamma).
	\end{align*}
	We also have the following result.
	\begin{theorem}{\label{at2:3}}
		Let $\Gamma$ be a compact connected quantum graph and $1 < p < \infty$. Then the embedding $W^{1,p}(\Gamma) \hookrightarrow C(\Gamma)$ is compact.
	\end{theorem}

	We shall use the $(S_+)$ property of the $p$-Laplacian in the following form.
	\begin{lemma}\label{n:splus}
		Let $X$ be a closed subspace of $W^{1,p}(\Gamma)$, $1<p<\infty$, and let $u_n\rightharpoonup u$ in $X$ with
		\begin{align*}
			\limsup_{n\to\infty}\int_\Gamma|u_n'|^{p-2}u_n'(u_n'-u')\,dx\le0 .
		\end{align*}
		Then $u_n\to u$ in $X$.
	\end{lemma}
	\begin{proof}
		The weak convergence $u_n'\rightharpoonup u'$ in $L^p(\Gamma)$ implies $\int_\Gamma|u'|^{p-2}u'(u_n'-u')\,dx\to0$, and hence
		\begin{align*}
			\limsup_{n\to\infty}\int_\Gamma\bigl(|u_n'|^{p-2}u_n'-|u'|^{p-2}u'\bigr)(u_n'-u')\,dx\le0 .
		\end{align*}
		By H\"older's inequality, the integral on the left is at least $\bigl(\|u_n'\|_p^{p-1}-\|u'\|_p^{p-1}\bigr)\bigl(\|u_n'\|_p-\|u'\|_p\bigr)\ge0$, so that $\|u_n'\|_p\to\|u'\|_p$. By the uniform convexity of $L^p(\Gamma)$, weak convergence and convergence of the norms imply $u_n'\to u'$ in $L^p(\Gamma)$. By Theorem~\ref{at:3}, $u_n\to u$ in $L^p(\Gamma)$, and therefore $u_n\to u$ in $X$. \qed
	\end{proof}

	\begin{definition}
		A quantum graph is a metric graph $\Gamma$ equipped with a
		differential operator $H$ acting on functions defined on the edges of $\Gamma$, together with prescribed vertex conditions. In this work, we consider the operator $H$ given by
		\begin{align*}
			H(u)(x) = -\Delta_p u(x) = -\bigl(|u'(x)|^{p-2}u'(x)\bigr)',
		\end{align*}
		where the derivative is taken along each edge of $\Gamma$.
	\end{definition}
	Let $V_D$ be a (possibly empty) subset of $V(\Gamma)$. We consider the following vertex conditions:
	\begin{align*}
		\begin{cases}
			u \text{ is continuous on } \Gamma, \\[0.3em]
			u(v) = 0 \quad \text{for all } v \in V_D, \\[0.3em]
			\displaystyle\sum_{e \in E_v(\Gamma)} |u_e'(v)|^{p-2}u_e'(v)=0,  \quad \text{for all } v \in V(\Gamma)\setminus V_D,
		\end{cases}
	\end{align*}
	where the derivatives are taken in the outward direction from the vertex. When $V_D\neq\emptyset$, these are called the Dirichlet-Kirchhoff vertex conditions. When $V_D=\emptyset$, they reduce to
	\begin{align*}
		\begin{cases}
			u \text{ is continuous on } \Gamma, \\[0.3em]
			\displaystyle\sum_{e \in E_v(\Gamma)} |u_e'(v)|^{p-2}u_e'(v)=0, \quad \text{for all } v \in V(\Gamma),
		\end{cases}
	\end{align*}
	which are called the Neumann-Kirchhoff vertex conditions. Throughout this work, the integral over $\Gamma$ is defined by
	\begin{align*}
		\int_{\Gamma} u(x)\,dx  = \sum_{e \in E(\Gamma)}\int_{0}^{\ell_e} u_e(x)\,dx.
	\end{align*}

	Let $\Gamma$ be a compact connected quantum graph and let $V_D$ be a nonempty subset of $V(\Gamma)$. The Dirichlet-Kirchhoff eigenvalue problem reads (for $V_D=\emptyset$ the same definitions are used with $X(\Gamma,\emptyset)=W^{1,p}(\Gamma)$):
	\begin{equation}{\label{a:6}}
		\begin{cases}
			-\bigl(|u'|^{p-2}u'\bigr)' = \lambda|u|^{p-2}u
			&\text{on the edges of } \Gamma, \\[6pt]
			u(v) = 0 & \forall\, v \in V_D, \\[8pt]
			\displaystyle\sum_{e \in E_v(\Gamma)} |u_e'(v)|^{p-2}u_e'(v)=0, &v\in V(\Gamma)\setminus V_D.
		\end{cases}
	\end{equation}
	The weak formulation of \eqref{a:6} is
	\begin{equation}\label{a:7}
		\begin{cases}
			\displaystyle\int_{\Gamma}|u'|^{p-2}u'\phi'\,dx
			= \lambda\int_{\Gamma}|u|^{p-2}u\phi\,dx,
			& \forall\,\phi \in X(\Gamma,V_D), \\[6pt]
			u \in X(\Gamma,V_D).
		\end{cases}
	\end{equation}
	A number $\lambda \in \mathbb{R}$ is an eigenvalue of the $p$-Laplacian if \eqref{a:6}, or equivalently its weak formulation \eqref{a:7}, has a nontrivial solution $u \in X(\Gamma, V_D)$. Any such solution is called an eigenfunction corresponding to $\lambda$.

	The first eigenvalue of the $p$-Laplacian is given by \eqref{a1:3}, and an eigenfunction associated with $\lambda_{1,p}(\Gamma, V_D)$ exists by a standard compactness argument. When $V_D \neq \emptyset$, the norm in $W^{1,p}(\Gamma)$ is equivalent to the seminorm
	\begin{align*}
		\left(\int_{\Gamma}|u'(x)|^{p}\,dx\right)^{1/p}
		= \left(\sum_{e \in E(\Gamma)}\|u'_e\|^{p}_{L^{p}(0,\ell_e)}\right)^{1/p}.
	\end{align*}
	We set $X(\Gamma,\emptyset)=W^{1,p}(\Gamma)$ and, for every $V_D$, endow $X(\Gamma,V_D)$ with the norm of $W^{1,p}(\Gamma)$, which we denote by $\|\cdot\|$. In all arguments below, boundedness in $X(\Gamma,V_D)$ is deduced from a bound on $\int_\Gamma|u'|^p\,dx$ together with $\|u\|_p=1$, so the equivalence above is not needed and the arguments also apply when $V_D=\emptyset$.

By \cite{del2016first}, for $V_D\neq\emptyset$ the infimum \eqref{a1:3} is positive and is attained by a nonnegative eigenfunction. Throughout the paper we assume
	\begin{equation}\label{ah:100}
		\Gamma\setminus V_D\ \text{ is connected.}\tag{H}
	\end{equation}
	Hypothesis \eqref{ah:100} holds when $V_D=\emptyset$ and, more generally, when $V_D\subseteq T(\Gamma)$, the set of vertices of degree one, because removing such vertices does not disconnect $\Gamma$. This is the condition under which simplicity of the first eigenvalue is proved in \cite{del2016first}. Hypothesis \eqref{ah:100} also allows, for instance, a cycle with one Dirichlet vertex. Under \eqref{ah:100} the first eigenvalue is simple (Lemma~\ref{al:100}). Without it, simplicity may fail, as the next remark shows.
	\begin{remark}{\label{ar:3}}
		Without \eqref{ah:100}, the first eigenvalue need not be simple. To illustrate this, let $\Gamma$ be a star graph with three edges and four vertices, that is, $V(\Gamma) = \{v_0, v_1, v_2, v_3\}$, where $v_0$ is the central vertex and $E(\Gamma) = \{[v_0,v_1],\,[v_0,v_2],\,[v_0,v_3]\}$, with each edge having the same length $L$. Let $V_D = V(\Gamma)$, that is, the Dirichlet conditions are imposed at all vertices.
		\begin{figure}[ht]
			\begin{center}
				\begin{tikzpicture}[scale=1.1]

					\draw[thick] (0,0) circle (0.08);
					\node[left] at (-0.05,0.05) {$v_0$};

					\draw[line width=2pt] (0,0) -- (0,1.8);
					\draw[thick] (0,1.8) circle (0.08);
					\node[left] at (0,1.8) {$v_2$};

					\draw[line width=2pt] (0,0) -- (-1.6,-1.6);
					\draw[thick] (-1.6,-1.6) circle (0.08);
					\node[left] at (-1.6,-1.6) {$v_1$};

					\draw[line width=2pt] (0,0) -- (1.8,-0.7);
					\draw[thick] (1.8,-0.7) circle (0.08);
					\node[right] at (1.8,-0.7) {$v_3$};

					\node at (3.5,0) {$=$};

					\draw[line width=2pt] (6,0) -- (6,1.8);
					\draw[thick] (6,0) circle (0.08);
					\node[left] at (6,-0.1) {$v_0$};

					\draw[thick] (6,1.8) circle (0.08);
					\node[left] at (6,1.8) {$v_2$};

					\draw[line width=2pt] (5,0) -- (3.5,-1.6);

					\draw[thick] (5,0) circle (0.08);
					\node[left] at (5,0) {$v_0$};

					\draw[thick] (3.5,-1.6) circle (0.08);
					\node[left] at (3.5,-1.6) {$v_1$};

					\draw[line width=2pt] (7,0) -- (8.6,-0.7);

					\draw[thick] (7,0) circle (0.08);
					\node[right] at (7.05,0) {$v_0$};

					\draw[thick] (8.6,-0.7) circle (0.08);
					\node[right] at (8.6,-0.7) {$v_3$};

				\end{tikzpicture}
			\end{center}
			\caption{A star graph with a Dirichlet condition at the central vertex.}
			\label{afig:1}
		\end{figure}
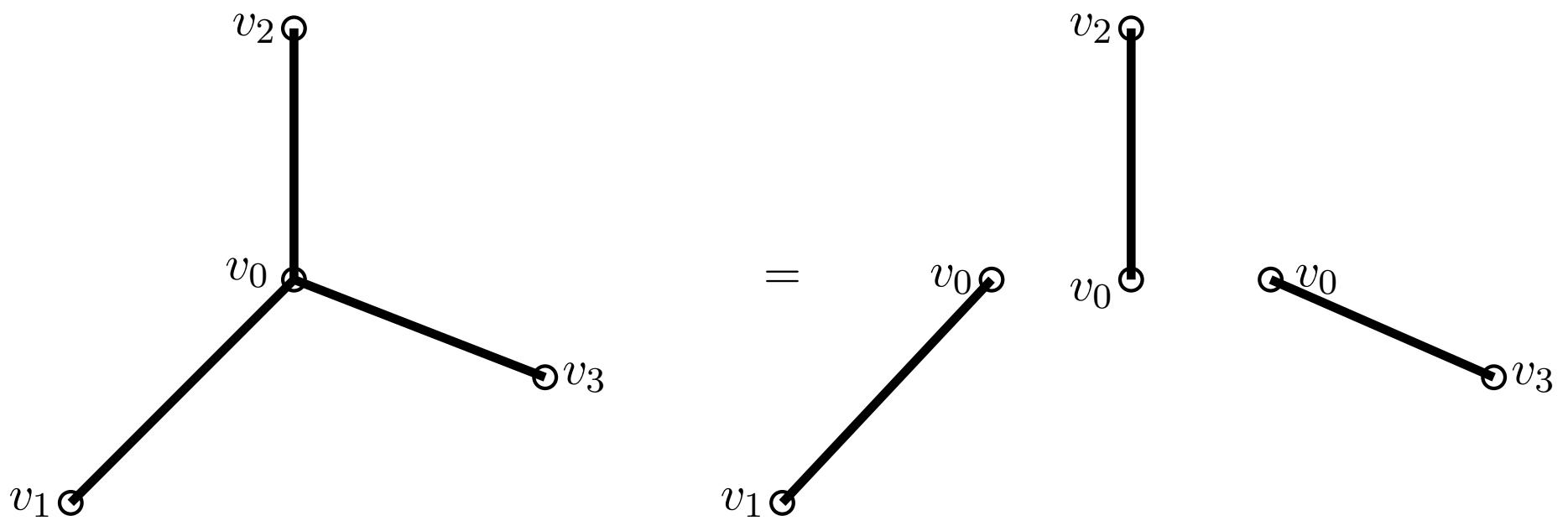

		In this setting, each edge becomes independent, and the first eigenvalue of the $p$-Laplacian is
		\begin{align*}
			\lambda_{1,p}(\Gamma,V_D)
			= \left(\frac{\pi_p}{L}\right)^p\frac{p}{p'},
			\quad \text{where } \frac{1}{p} + \frac{1}{p'} = 1,
		\end{align*}
		(see \cite{lang2011eigenvalues}). The corresponding eigenfunctions can be chosen as
		\begin{align*}
			u_i(x)=\frac{L}{\pi_p}\sin_p\Bigl(\frac{\pi_p}{L}x\Bigr)\ \text{ on }[v_0,v_i]=[0,L],\qquad u_i=0\ \text{ elsewhere},\qquad i=1,2,3,
		\end{align*}
		and $u_1, u_2, u_3$ are linearly independent eigenfunctions associated with $\lambda_{1,p}(\Gamma,V_D)$, so the first eigenvalue is not simple. This failure of simplicity arises because the Dirichlet condition at the central vertex $v_0$ decouples the graph into three independent intervals (Figure~\ref{afig:1}), so \eqref{ah:100} fails.
            
	\end{remark}

		\begin{lemma}\label{al:100}
			Let $\Gamma$ be a compact connected quantum graph, $1<p<\infty$, and let $V_D\subseteq V(\Gamma)$, possibly empty, satisfy \eqref{ah:100}. Write $\lambda_1=\lambda_{1,p}(\Gamma,V_D)$.
			\begin{enumerate}[label=(\roman*)]
				\item Let $m\in L^\infty(\Gamma)$ and let $u\in X(\Gamma,V_D)$, $u\ge0$, $u\not\equiv0$, satisfy
			\begin{equation}\label{a:101}
					\int_\Gamma|u'|^{p-2}u'\phi'\,dx=\int_\Gamma m(x)\,u^{p-1}\phi\,dx\qquad\forall\,\phi\in X(\Gamma,V_D).
				\end{equation}
				Then $u>0$ on $\Gamma\setminus V_D$.
				\item $\lambda_1$ is simple: it has an eigenfunction $\phi_{1,p}$ with $\|\phi_{1,p}\|_p=1$ and $\phi_{1,p}>0$ on $\Gamma\setminus V_D$, and every eigenfunction associated with $\lambda_1$ is a multiple of $\phi_{1,p}$.
				\item If $u$ is as in (i), then
				\begin{equation}\label{a:103}
					\int_\Gamma\bigl(m(x)-\lambda_1\bigr)\phi_{1,p}^p\,dx\le0 .
				\end{equation}
				In particular, an eigenfunction that does not change sign is associated with $\lambda_1$.
			\end{enumerate}
		\end{lemma}

		\begin{proof}
\textit{(i)}
On each edge $e$, testing \eqref{a:101} with
$\phi\in C_c^\infty(0,\ell_e)$ shows that
$w_e=|u_e'|^{p-2}u_e'$ satisfies
\[
w_e'=-m\,u_e^{p-1}\in L^\infty(0,\ell_e).
\]
Then $w_e$ is Lipschitz, $u_e\in C^1([0,\ell_e])$,
and an integration by parts in \eqref{a:101} with
$\phi$ supported near a vertex $v\notin V_D$
gives the Kirchhoff condition
$\sum_{e\in E_v(\Gamma)}
|u_e'(v)|^{p-2}u_e'(v)=0$,
the derivatives being taken in the outward direction
from $v$. On each edge,
\[
-(|u_e'|^{p-2}u_e')'
+\|m\|_\infty u_e^{p-1}\ge0 \quad\text{ and } \quad u_e\ge0.
\]
The strong maximum principle and
the boundary point lemma of V\'azquez
\cite[Theorem~5]{vazquez1984strong}, applied with
$\beta(s)=\|m\|_\infty s^{p-1}$, for which
\[
\int_0^1(\beta(s)s)^{-1/p}\,ds=\infty,
\]
show that either $u_e\equiv0$, or $u_e>0$ in
$(0,\ell_e)$ and the outward derivative of $u_e$
is positive at each endpoint where $u_e$ vanishes. Let
$Z=\{x\in\Gamma\setminus V_D:u(x)=0\}$.
The set $Z$ is closed in $\Gamma\setminus V_D$.
It is also open: a point of $Z$ inside an edge
lies on an edge where $u\equiv0$. If a vertex $v$
belongs to $Z$, then $v\notin V_D$, and every term
of the Kirchhoff sum at $v$ is nonnegative.
The terms corresponding to edges with
$u_e\not\equiv0$ are positive, which forces
$u_e\equiv0$ on every edge incident to $v$.
By \eqref{ah:100}, either $Z=\emptyset$ or
$Z=\Gamma\setminus V_D$. The latter implies
$u\equiv0$ because $V_D$ is finite.
Hence $u>0$ on $\Gamma\setminus V_D$.

\textit{(ii)} A minimizer of \eqref{a1:3} exists by Theorem~\ref{at:3}. If $u$ is a minimizer, then so is $|u|$. By the Lagrange multiplier rule, $|u|$ satisfies \eqref{a:101} with $m\equiv\lambda_1$.
By part \textit{(i)}, $|u|>0$ on $\Gamma\setminus V_D$. Since $u$ is continuous and $\Gamma\setminus V_D$ is connected, $u$ does not change sign. Every eigenfunction associated with $\lambda_1$ is a minimizer (test the equation with the eigenfunction itself), so it suffices to show that two positive minimizers $u,v$ with $\|u\|_p=\|v\|_p=1$ coincide. We use the convexity of $\sigma\mapsto\int_\Gamma|(\sigma^{1/p})'|^p\,dx$ on nonnegative functions, observed by Belloni and Kawohl \cite{belloni2002direct}. For $t\in(0,1)$ put $w_t=\bigl((1-t)u^p+tv^p\bigr)^{1/p}$. Inside each edge $u,v>0$ and, with $\theta=(1-t)u^p/w_t^p\in(0,1)$,
\begin{align*}
\frac{w_t'}{w_t}=\theta\,\frac{u'}{u}+(1-\theta)\,\frac{v'}{v},
\end{align*}
so the convexity of $s\mapsto|s|^p$ gives
\begin{equation}\label{a:102}
				|w_t'|^p\le(1-t)|u'|^p+t|v'|^p\qquad\text{a.e. on }\Gamma,
			\end{equation}
			with equality only where $u'/u=v'/v$. Thus $w_t\in X(\Gamma,V_D)$, $\|w_t\|_p=1$ and $\int_\Gamma|w_t'|^p\,dx\le\lambda_1$, so equality holds a.e.\ in \eqref{a:102} and $(u/v)'=0$ inside each edge. The quotient $u/v$ is therefore constant on each edge, and these constants agree at the vertices of $\Gamma\setminus V_D$, where $u,v>0$. By \eqref{ah:100}, $u=cv$, and $c=1$.

			\textit{(iii)} For $\varepsilon>0$ the function $\psi_\varepsilon=\phi_{1,p}^p/(u+\varepsilon)^{p-1}$ belongs to $X(\Gamma,V_D)$, since $\phi_{1,p},u\in W^{1,p}(\Gamma)\cap C(\Gamma)$, $u+\varepsilon\ge\varepsilon$ and $\phi_{1,p}=0$ on $V_D$. By Picone's identity \cite{allegretto}, $|\phi_{1,p}'|^p\ge|u'|^{p-2}u'\psi_\varepsilon'$ a.e.\ on $\Gamma$. Testing \eqref{a:101} with $\psi_\varepsilon$, we obtain
			\begin{align*}
				\lambda_1\int_\Gamma\phi_{1,p}^p\,dx=\int_\Gamma|\phi_{1,p}'|^p\,dx\ge\int_\Gamma m(x)\Bigl(\frac{u}{u+\varepsilon}\Bigr)^{p-1}\phi_{1,p}^p\,dx .
			\end{align*}
			By \textit{(i)}, $u>0$ a.e., so \eqref{a:103} follows by dominated convergence as $\varepsilon\to0$. If $u\ge0$ is an eigenfunction with eigenvalue $\lambda$, then \eqref{a:103} with $m\equiv\lambda$ gives $\lambda\le\lambda_1$, while $\lambda\ge\lambda_1$ by \eqref{a1:3}; for $u\le0$ we apply this to $-u$. \qed
		\end{proof}

		Throughout the paper $\phi_{1,p}$ denotes the eigenfunction of Lemma~\ref{al:100}\textit{(ii)}; for $V_D=\emptyset$ it is the constant $\ell(\Gamma)^{-1/p}$. We denote by $\lambda_{2,p}(\Gamma,V_D)$ the least eigenvalue larger than $\lambda_{1,p}(\Gamma,V_D)$,
		\begin{equation}\label{a:104}
			\lambda_{2,p}(\Gamma,V_D)=\inf\bigl\{\lambda:\lambda \text{ is an eigenvalue and }\lambda>\lambda_{1,p}(\Gamma,V_D)\bigr\}.
		\end{equation}
		Corollary~\ref{ar:2} shows that this infimum is attained and equals the minimax value $\alpha(0)$ of \eqref{a:28}.

	We study nonlinear eigenvalue problems associated with
	the $p$-Laplacian on the compact connected quantum graph $\Gamma$, namely problems \eqref{a1:1} and \eqref{a1:2}.

	\begin{definition}
		A function $u \in X(\Gamma, V_D)$ (respectively, $u \in W^{1,p}(\Gamma)$) is a weak solution of problem \eqref{a1:1} (respectively, \eqref{a1:2}) if
		\begin{equation}\label{a2:7}
			\int_{\Gamma}|u'|^{p-2}u'\phi'\,dx
			= \int_{\Gamma}\bigl(\lambda^+(u^{+})^{p-1}
			- \lambda^-(u^{-})^{p-1}\bigr)\phi\,dx,
			\qquad \forall\,\phi \in X(\Gamma, V_D)
		\end{equation}
		(\text{respectively, } $\phi\in W^{1,p}(\Gamma)$).
	\end{definition}

	The energy functional associated with problems \eqref{a1:1} and \eqref{a1:2} is
	\begin{align*}
		I(u) = \frac{1}{p}\int_{\Gamma}|u'|^{p}\,dx
		- \frac{\lambda^+}{p}\int_{\Gamma}(u^{+})^{p}\,dx
		- \frac{\lambda^-}{p}\int_{\Gamma}(u^{-})^{p}\,dx.
	\end{align*}
	The functional $I$ is of class $C^1$, and its critical points are weak solutions of problems \eqref{a1:1} and \eqref{a1:2}.

	\section{Construction of the first nontrivial curve}{\label{as:4}}

	In this section we construct the curve $\mathfrak C$ by a minimax argument and show that it is the first nontrivial curve of $\Sigma_{p,\Gamma}$. The points of $\Sigma_{p,\Gamma}$ are characterized as critical values of a suitably restricted functional. As a consequence, we obtain a variational characterization of $\lambda_{2,p}(\Gamma,V_D)$.

	\subsection{The minimax value}

	Fix $b \in \mathbb{R}$ and consider the $C^{1}$ functional $I_{b} : X(\Gamma, V_D) \to \mathbb{R}$ given by
	\begin{align*}
		I_b(u) = \int_{\Gamma} |u'|^{p}\,dx - b \int_{\Gamma} (u^{+})^{p}\,dx,
	\end{align*}
	restricted to the manifold
	\begin{align*}
		M = \left\{u \in X(\Gamma,V_D) : J(u) = \int_{\Gamma} |u|^{p}\,dx = 1\right\}.
	\end{align*}
	Let $\tilde{I}_b$ denote the restriction of $I_b$ to $M$. As $I_b$ is a $C^1$ functional, the Lagrange multiplier rule shows that $u \in M$ is a critical point of $\tilde{I}_b$ if and only if there exists $c \in \mathbb{R}$ with $I'_b(u) = cJ'(u)$, that is,
	\begin{equation}\label{a:14}
		\int_{\Gamma} |u'|^{p-2}u'\phi'\,dx - b\int_{\Gamma} (u^{+})^{p-1}\phi\,dx
		= c\int_{\Gamma} |u|^{p-2}u\phi\,dx,
		\qquad \forall\, \phi \in X(\Gamma,V_D).
	\end{equation}
	Equivalently, $u \in M$ is a nontrivial weak solution of \eqref{a1:1} with $(\lambda^+,\lambda^-)=(b+c,c)$, so that $(b+c, c) \in \Sigma_{p,\Gamma}$. Taking $\phi = u$ in~\eqref{a:14} shows that the Lagrange multiplier $c$ equals the critical value $\tilde{I}_b(u)$. This proves the following lemma.
	\begin{lemma}\label{all1:1}
		For $b \in \mathbb{R}$, the intersection of $\Sigma_{p,\Gamma}$ with the line parallel to the diagonal passing through $(b, 0)$ consists precisely of points of the form
		$\bigl(b + \tilde{I}_b(u),\, \tilde{I}_b(u)\bigr)$,
		where $u$ is a critical point of $\tilde{I}_b$.
	\end{lemma}
	\begin{prop}\label{ap:1}
		Let $b\ge0$. Then $\phi_{1,p}$ is a global minimizer of $\tilde I_b$, with
		$\tilde{I}_b(\phi_{1,p}) = \lambda_{1,p}(\Gamma,V_D) - b$.
		The corresponding point in $\Sigma_{p,\Gamma}$ is
		$\bigl(\lambda_{1,p}(\Gamma,V_D),\, \lambda_{1,p}(\Gamma,V_D) - b\bigr)$,
		which lies on the vertical line passing through
		$\bigl(\lambda_{1,p}(\Gamma,V_D),\, \lambda_{1,p}(\Gamma,V_D)\bigr)$.
	\end{prop}
	\begin{proof}
		For $u\in M$, $\tilde I_b(u)\ge\lambda_{1,p}(\Gamma,V_D)\int_\Gamma|u|^p\,dx-b\int_\Gamma(u^+)^p\,dx\ge\lambda_{1,p}(\Gamma,V_D)-b$, with equality for $u=\phi_{1,p}$. \qed
\end{proof}

	We show that $-\phi_{1,p}$ is a strict local minimum of the constrained functional $\tilde{I}_b$.
	\begin{lemma}{\label{al:2}}
		Let $\Gamma$ be a compact connected quantum graph and let
		$\psi_n \not\equiv 0$, $\psi_n \in X(\Gamma,V_D)$ be such that $\psi_{n} \ge 0$ a.e.\ on $\Gamma$, and
		$\ell\{x \in \Gamma : \psi_n(x)>0\} \to 0$. Then
		\begin{align*}
			\frac{\int_\Gamma |\psi'_n|^p \, dx}{\int_\Gamma {\psi_n}^p \, dx}
			\rightarrow +\infty
			\quad \text{as } n \rightarrow \infty.
		\end{align*}
	\end{lemma}
	\begin{proof}
		Put $u_n=\psi_n/\|\psi_n\|_p$. If $\int_\Gamma|u_n'|^p\,dx$ were bounded along a subsequence, then, by Theorem~\ref{at:3}, a further subsequence would converge in $L^p(\Gamma)$ to some $u\ge0$ with $\|u\|_p=1$. Then $\ell\{u>\epsilon\}=\eta>0$ for some $\epsilon>0$, and $\ell\{u_n>\epsilon/2\}>\eta/2$ for large $n$ along that subsequence, contrary to $\ell\{\psi_n>0\}\to0$. \qed
\end{proof}

	\begin{prop}{\label{ap:2}}
		Let $b\ge0$. Then $-\phi_{1,p}$ is a strict local minimum of $\tilde{I}_b$ with $\tilde{I}_b(-\phi_{1,p}) = \lambda_{1,p}(\Gamma,V_D)$. The associated point in $\Sigma_{p,\Gamma}$ is $\bigl(\lambda_{1,p}(\Gamma,V_D)+b,\,\lambda_{1,p}(\Gamma,V_D)\bigr)$, which lies on the horizontal line passing through $\bigl(\lambda_{1,p}(\Gamma,V_D),\,\lambda_{1,p}(\Gamma,V_D)\bigr)$.
	\end{prop}

	\begin{proof}
		Suppose for contradiction that there exists a sequence $u_n \in M$ with $u_n \neq -\phi_{1,p}$, $u_n \rightarrow -\phi_{1,p}$ in $X(\Gamma,V_D)$, and $\tilde{I}_b(u_n) \le \lambda_{1,p}(\Gamma,V_D)$. We claim that $u_n$ changes sign for sufficiently large $n$. Since $u_n \rightarrow -\phi_{1,p}$ and $\phi_{1,p} > 0$ on $\Gamma\setminus V_D$, the function $u_n$ must be negative for some $x \in \Gamma$ when $n$ is large enough. If $u_n \le 0$ a.e.\ on $\Gamma$, then
		\begin{align*}
			\tilde{I}_b(u_n)
			= \int_{\Gamma} |u'_n|^p \, dx - b \int_\Gamma (u^{+}_{n})^p \,dx
			= \int_{\Gamma} |u'_n|^p \, dx
			> \lambda_{1,p}(\Gamma,V_D),
		\end{align*}
		the strict inequality holding by Lemma~\ref{al:100}\textit{(ii)}, as $u_n \neq \pm\phi_{1,p}$ and $\pm\phi_{1,p}$ are the only minimizers of \eqref{a1:3} in $M$. This contradicts $\tilde{I}_b(u_n) \le \lambda_{1,p}(\Gamma,V_D)$ and proves that $u_n$ changes sign. Setting $t_n = {\int_\Gamma |(u^{+}_n)'|^p \, dx}/{\int_\Gamma {(u^{+}_n)}^p \, dx}$,
		we compute
		\begin{align*}
			\tilde{I}_b(u_n)
			&= t_n \int_\Gamma {(u^{+}_n)}^p \, dx
			+ \int_\Gamma |(u^{-}_n)'|^p \, dx
			- b\int_\Gamma {(u^{+}_n)}^p \, dx \\[0.2em]
			&\ge (t_n - b)\int_\Gamma {(u^{+}_n)}^p \, dx
			+ \lambda_{1,p}(\Gamma,V_D) \int_\Gamma |(u^{-}_n)|^p \, dx.
		\end{align*}
		As $u_n \in M$ and $u_n = u^{+}_n - u^{-}_n$,
		\begin{align*}
			\tilde{I}_b(u_n) \le \lambda_{1,p}(\Gamma,V_D)
			= \lambda_{1,p}(\Gamma,V_D) \int_\Gamma (u^{+}_n)^p \, dx
			+ \lambda_{1,p}(\Gamma,V_D) \int_\Gamma |(u^{-}_n)|^p \, dx.
		\end{align*}
		Combining both inequalities and using $\int_\Gamma (u^+_n)^p \, dx > 0$, we deduce
		\begin{align*}
			\bigl(\lambda_{1,p}(\Gamma, V_D) - t_n + b\bigr)
			\int_\Gamma (u^{+}_n)^p \, dx \ge 0,
		\end{align*}
		and therefore
		\begin{equation}{\label{a:16}}
			t_n - b \le \lambda_{1,p}(\Gamma, V_D).
		\end{equation}
		Finally, $u_n \rightarrow -\phi_{1,p}$ uniformly on $\Gamma$ (Theorem~\ref{at2:3}) and $\phi_{1,p}>0$ on $\Gamma\setminus V_D$, the removed set being finite, so $\ell\{x : u_n(x) > 0\} \rightarrow 0$ and, by Lemma~\ref{al:2},
		$t_n \rightarrow \infty$ as $n \rightarrow \infty$, which contradicts
		\eqref{a:16}. Hence $-\phi_{1,p}$ is a strict local minimum of
		$\tilde{I}_b$. \qed
	\end{proof}

	To find a further critical point, we apply the Mountain Pass Theorem for a $C^1$ functional restricted to a $C^1$ manifold \cite[Theorem 2.1]{cuesta2003minimax}, with the Palais--Smale condition understood in the norm \eqref{a:17}. The space $X(\Gamma,V_D)$ is uniformly convex, as required there. We verify the Palais--Smale condition and the geometric assumptions for $\tilde{I}_b$ restricted to $M$. The norm of the derivative of the restriction $\tilde{I}_b$ of $I_b$ at $u \in M$ is
	\begin{equation}{\label{a:17}}
		\|\tilde{I}'_b(u)\|_{*}
		= \min\bigl\{\|I'_b(u) - c J'(u)\|_{(X(\Gamma,V_D))'} : c \in \mathbb{R}\bigr\},
	\end{equation}
	where $\|\cdot\|_{(X(\Gamma,V_D))'}$ denotes the norm on the dual space of $X(\Gamma,V_D)$.
	\begin{definition}{\label{ad:4}}
		The functional $I_b$ is said to satisfy the Palais--Smale condition $(P.S.)$ on $M$ if every sequence $u_n \in M$ such that $I_b(u_n)$ is bounded and $\|\tilde{I}'_b(u_n)\|_{*} \rightarrow 0$ admits a convergent subsequence.
	\end{definition}

	For $p>1$, the map $u\mapsto\int_\Gamma|u|^p\,dx$ is of class $C^1$ on $X(\Gamma,V_D)$ with derivative $p\int_\Gamma|u|^{p-2}u\,\cdot\,dx$, which is nonzero on $M$. Hence $M$ is a complete $C^1$ submanifold of $X(\Gamma,V_D)$, and for $u\in M$ the tangent space is
	\begin{align*}
		T_uM=\Bigl\{\phi\in X(\Gamma,V_D):\int_\Gamma|u|^{p-2}u\,\phi\,dx=0\Bigr\}.
	\end{align*}
	The number $\|\tilde I_b'(u)\|_*$ defined in \eqref{a:17} is the norm of the restriction of $I_b'(u)$ to $T_uM$, and $u$ is a critical point of $\tilde I_b$ if and only if $\|\tilde I_b'(u)\|_*=0$, which by the Lagrange multiplier rule is equivalent to $u$ solving \eqref{a1:1} with the pair $(b+\tilde I_b(u),\tilde I_b(u))$. Palais--Smale sequences and the Palais--Smale condition are understood with respect to $\|\cdot\|_*$, as in Definition \ref{ad:4}.

	\begin{lemma}{\label{al:3}}
		The functional $I_b$ satisfies the $(P.S.)$ condition on $M$.
	\end{lemma}
	\begin{proof}
Let $u_n\in M$ satisfy $|I_b(u_n)|\le K$ and $\|\tilde I_b'(u_n)\|_*\to0$, and let $c_n$ attain the minimum in \eqref{a:17} for each $u_n$. Since $\|u_n\|_p=1$ and
\[
\int_\Gamma |u_n'|^p\,dx\le K+|b|,
\]
the sequence $(u_n)$ is bounded in $X(\Gamma,V_D)$. Moreover, testing $I_b'(u_n)-c_nJ'(u_n)$ with $u_n$ shows that $(c_n)$ is bounded. By Theorem~\ref{at:3}, passing to a subsequence, we have
$u_n\rightharpoonup u_0 \text{ in }X(\Gamma,V_D), \text{ and }
u_n\to u_0 \text{ in }L^p(\Gamma)$.
Now, testing $I_b'(u_n)-c_nJ'(u_n)$ with $u_n-u_0$ and using H{\"o}lder's inequality, we obtain
\[
\begin{aligned}
\int_\Gamma (u_n^+)^{p-1}|u_n-u_0|\,dx
+\int_\Gamma |u_n|^{p-1}|u_n-u_0|\,dx
\le 2\|u_n-u_0\|_p\to0.
\end{aligned}
\]
Consequently, $\int_\Gamma |u_n'|^{p-2}u_n'(u_n'-u_0')\,dx\to0$.
Hence, by Lemma~\ref{n:splus}, we conclude that $u_n\to u_0 \text{ in }X(\Gamma,V_D)$. \qed
\end{proof}

	Let $c_e>0$ be such that $\|w\|_p\le c_e\|w\|$ for $w\in X(\Gamma,V_D)$ (Theorem~\ref{at:3}), and for $x\in M$ and $r>0$ put $B_r(x)=\{u\in M:\|u-x\|<r\}$.
	\begin{lemma}\label{n:proj}
		Let $x\in M$, $0<r\le 1/(2c_e)$ and $u,v\in B_r(x)$, and put $z_t=(1-t)u+tv$ for $t\in[0,1]$. Then $w_t=z_t/\|z_t\|_p$ is a path in $M$ from $u$ to $v$, and
		\begin{equation}\label{a10:4}
			\|w_t-x\|\le C_0r,\qquad C_0=2\bigl(1+c_e\|x\|\bigr),\qquad t\in[0,1].
		\end{equation}
		Consequently, the path components of an open subset of $M$ are open.
	\end{lemma}
	\begin{proof}
		By convexity $\|z_t-x\|<r$, and since $\|u\|_p=\|v\|_p=\|x\|_p=1$, we have $\bigl|\,\|z_t\|_p-1\,\bigr|\le\|z_t-x\|_p\le c_er\le\frac12$. Hence $w_t$ is well defined, lies in $M$ and depends continuously on $t$, and
		\begin{align*}
			\|w_t-x\|=\frac{\|z_t-x+(1-\|z_t\|_p)x\|}{\|z_t\|_p}\le 2\bigl(r+c_er\|x\|\bigr)=C_0r .
		\end{align*}
		If $U\subseteq M$ is open and $x\in U$, choose $\rho>0$ with $B_\rho(x)\subseteq U$ and $r=\min\{\rho/C_0,1/(2c_e)\}$; then every point of $B_r(x)$ is joined to $x$ by a path in $B_{C_0r}(x)\subseteq U$, and the last assertion follows. \qed
	\end{proof}

We shall use the following three facts of critical point theory on $M$. The deformation lemma on a manifold $\{\psi=1\}$ in \cite[Lemma~5.15]{willem1996minimax} requires $\psi\in C^2$, while $u\mapsto\int_\Gamma|u|^p\,dx$ is of class $C^2$ only for $p\ge2$, and it gives no bound on the displacement, which is used in Section~\ref{as:10}. We therefore deduce the first fact from \cite[Lemma~2.3]{willem1996minimax} for $C^1$ functionals on a Banach space, passing through a $0$-homogeneous extension of $\tilde I_b$. The second and third facts are deduced from the first.

	\begin{lemma}[Quantitative deformation lemma]\label{al:deform}
		Let $b\in\mathbb R$, $c\in\mathbb R$, $0<\varepsilon\le1$ and $S\subseteq M$. Put $K=(3+|b|+|c|)^{1/p}$ and $C_*=2(1+c_eK)$, and let $0<\delta\le1/(16c_e)$ satisfy
		\begin{equation}\label{a:defhyp}
			\|\tilde I_b'(u)\|_*\ge\frac{10\,\varepsilon}{\delta}\qquad\text{for all }u\in M\text{ with }\operatorname{dist}(u,S)\le2C_*\delta\text{ and }|\tilde I_b(u)-c|\le2\varepsilon .
		\end{equation}
		Then there is a continuous map $\eta:[0,1]\times M\to M$ such that
		\begin{equation}\label{a:defconc}
			\begin{aligned}
				&\eta(0,u)=u,\qquad \eta(t,u)=u\ \text{ if }|\tilde I_b(u)-c|>2\varepsilon\text{ or }\operatorname{dist}(u,S)>2\delta,\\
				&\|\eta(t,u)-u\|\le C_*\delta,\qquad t\mapsto\tilde I_b(\eta(t,u))\ \text{ is non-increasing},\\
				&\eta\bigl(1,\{\tilde I_b\le c+\varepsilon\}\cap S\bigr)\subseteq\{\tilde I_b\le c-\varepsilon\}.
			\end{aligned}
		\end{equation}
	\end{lemma}

	\begin{proof}
For $u\neq0$, define $R(u)=I_b(u)/J(u)$.
Then $R\in C^1(X(\Gamma,V_D)\setminus\{0\})$ and
\[
R(tu)=R(u), \qquad R'(tu)=t^{-1}R'(u)
\quad\text{for all }t>0.
\]
Since $R=\tilde I_b$ on $M$, for every $u\in M$,
we have $R'(u)=I_b'(u)-I_b(u)J'(u)$,
and hence, by \eqref{a:17},
$\|R'(u)\|\ge\|\tilde I_b'(u)\|_*$.
Choose $\chi\in C^\infty(\mathbb R)$ such that
$\chi=0$ on $(-\infty,4^{-p}]$ and
$\chi=1$ on $[2^{-p},\infty)$.
Define $\varphi:X(\Gamma,V_D)\to\mathbb R$ by
\[
\varphi(u)=
\begin{cases}
\chi(J(u))R(u), & u\neq0,\\
0, & u=0.
\end{cases}
\]
Then $\varphi\in C^1(X(\Gamma,V_D),\mathbb R)$ and $\varphi=R$ on
$W=\{u:\|u\|_p>1/2\}$. We verify the hypothesis of \cite[Lemma~2.3]{willem1996minimax} for $\varphi$, $S$, $c$, $\varepsilon$, $\delta$. Let $w\in X(\Gamma,V_D)$ satisfy $\operatorname{dist}(w,S)\le2\delta$ and $|\varphi(w)-c|\le2\varepsilon$, and choose $s\in S$ with $\|w-s\|\le3\delta$. Then $\bigl|\,\|w\|_p-1\,\bigr|\le\|w-s\|_p\le3c_e\delta\le\frac3{16}$, so $w\in W$, and $\hat w=w/\|w\|_p\in M$ satisfies $\tilde I_b(\hat w)=R(w)=\varphi(w)\le c+2$. Hence
\begin{align*}
\|\hat w\|^p=\tilde I_b(\hat w)+b\int_\Gamma(\hat w^+)^p\,dx+1\le K^p,
\end{align*}
so that $\|\hat w-w\|=\bigl|\,\|w\|_p-1\,\bigr|\,\|\hat w\|\le3c_eK\delta$ and $\operatorname{dist}(\hat w,S)\le3\delta+3c_eK\delta\le2C_*\delta$. By \eqref{a:defhyp},
\begin{align*}
\|\varphi'(w)\|=\|R'(w)\|=\|w\|_p^{-1}\|R'(\hat w)\|\ge\frac{16}{19}\cdot\frac{10\,\varepsilon}{\delta}>\frac{8\varepsilon}{\delta}.
\end{align*}
Let $\tilde\eta\in C([0,1]\times X(\Gamma,V_D),X(\Gamma,V_D))$
be the map given by \cite[Lemma~2.3]{willem1996minimax},
which satisfies
\begin{align*}
&\tilde\eta(t,u)=u
\quad\text{if }t=0,\
|\varphi(u)-c|>2\varepsilon
\text{ or }\operatorname{dist}(u,S)>2\delta,\\
&\tilde\eta(1,\{\varphi\le c+\varepsilon\}\cap S)
\subseteq\{\varphi\le c-\varepsilon\},\\
&\|\tilde\eta(t,u)-u\|\le\delta,
\end{align*}
and $t\mapsto\varphi(\tilde\eta(t,u))$ is
non-increasing. For $u\in M$, we have
$\|\tilde\eta(t,u)\|_p
\ge1-c_e\delta\ge\frac{15}{16}$,
and we put
\[
\eta(t,u)=
\frac{\tilde\eta(t,u)}
{\|\tilde\eta(t,u)\|_p}\in M.
\]
Then $\eta$ is continuous,
$\tilde\eta(t,u)\in W$, and
$\tilde I_b(\eta(t,u))
=\varphi(\tilde\eta(t,u))$.
This gives the first, fourth and fifth properties
in \eqref{a:defconc}. If $\eta(t,u)\neq u$, then
$\tilde I_b(\eta(t,u))
\le\tilde I_b(u)\le c+2\varepsilon$,
so $\|\eta(t,u)\|\le K$ as above. Consequently,
\begin{align*}
\|\eta(t,u)-u\|
\le\|\tilde\eta(t,u)-u\|
+\bigl|\,\|\tilde\eta(t,u)\|_p-1\,\bigr|
\|\eta(t,u)\|
\le\delta+c_eK\delta
\le C_*\delta. \qquad \qed
\end{align*}
\end{proof}

	\begin{lemma}\label{al:unifgrad}
		Let $[a,a']\subset\mathbb R$ contain no critical value of $\tilde I_0$. Then there are $\varepsilon>0$ and $\beta>0$ such that $[a-2\varepsilon,a'+2\varepsilon]$ contains no critical value and
		\begin{align*}
			\|\tilde I_0'(u)\|_*\ge\beta\qquad\text{for all }u\in M\text{ with }a-2\varepsilon\le\tilde I_0(u)\le a'+2\varepsilon .
		\end{align*}
	\end{lemma}

	\begin{proof}
		The set of critical values of $\tilde I_0$ in a compact interval is closed. Indeed, if $c_n\to c$ are critical values with critical points $u_n\in M$, then $\tilde I_0(u_n)=c_n$ is bounded and $\|\tilde I_0'(u_n)\|_*=0$, so $(u_n)$ is a Palais--Smale sequence and Lemma \ref{al:3} provides a subsequence converging to a critical point at level $c$. Hence the complement of the set of critical values is open and contains the compact set $[a,a']$, so it contains $[a-2\varepsilon,a'+2\varepsilon]$ for some $\varepsilon>0$. If no $\beta>0$ existed for this $\varepsilon$, there would be $u_n\in M$ with $\tilde I_0(u_n)\in[a-2\varepsilon,a'+2\varepsilon]$ and $\|\tilde I_0'(u_n)\|_*\to0$; passing to a subsequence with $\tilde I_0(u_n)\to c$, Lemma \ref{al:3} would produce a critical point at the level $c\in[a-2\varepsilon,a'+2\varepsilon]$, a contradiction. \qed
	\end{proof}

	\begin{lemma}\label{al:noncrit}
		Let $[a,a']$ contain no critical value of $\tilde I_0$. Then for every $u\in M$ with $\tilde I_0(u)\le a'$ there is a continuous path $\gamma:[0,1]\to M$ with
		\begin{align*}
			\gamma(0)=u,\qquad \tilde I_0(\gamma(1))\le a,\qquad t\mapsto\tilde I_0(\gamma(t))\ \text{non-increasing}.
		\end{align*}
		In particular $\gamma[0,1]\subseteq\{\tilde I_0\le\tilde I_0(u)\}$.
	\end{lemma}

	\begin{proof}
		We may assume $\tilde I_0(u)>a$. Let $\varepsilon$ and $\beta$ be as in Lemma \ref{al:unifgrad}. Let $N$ be an integer with $N\ge(a'-a)/(2\varepsilon)$, so large that $(a'-a)/(2N)\le1$ and $10(a'-a)/(2N\beta)\le1/(16c_e)$, and replace $\varepsilon$ by $(a'-a)/(2N)\le\varepsilon$; the interval $[a-2\varepsilon,a'+2\varepsilon]$ then shrinks, so the conclusion of Lemma~\ref{al:unifgrad} persists with the same $\beta$. Put $\delta=10\varepsilon/\beta$ and $S=M$. Then, with $b=0$, \eqref{a:defhyp} holds at every level $c$ for which $[c-2\varepsilon,c+2\varepsilon]\subseteq[a-2\varepsilon,a'+2\varepsilon]$. Set $c_j=a'-(2j-1)\varepsilon$ for $j=1,\dots,N$, so that each slab $|\tilde I_0-c_j|\le2\varepsilon$ lies in $[a-2\varepsilon,a'+2\varepsilon]$ and Lemma \ref{al:deform} applies at $c=c_j$, giving a deformation $\eta_j$. Since $c_1+\varepsilon=a'$ and $c_j+\varepsilon=c_{j-1}-\varepsilon$ for $j\ge2$,
		\begin{align*}
			\eta_N(1,\cdot)\circ\dots\circ\eta_1(1,\cdot)\ \text{maps}\ \{\tilde I_0\le a'\}\ \text{into}\ \{\tilde I_0\le a'-2N\varepsilon\}=\{\tilde I_0\le a\}.
		\end{align*}
		Put $u_0=u$ and $u_j=\eta_j(1,u_{j-1})$, and define $\gamma$ on $[(j-1)/N,j/N]$ by $\gamma(t)=\eta_j(Nt-j+1,u_{j-1})$. Then $\gamma$ is continuous, $\gamma(0)=u$, $\tilde I_0(\gamma(1))=\tilde I_0(u_N)\le a$, and $\tilde I_0$ is non-increasing along each piece by \eqref{a:defconc}, hence along $\gamma$. \qed
	\end{proof}

	\begin{lemma}{\label{al:4}}
		Let $\delta_0 > 0$ be such that
		\begin{align*}
			\tilde{I_b}(u) > \tilde{I_b}(-\phi_{1,p}), \quad \forall u \in B_{\delta_0}(-\phi_{1,p})\setminus\{-\phi_{1,p}\},
		\end{align*}
		with $B_{\delta_0}(-\phi_{1,p})$ as in Lemma~\ref{n:proj}. Then, for any $0<\delta<\delta_0$,
		\begin{equation}{\label{a:21}}
			\inf \{ \tilde{I_{b}}(u) : u \in M , ||u-(-\phi_{1,p})|| = \delta \} > \tilde{I_b}(-\phi_{1,p}).
		\end{equation}

	\end{lemma}

	\begin{proof}
		The argument is that of \cite{cuesta1999beginning}, and we include it for completeness. Suppose, for contradiction, that there exists a $\delta \in (0,\delta_0)$ for which equality holds in \eqref{a:21}.
	Then there exists a sequence $u_n \in M$ with $||u_n-(-\phi_{1,p})|| = \delta$ and $\tilde{I_{b}}(u_n) \le \lambda_{1,p}(\Gamma, V_D) + 1/2n^2$. Consider
	\begin{align*}
		M_{\eta} = \{ u \in M : \delta - \eta \le ||u-(-\phi_{1,p})||  \le \delta +\eta \},
	\end{align*}
	where $\eta >0$ is chosen so that $0<\delta - \eta < \delta + \eta < \delta_0 $. By the contradiction hypothesis and $\tilde{I_b}(u) > \tilde{I_b}(-\phi_{1,p})$, we have $\inf \{ \tilde{I_{b}}(u) : u\in M_{\eta}\} = \lambda_{1,p}(\Gamma,V_D)$. Ekeland's variational principle applied to $\tilde{I_{b}}$ on $M_{\eta}$ provides, for each $n \in \mathbb{N}$, some $w_n \in M_{\eta}$ such that
	\begin{align}
		\tilde{I_b}(w_n) &\le \tilde{I_b}(u_n), \label{a:23}\\[0.3em]
		||w_n - u_n|| &\le 1/n,\label{a:24}\\[0.3em]
		\tilde{I_b}(w_n) &\le \tilde{I_b}(u)
		+ \frac{1}{n} ||u-w_n||, \quad \forall u\in M_{\eta}.\label{a:25}
	\end{align}
	We show that $w_n$ is a $(P.S.)$ sequence for $\tilde{I_{b}}$ on $M$. By Definition \ref{ad:4}, if $\tilde{I_{b}}(w_n)$ is bounded (by \eqref{a:23}) and $||{I'_b}(w_n)||_{*} \rightarrow 0$, then Lemma \ref{al:3} provides a subsequence $w_{n_k} \to w$ strongly in $X(\Gamma, V_D)$. We claim that $||{I'_b}(w_n)||_{*} \rightarrow 0$. Fix $n> 1/\eta$ and take $\phi \in X(\Gamma,V_D)$ tangent to $M$ at $w_n$, that is, satisfying
	\begin{equation}{\label{a:26}}
		\int_{\Gamma} |w_n|^{p-2} w_n \phi \, dx  = 0,
	\end{equation}
	and consider, for each $t \in \mathbb{R}$,
	\begin{align*}
		u_t = \frac{w_n + t\phi}{||w_n + t\phi||_p}.
	\end{align*}
	Here $u_t$ is normalized by dividing by the $p$-norm of $w_n + t\phi$. For sufficiently small $|t|$, $u_t \in M_{\eta}$.
	Indeed
	\begin{align*}
		\lim _{t \to 0} ||u_t - (-\phi_{1,p})|| = ||w_n - (-\phi_{1,p})|| .
	\end{align*}
	By \eqref{a:24},
	\begin{align*}
		||w_n - (-\phi_{1,p})||
		\le ||w_n - u_n|| + ||u_n - (-\phi_{1,p})|| \\[0.5em]
		\le \frac{1}{n} + \delta < \eta +\delta,
	\end{align*}
	and on the other hand
	\begin{align*}
		||w_n - (-\phi_{1,p})||
		\ge ||(-\phi_{1,p})- (u_n)||  - ||u_n - w_n|| \\[0.5em]
		\ge  \delta - \frac{1}{n} > \delta -\eta.
	\end{align*}
	Replacing $u$ by $u_t$ in \eqref{a:25} and
	writing $\beta(t) = ||w_n + t\phi||_p$, we obtain, for $t>0$,
	\begin{equation}{\label{a4:26}}
		\frac{I_b(w_n)-I_b(w_n+t\phi)}{t}
		\le
		\frac{1}{t}\left(\frac{1}{\beta(t)^p}-1\right) I_b(w_n+t\phi)
		+\frac{1}{nt}\frac{1}{\beta(t)}
		\|w_n(1-\beta(t))+t\phi\|.
	\end{equation}

	The first term on the right-hand side of \eqref{a4:26} contains $({\beta(t)^p} - 1)/t$, which is equal to the derivative of $\beta(t)^p$ near $t=0$. Using \eqref{a:26},
	\begin{align*}
		\left.\frac{d}{dt}\beta(t)^p\right|_{t=0} =
		\left.\frac{d}{dt}\left(\int_{\Gamma} |w_n + t\phi|^p \, dx\right)\right|_{t=0}
		=p \int_{\Gamma} |w_n|^{p-2} w_n \phi \, dx = 0 .
	\end{align*}

	Hence $({\beta(t)^p} - 1)/t \to 0$ as $t \to 0$. Similarly, the second term on the right-hand side of \eqref{a4:26} involves $(1- \beta(t))/t$, which also tends to zero as $t \to 0$ by \eqref{a:26}. Letting $t \to 0$ in \eqref{a4:26} and repeating the argument with $-\phi$ in place of $\phi$, we obtain
	\begin{equation}{\label{a:27}}
		\left|\langle I'_b(w_n), \phi \rangle \right|\le \frac{1}{n} ||\phi||,
	\end{equation}
	for all $\phi\in X(\Gamma,V_D)$ tangent to $M$ at $w_n$. For arbitrary $\phi$ in $X(\Gamma,V_D)$, choose $a_n$ such that $\phi -a_n w_n$ satisfies \eqref{a:26},
		that is, $a_n = \int_{\Gamma} |w_n|^{p-2}w_n \phi \, dx$. Replacing $\phi$ by $\phi-a_n w_n$ in \eqref{a:27}, we deduce
		\begin{align*}
			\left|\langle I'_b(w_n),\phi \rangle - \langle I'_b(w_n),w_n\rangle \int_{\Gamma} |w_n|^{p-2}w_n \phi \, dx \right| \le \frac{1}{n} ||\phi - a_n w_n||
		\end{align*}
		With $c_n = \langle I'_b(w_n), w_n \rangle$ and $\| a_n w_n \| \le K \| \phi \|$, we obtain, for some $\delta_n \to 0$,
		\begin{align*}
			\left|\langle I'_b(w_n),\phi \rangle - c_n \int_{\Gamma} |w_n|^{p-2}w_n \phi \, dx \right| \le \delta_n  ||\phi||.
		\end{align*}
		Thus $||{\tilde{I'}_b}(w_n)||_{*} \rightarrow 0$, and hence $w_n$ is a $(P.S.)$ sequence for ${\tilde{I}_{b}}$ on $M$. By the $(P.S.)$ condition for $\tilde{I}_b$ on $M$, there exists a subsequence $w_{n_{k}} \to w$ strongly in $X(\Gamma,V_D)$. Then $w \in M$ with $\| w-(-\phi_{1,p})\| = \delta$ and $\tilde{I_{b}}(w) = \lambda_{1,p}(\Gamma,V_D)$, which contradicts $\tilde{I_{b}}(u)> \tilde{I_{b}}(-\phi_{1,p})$. \qed

	\end{proof}

	A version of the Mountain Pass Theorem now provides a further critical point of $\tilde{I_{b}}$.

	\begin{prop}{\label{ap:3}}
		Let $b\ge0$.
		Then the set
    \[
        \mathbb{P} = \{P \in C([-1,1], M) : P(-1) = -\phi_{1,p}, \text{ } P(1) = \phi_{1,p} \}
	\]
        is non-empty and
		\begin{equation}{\label{a:28}}
			\alpha(b) = \inf _{P \in \mathbb{P}} \max_{u \in P[-1,+1]} \tilde{I_b}(u)
		\end{equation}
		is a critical value of $\tilde{I_{b}}.$ Moreover,
		$\alpha(b) > \max \{ \tilde{I_{b}}(-\phi_{1,p}),\tilde{I_{b}}(\phi_{1,p})\} = \lambda_{1,p}(\Gamma,V_D).$
	\end{prop}

	\begin{proof}
		To see that $\mathbb{P}$ is non-empty, take $\psi \in X(\Gamma,V_D)$ such that $\psi \notin span (\phi_{1,p})$ and consider the path through $\psi$ and $\phi_{1,p}$
		\begin{align*}
			{P}(t) = \frac{t\phi_{1,p} + (1-|t|)\psi}{\| t\phi_{1,p} + (1-|t|)\psi \|_p}\qquad t\in[-1,1],
		\end{align*}
		which is well defined since $t\phi_{1,p}+(1-|t|)\psi\neq0$. By Proposition~\ref{ap:2} and Lemma~\ref{al:4} there is $\delta\in(0,2\|\phi_{1,p}\|)$ with
			\begin{align*}
				\inf \{ \tilde{I_{b}}(u) : u \in M , \|u+\phi_{1,p}\| = \delta \} > \lambda_{1,p}(\Gamma,V_D)= \max \{ \tilde{I_{b}}(-\phi_{1,p}),\tilde{I_{b}}(\phi_{1,p})\}.
		\end{align*}
		The $(P.S)$ condition holds by Lemma \ref{al:3}, and the mountain pass geometry by the inequality above. Hence $\alpha (b)$ is a critical value of $\tilde{I_{b}}$. Since every path in $\mathbb P$ meets the sphere $\{\|u+\phi_{1,p}\|=\delta\}$, the inequality above also gives $\alpha(b) > \max \{ \tilde{I_{b}}(-\phi_{1,p}),\tilde{I_{b}}(\phi_{1,p})\} = \lambda_{1,p}(\Gamma,V_D)$.     \qed
	\end{proof}

	With the critical value $\alpha(b)$ defined in \eqref{a:28}, we have proved the following.

	\begin{prop}{\label{at4:5}}
		For every $b\ge 0$, the point $(b +\alpha(b),\alpha(b))$, where $\alpha(b) >  \lambda_{1,p}(\Gamma, V_D)$  is defined by the minimax formula \eqref{a:28}, lies in $\Sigma_{p,\Gamma}$.
	\end{prop}

	Since $\Sigma_{p,\Gamma}$ is symmetric with respect to the diagonal, the entire curve obtained in Proposition \ref{at4:5} is given by
	\begin{equation}{\label{a:29}}
		\mathfrak{C} = \{ (b +\alpha(b),\alpha(b)), (\alpha(b),b+ \alpha(b) ) : b\ge 0 \}.
	\end{equation}

	\subsection{Isolatedness of the trivial lines and minimality of $\mathfrak C$}{\label{as:5}}

	We first establish that the lines $\mathbb{R} \times \{\lambda_{1,p}(\Gamma,V_D)\}$ and
	$\{\lambda_{1,p}(\Gamma,V_D)\} \times \mathbb{R}$ are isolated in $\Sigma_{p,\Gamma}$. We then record some topological properties of the functional $\tilde{I}_b$ along with several auxiliary lemmas. The main purpose of this section is to show that the points of $\mathfrak{C}$ defined in \eqref{a:29} are the first nontrivial points for $b \ge 0$, lying in the intersection of $\Sigma_{p,\Gamma}$ with the line $(b,0) + c(1,1)$ for $c \in \mathbb{R}$.

	\begin{prop}{\label{ap:4}}
		The lines $\mathbb{R} \times \{\lambda_{1,p}(\Gamma,V_D)\}$ and $\{\lambda_{1,p}(\Gamma,V_D)\} \times \mathbb{R}$ are isolated in $\Sigma_{p,\Gamma}$, that is, there exists no sequence $(\lambda^{+}_n, \lambda^{-}_n) \in \Sigma_{p,\Gamma}$ with $\lambda^{+}_n, \lambda^{-}_n > \lambda_{1,p}(\Gamma,V_D)$ such that $(\lambda^{+}_n,\lambda^{-}_n) \to (\lambda^+,\lambda^-)$ with $\lambda^{+} = \lambda_{1,p}(\Gamma,V_D)$ or
		$\lambda^{-} = \lambda_{1,p}(\Gamma,V_D)$.
	\end{prop}

	\begin{proof}
		Suppose, for contradiction, that there exists a sequence
		$(\lambda^{+}_n,\lambda^{-}_n) \in \Sigma_{p,\Gamma}$ with
		$\lambda^{+}_n, \lambda^{-}_n > \lambda_{1,p}(\Gamma,V_D)$ and $(\lambda^{+}_n,\lambda^{-}_n) \to (\lambda^+,\lambda^-)$ with $\lambda^{+} = \lambda_{1,p}(\Gamma,V_D)$ or
		$\lambda^{-} = \lambda_{1,p}(\Gamma,V_D)$. Let $u_n \in X(\Gamma,V_D)$ be a solution of
		\begin{equation}{\label{a:30}}
			-\bigl(|u'_n|^{p-2}u'_n\bigr)'
			= \lambda^+_n (u^+_n)^{p-1} - \lambda^-_n (u^-_n)^{p-1}
			\quad \text{on the edges of } \Gamma,
		\end{equation}
		normalized by $\|u_n\|_p = 1$. Then
		\begin{align*}
			\int_\Gamma |u'_n|^p\,dx
			= \lambda^+_n \int_\Gamma (u^+_n)^p \, dx
			+ \lambda^-_n \int_\Gamma (u^-_n)^p \, dx
			\le \lambda^+_n + \lambda^-_n
		\end{align*}
		which implies that $\{u_n\}$ is bounded in $X(\Gamma,V_D)$ and, up to
		a subsequence, $u_{n_k} \rightharpoonup u$ weakly in $X(\Gamma,V_D)$ and
		$u_{n_k} \to u$ strongly in $L^p(\Gamma)$. Testing the equation of $u_{n_k}$ with $u_{n_k}-u$ and using H\"older's inequality for the terms of order zero, we obtain $\limsup_k\int_\Gamma|u_{n_k}'|^{p-2}u_{n_k}'(u_{n_k}'-u')\,dx\le0$, so that $u_{n_k}\to u$ in $X(\Gamma,V_D)$ by Lemma~\ref{n:splus}. In the limit, the weak formulation shows that $u$ satisfies
		\begin{equation}{\label{a:31}}
			-\bigl(|u'|^{p-2}u'\bigr)'
			= \lambda^+ (u^+)^{p-1} - \lambda^- (u^-)^{p-1}
			\quad \text{on the edges of } \Gamma.
		\end{equation}

			Since $\lambda_n^\pm>\lambda_{1,p}(\Gamma,V_D)$, the function $u_n$ changes sign: if, say, $u_n\ge0$, then $u_n$ is a nonnegative eigenfunction with eigenvalue $\lambda_n^+>\lambda_{1,p}(\Gamma,V_D)$, which is impossible by Lemma~\ref{al:100}\textit{(iii)}, and the case $u_n\le0$ is the same. Suppose that $\lambda^+=\lambda_{1,p}(\Gamma,V_D)$; the case $\lambda^-=\lambda_{1,p}(\Gamma,V_D)$ is symmetric. Multiplying \eqref{a:31} by $u^+$ and integrating over $\Gamma$, we obtain
			\begin{equation}\label{a:32}
				\int_\Gamma |(u^+)'|^p\,dx = \lambda_{1,p}(\Gamma,V_D)\int_\Gamma (u^+)^{p}\,dx.
			\end{equation}
			If $u^+\not\equiv0$, then $u^+/\|u^+\|_p$ attains the infimum \eqref{a1:3} and is therefore an eigenfunction associated with $\lambda_{1,p}(\Gamma,V_D)$. By Lemma~\ref{al:100}\textit{(ii)}, $u^+>0$ on $\Gamma\setminus V_D$, so $u^-=0$ and $u=\phi_{1,p}$. If $u^+\equiv0$, then $-u$ is a nonnegative eigenfunction with eigenvalue $\lambda^-$, so $\lambda^-=\lambda_{1,p}(\Gamma,V_D)$ and $u=-\phi_{1,p}$ by Lemma~\ref{al:100}. Hence $u_{n_k}\to\phi_{1,p}$ or $u_{n_k}\to-\phi_{1,p}$ in $X(\Gamma,V_D)$ and, by Theorem~\ref{at2:3}, uniformly on $\Gamma$. Since $\phi_{1,p}>0$ on $\Gamma\setminus V_D$ and $V_D$ is finite, $\ell\{x\in\Gamma:\phi_{1,p}(x)<\delta\}\to0$ as $\delta\to0$, and therefore
			\begin{equation}\label{a:33}
				\text{either } \ell\{x \in \Gamma : u_{n_k}(x) < 0\} \to 0
				\quad \text{or} \quad
				\ell\{x \in \Gamma : u_{n_k}(x) > 0\} \to 0.
			\end{equation}
			On the other hand, multiplying \eqref{a:30} by $u^+_n$ and integrating, we obtain
			\begin{equation}\label{a:34}
				\int_\Gamma |(u^+_n)'|^p\,dx = \lambda^+_{n}\int_\Gamma (u^+_n)^{p}\,dx.
			\end{equation}
			Fix $q\in(p,\infty)$ and let $K>0$ be such that $\|u\|_q\le K\|u\|_{W^{1,p}(\Gamma)}$ (Theorem~\ref{at:3}). By H\"older's inequality and \eqref{a:34},
			\begin{align*}
				\int_\Gamma (u^+_n)^{p}\,dx &\le \bigl(\ell\{x \in \Gamma : u_n(x) > 0\}\bigr)^{1-\frac{p}{q}}\|u^+_n\|^p_q\\
				&\le \bigl(\ell\{x \in \Gamma : u_n(x) > 0\}\bigr)^{1-\frac{p}{q}}K^p\,(1+\lambda^+_n)\int_\Gamma (u^+_n)^{p}\,dx,
			\end{align*}
			where \eqref{a:34} was used in the form $\|u_n^+\|_{W^{1,p}(\Gamma)}^p=(1+\lambda_n^+)\int_\Gamma(u_n^+)^p\,dx$. As $u_n^+\not\equiv0$,
			\begin{align*}
				(1+\lambda^+_{n})^{-1}K^{-p} \le \bigl(\ell\{x \in \Gamma : u_n(x) > 0\}\bigr)^{1-\frac{p}{q}},
			\end{align*}
			and in the same way
			\begin{align*}
				(1+\lambda^-_{n})^{-1}K^{-p} \le \bigl(\ell\{x \in \Gamma : u_n(x) < 0\}\bigr)^{1-\frac{p}{q}}.
			\end{align*}
			The sequences $\lambda_n^\pm$ are bounded, so these two inequalities contradict \eqref{a:33}. This proves that the trivial lines are isolated in $\Sigma_{p,\Gamma}$. \qed
	\end{proof}

		\begin{lemma}\label{al:101}
			The set $\Sigma_{p,\Gamma}$ is closed in $\mathbb R^2$.
		\end{lemma}
		\begin{proof}
			Let $(\lambda_n^+,\lambda_n^-)\in\Sigma_{p,\Gamma}$ converge to $(\lambda^+,\lambda^-)$, and let $u_n$ be solutions of \eqref{a1:1} with $\|u_n\|_p=1$. Testing \eqref{a2:7} with $u_n$ gives $\int_\Gamma|u_n'|^p\,dx\le|\lambda_n^+|+|\lambda_n^-|$, so $(u_n)$ is bounded in $X(\Gamma,V_D)$. Up to a subsequence, $u_n\rightharpoonup u$ in $X(\Gamma,V_D)$ and $u_n\to u$ in $L^p(\Gamma)$, so $\|u\|_p=1$. Testing \eqref{a2:7} with $u_n-u$ and using the $(S_+)$ property as in the proof of Lemma~\ref{al:3}, we get $u_n\to u$ in $X(\Gamma,V_D)$, and passing to the limit in \eqref{a2:7} shows that $u$ is a nontrivial solution of \eqref{a1:1} with $(\lambda^+,\lambda^-)$. \qed
		\end{proof}

	\begin{lemma}{\label{al:5}}
		Let $\mathcal{U}\subseteq M$ be open. Then every connected component $\mathcal{U}'$ of $\mathcal{U}$ is open and arc-connected, and $\partial\mathcal{U}'\cap\mathcal{U}=\emptyset$.
	\end{lemma}

	\begin{proof}
		By Lemma~\ref{n:proj}, for $u\in M$ and $\rho>0$ the path component of $B_\rho(u)$ containing $u$ is an open path-connected neighbourhood of $u$. Hence $M$ is locally path-connected, the components of $\mathcal U$ are open and coincide with its path components, and they are arc-connected, since a path-connected Hausdorff space is arc-connected. If $u_0\in\partial\mathcal U'\cap\mathcal U$, let $\mathcal V\subseteq\mathcal U$ be an open path-connected neighbourhood of $u_0$. Then $\mathcal V$ meets $\mathcal U'$, so $\mathcal U'\cup\mathcal V$ is connected and $\mathcal V\subseteq\mathcal U'$. Thus $u_0\in\mathcal U'$, which is impossible, because $\mathcal U'$ is open and $u_0\in\partial\mathcal U'$. \qed
	\end{proof}

	In the next lemma, we identify critical points of $\tilde{I}_b$ in sets of the form
	\begin{align*}
		\mathcal{U} = \{u \in M : \tilde{I}_b(u) < r\}.
	\end{align*}

	\begin{lemma}{\label{al:6}}
		Every connected component of $\mathcal{U}$ contains a critical point of $\tilde{I}_b$.
	\end{lemma}

	\begin{proof}
		Let $\mathcal{U}_1$ be a connected component of $\mathcal{U}$ and set $h = \inf\{\tilde{I}_b(u) : u \in \overline{\mathcal{U}_1}\}$, where $\overline{\mathcal{U}_1}$ denotes the closure of $\mathcal{U}_1$. We claim that this infimum is attained at some $\bar{u} \in \overline{\mathcal{U}_1}$ satisfying $\tilde{I}_b(\bar{u}) = h < r$, and that $\bar{u} \in \mathcal{U}$. If so, either $\bar{u} \in \mathcal{U}_1$, and then $\bar{u}$ is a critical point of $\tilde{I}_b$, or $\bar{u} \in \partial\mathcal{U}_1$, which is impossible by Lemma~\ref{al:5}.

		To show that $h$ is attained, let $\{u_n\} \subset \mathcal{U}_1$ be a minimizing sequence satisfying $\tilde{I}_b(u_n) \le h + \frac{1}{2n^2}$. Applying Ekeland's variational principle to $\tilde{I}_b$ on
		$\overline{\mathcal{U}_1}$, for each $n \in \mathbb{N}$ we find $w_n \in \overline{\mathcal{U}_1}$ such that
		\begin{align}
			\tilde{I}_b(w_n) &\le \tilde{I}_b(u_n), \label{a:36} \\[0.3em]
			\|w_n - u_n\| &\le 1/n, \label{a:37} \\[0.3em]
			\tilde{I}_b(w_n) &\le \tilde{I}_b(u)
			+ \frac{1}{n}\|u - w_n\|,
			\quad \forall\, u \in \overline{\mathcal{U}_1}. \label{a:38}
		\end{align}
		We show that $\{w_n\}$ is a $(P.S.)$ sequence for $\tilde{I}_b$ on $M$. The sequence $\tilde{I}_b(w_n)$ is bounded by \eqref{a:36}, and once $\|\tilde{I}'_b(w_n)\|_{*} \to 0$ is known, Definition~\ref{ad:4} applies and Lemma~\ref{al:3} yields a subsequence $w_{n_k} \to w$ strongly in $X(\Gamma,V_D)$. To verify $\|\tilde{I}'_b(w_n)\|_{*} \to 0$, fix $n$ and let $\phi \in X(\Gamma,V_D)$ be tangent to $M$ at $w_n$, that is,
		\begin{align*}
			\int_{\Gamma} |w_n|^{p-2}w_n\,\phi\,dx = 0,
		\end{align*}
		and define
		\begin{align*}
			u_t = \frac{w_n + t\phi}{\|w_n + t\phi\|_p}.
		\end{align*}
		From \eqref{a:36} and the minimizing property, for $n$ sufficiently large,
		\begin{equation}{\label{a5:38}}
			\tilde{I}_b(w_n) \le \tilde{I}_b(u_n) \le h + \frac{1}{2n^2} < r.
		\end{equation}
		Thus $w_n \in \mathcal{U}_1$: if $w_n \in \partial \mathcal{U}_1$, then Lemma~\ref{al:5} would force $\tilde{I}_b(w_n) = r$, contradicting \eqref{a5:38}. We may therefore substitute $u_t$ for $u$ in \eqref{a:38} for $|t|$ sufficiently small, giving
		\begin{align*}
			\tilde{I}_b(w_n) \le \tilde{I}_b(u_t)
			+ \frac{1}{n}\|u_t - w_n\|.
		\end{align*}
		The calculations in the proof of Lemma~\ref{al:4} now show that $\{w_n\}$ is a $(P.S.)$ sequence.
		By \eqref{a:37}, $\{u_n\}$ then admits a convergent subsequence, the infimum $h$ is attained, and the result follows. \qed
	\end{proof}

	\textbf{Proof of Theorem~\ref{at:1}\textit{(i)}}.
	We first show that $(b+\omega,\omega)$ is not a nontrivial point of $\Sigma_{p,\Gamma}$ when $\omega\le\lambda_{1,p}(\Gamma,V_D)$. Let $u$ be a nontrivial solution of \eqref{a1:1} with $(\lambda^+,\lambda^-)=(b+\omega,\omega)$. If $u\ge0$, then $u$ is a nonnegative eigenfunction with eigenvalue $b+\omega$, and Lemma~\ref{al:100}\textit{(iii)} gives $b+\omega=\lambda_{1,p}(\Gamma,V_D)$. If $u\le0$, then in the same way $\omega=\lambda_{1,p}(\Gamma,V_D)$. If $u$ changes sign, testing \eqref{a2:7} with $u^-$ gives $\int_\Gamma|(u^-)'|^p\,dx=\omega\int_\Gamma(u^-)^p\,dx$, hence $\omega\ge\lambda_{1,p}(\Gamma,V_D)$. In all cases the point lies on a trivial line.

	Suppose for contradiction that there exists $\omega$ with
	$\lambda_{1,p}(\Gamma,V_D) < \omega < \alpha(b)$ such that
	$(b+\omega, \omega) \in \Sigma_{p,\Gamma}$. We may choose such an $\omega$ to be minimal, because the lines
	$\mathbb{R} \times \{\lambda_{1,p}(\Gamma, V_D)\}$ and
	$\{\lambda_{1,p}(\Gamma, V_D)\} \times \mathbb{R}$ are isolated in $\Sigma_{p,\Gamma}$ and $\Sigma_{p,\Gamma}$ is closed (Proposition~\ref{ap:4} and Lemma~\ref{al:101}). Thus $\omega$ is a critical value of $\tilde{I}_b$ with $\lambda_{1,p}(\Gamma,V_D) < \omega < \alpha(b)$, and there is no critical value in the interval $(\lambda_{1,p}(\Gamma,V_D), \omega)$. A path from $-\phi_{1,p}$ to $\phi_{1,p}$ along which $\tilde{I}_b \le \omega$ would contradict the definition of $\alpha(b)$
	(see \eqref{a:28}).

	Let $u \in M$ be a critical point of $\tilde{I}_b$ at level $\omega$, so that $u$ satisfies
	\begin{equation}{\label{a:39}}
		\begin{cases}
			-\bigl(|u'|^{p-2}u'\bigr)'
			= (b+\omega)(u^+)^{p-1} - \omega(u^-)^{p-1}
			&\text{on the edges of } \Gamma, \\[6pt]
			u(v) = 0 & \forall\, v \in V_D, \\[8pt]
			\displaystyle\sum_{e \in E_v(\Gamma)} |u_e'(v)|^{p-2}u_e'(v)=0, &v\in V(\Gamma)\setminus V_D.
		\end{cases}
	\end{equation}
	The weak formulation of \eqref{a:39} reads
	\begin{equation}{\label{a:40}}
		\int_{\Gamma} |u'|^{p-2}u'\phi'\,dx
		= (b+\omega)\int_{\Gamma}(u^+)^{p-1}\phi\,dx
		- \omega\int_\Gamma (u^-)^{p-1}\phi\,dx,
		\quad \forall\, \phi \in X(\Gamma,V_D).
	\end{equation}
	With $\phi = u^+$ and $\phi = u^-$ in turn, \eqref{a:40} becomes
	\begin{equation}{\label{a:41}}
		\int_\Gamma |(u^+)'|^p\,dx = (b+\omega)\int_\Gamma (u^+)^p\,dx,
		\qquad
		\int_\Gamma |(u^-)'|^p\,dx = \omega\int_\Gamma (u^-)^p\,dx,
	\end{equation}
	respectively. Thus
	\begin{align*}
		\tilde{I}_b(u) = \omega,
	\end{align*}
	\begin{align*}
		\tilde{I}_b\!\left(\frac{u^+}{\|u^+\|_p}\right)
		= \frac{\int_\Gamma |(u^+)'|^p\,dx - b\int_\Gamma (u^+)^p\,dx}
		{\|u^+\|_p^p}
		= \frac{\omega\int_\Gamma (u^+)^p\,dx}{\|u^+\|_p^p}
		= \omega,
	\end{align*}
	\begin{align*}
		\tilde{I}_b\!\left(\frac{-u^-}{\|u^-\|_p}\right)
		= \frac{\int_\Gamma |(-u^-)'|^p\,dx}{\|u^-\|_p^p}
		= \frac{\omega\int_\Gamma (u^-)^p\,dx}{\|u^-\|_p^p}
		= \omega,
	\end{align*}
	\begin{align*}
		\tilde{I}_b\!\left(\frac{u^-}{\|u^-\|_p}\right)
		= \frac{\int_\Gamma |(u^-)'|^p\,dx - b\int_\Gamma (u^-)^p\,dx}
		{\|u^-\|_p^p}
		= \frac{(\omega-b)\int_\Gamma (u^-)^p\,dx}{\|u^-\|_p^p}
		= \omega - b.
	\end{align*}
	We have four distinguished points, so we consider three paths in $M$: $\gamma_1(t)$ from $u$ to $\frac{u^+}{\|u^+\|_p}$, $\gamma_2(t)$ from $\frac{u^+}{\|u^+\|_p}$ to $\frac{u^-}{\|u^-\|_p}$, and $\gamma_3(t)$ from $\frac{-u^-}{\|u^-\|_p}$ to $u$, defined by
	\begin{align*}
		\gamma_1(t) = \frac{tu + (1-t)u^+}{\|tu + (1-t)u^+\|_p}, \quad
		\gamma_2(t) = \frac{tu^+ + (1-t)u^-}{\|tu^+ + (1-t)u^-\|_p}, \quad
		\gamma_3(t) = \frac{-tu^- + (1-t)u}{\|-tu^- + (1-t)u\|_p}.
	\end{align*}
	Using $tu + (1-t)u^+ = u^+ - tu^-$ and $-tu^- + (1-t)u = (1-t)u^+ - u^-$, together with \eqref{a:41}, we compute for all $t \in [0,1]$:
	\begin{align*}
		\tilde{I}_b\!\left(\frac{tu+(1-t)u^+}{\|tu+(1-t)u^+\|_p}\right)
		= \frac{\int_\Gamma|(u^+-tu^-)'|^p\,dx
			- b\int_\Gamma(u^+)^p\,dx}
		{\|u^+-tu^-\|_p^p} \\[0.3em]
	= \frac{\omega\int_\Gamma(u^+-tu^-)^p\,dx}
	{\|u^+-tu^-\|_p^p} = \omega.
\end{align*}
\begin{align*}
	\tilde{I}_b\!\left(\frac{tu^++(1-t)u^-}{\|tu^++(1-t)u^-\|_p}\right)
&= \frac{\omega\left(\int_\Gamma(tu^++(1-t)u^-)^p\,dx\right)
	- b(1-t)^p\int_\Gamma(u^-)^p\,dx}
{\|tu^++(1-t)u^-\|_p^p} \\[0.3em]
&= \omega - \frac{b(1-t)^p\int_\Gamma(u^-)^p\,dx}
{\|tu^++(1-t)u^-\|_p^p} \le \omega.
\end{align*}
\begin{align*}
\tilde{I}_b\!\left(\frac{-tu^-+(1-t)u}{\|-tu^-+(1-t)u\|_p}\right)
&= \frac{\int_\Gamma|((1-t)u^+-u^-)'|^p\,dx
	- b(1-t)^p\int_\Gamma(u^+)^p\,dx}
{\|(1-t)u^+-u^-\|_p^p} \\[0.3em]
&= \frac{\omega\int_\Gamma((1-t)u^+-u^-)^p\,dx}
{\|(1-t)u^+-u^-\|_p^p} = \omega.
\end{align*}
Traversing $\gamma_1(t)$ from $u$ to $\frac{u^+}{\|u^+\|_p}$ and then $\gamma_2(t)$ from $\frac{u^+}{\|u^+\|_p}$ to $\frac{u^-}{\|u^-\|_p}$, we arrive at $\frac{u^-}{\|u^-\|_p}$, which lies at level $\omega - b$, while staying at levels at most $\omega$ throughout. Let $\mathcal{U} = \{\psi \in M : \tilde{I}_b(\psi) < \omega - b\}$. Then $\phi_{1,p} \in \mathcal{U}$, because $\tilde{I}_b(\phi_{1,p}) = \lambda_{1,p}(\Gamma,V_D) - b < \omega - b$, and $-\phi_{1,p} \in \mathcal{U}$ whenever $\omega - b > \lambda_{1,p}(\Gamma,V_D)$. By the choice of $\omega$, the critical points of $\tilde{I}_b$ in $\mathcal{U}$ are $\phi_{1,p}$ and, if $\omega-b>\lambda_{1,p}(\Gamma,V_D)$, also $-\phi_{1,p}$. Let $v\in\mathcal U$ be a critical point with value $c<\omega-b\le\omega$. If $v\ge0$, then $v$ is a nonnegative eigenfunction with eigenvalue $b+c$, so $v=\phi_{1,p}$ by Lemma~\ref{al:100}. If $v\le0$, then $-v$ is a nonnegative eigenfunction with eigenvalue $c$, so $v=-\phi_{1,p}$. If $v$ changes sign, testing with $v^-$ gives $c\ge\lambda_{1,p}(\Gamma,V_D)$, and $c=\lambda_{1,p}(\Gamma,V_D)$ would make $v^-/\|v^-\|_p$ a first eigenfunction vanishing on the nonempty open set $\{v>0\}$, contrary to Lemma~\ref{al:100}\textit{(ii)}. This leaves $\lambda_{1,p}(\Gamma,V_D)<c<\omega$, contrary to the choice of $\omega$.

Recall that $\tilde{I}_b\!\left(\frac{u^-}{\|u^-\|_p}\right) = \omega - b$. The function $\frac{u^-}{\|u^-\|_p}$ does not change sign and vanishes on a set of positive measure, hence it is not a critical point of $\tilde{I}_b$ (by Lemma~\ref{al:100}, a critical point of constant sign is $\pm\phi_{1,p}$, which does not vanish on $\Gamma\setminus V_D$). Put $v=\frac{u^-}{\|u^-\|_p}$. Since $\|\tilde I_b'(v)\|_*>0$, there is $\varphi\in T_vM$ with $\langle I_b'(v),\varphi\rangle<0$, and $\xi(t)=(v+t\varphi)/\|v+t\varphi\|_p$ is a $C^1$ path in $M$ with $\xi(0)=v$ and $\frac{d}{dt}\tilde I_b(\xi(t))\big|_{t=0}=\langle I_b'(v),\varphi\rangle<0$, because the derivative of $t\mapsto\|v+t\varphi\|_p$ vanishes at $t=0$. Hence, for small $t>0$, $\xi([0,t])$ joins $v$ to a point $\psi=\xi(t)$ with $\tilde{I}_b(\psi) < \omega - b$, and $\tilde I_b\le\omega-b$ along $\xi([0,t])$. By Lemma~\ref{al:6}, the connected component of $\mathcal{U}$ containing $\psi$ contains $\phi_{1,p}$ or $-\phi_{1,p}$, and only the first case can occur if $\omega-b\le\lambda_{1,p}(\Gamma,V_D)$, because then $-\phi_{1,p}\notin\mathcal U$. Suppose it contains $\phi_{1,p}$. Then there exists a path $\gamma_4(t)$ from
$\frac{u^-}{\|u^-\|_p}$ to $\phi_{1,p}$ along which $\tilde{I}_b\le\omega-b$. Concatenating $\gamma_1(t)$, $\gamma_2(t)$, and $\gamma_4(t)$ yields a path in $M$ from $u$ to $\phi_{1,p}$ at levels at most $\omega$.

Consider the symmetric path $-\gamma_4(t)$, going from
$\frac{-u^-}{\|u^-\|_p}$ to $-\phi_{1,p}$. For any $u \in M$, $\left|\tilde{I}_b(u) - \tilde{I}_b(-u)\right|
\le b$. Therefore
\begin{align*}
\tilde{I}_b(-\gamma_4(t)) \le \tilde{I}_b(\gamma_4(t)) + b
\le \omega - b + b = \omega, \quad \forall\, t \in [0,1].
\end{align*}
Joining $\gamma_3(t)$ and $-\gamma_4(t)$ then gives a path from $u$ to $-\phi_{1,p}$. Altogether, this is a path $\gamma(t)$ in $M$ from $-\phi_{1,p}$ to $\phi_{1,p}$ with $\tilde{I}_b(\gamma(t)) \le \omega$ for all $t \in [0,1]$, the desired contradiction.

If instead the component contains $-\phi_{1,p}$, let $\gamma_4$ be a path from $\frac{u^-}{\|u^-\|_p}$ to $-\phi_{1,p}$ along which $\tilde I_b\le\omega-b$. Then $-\gamma_4$ joins $\frac{-u^-}{\|u^-\|_p}$ to $\phi_{1,p}$ with $\tilde I_b\le\omega$, and traversing $\gamma_4$, $\gamma_2$, $\gamma_1$, $\gamma_3$ backwards and then $-\gamma_4$ gives a path from $-\phi_{1,p}$ to $\phi_{1,p}$ along which $\tilde I_b\le\omega$, the same contradiction. \qed

The following result, which provides a variational characterization of the second eigenvalue $\lambda_{2,p}(\Gamma, V_D)$ of problem \eqref{a1:1}, follows directly from Theorem \ref{at:1}.
\begin{corollary}{\label{ar:2}}
The infimum in \eqref{a:104} is attained, and one has
\begin{align*}
\lambda_{2,p}(\Gamma,V_D) = \alpha(0)=\inf _{P \in \mathbb{P}} \max_{u \in P[-1,+1]} \int_\Gamma |u'|^p \, dx.
\end{align*}
where $\mathbb{P}$ denotes the family of all continuous paths joining $-\phi_{1,p}$ and $\phi_{1,p}$.
\end{corollary}

\begin{proof}
By Proposition~\ref{at4:5} with $b=0$, $(\alpha(0),\alpha(0))\in\Sigma_{p,\Gamma}$, so $\alpha(0)$ is an eigenvalue with $\alpha(0)>\lambda_{1,p}(\Gamma,V_D)$ and $\lambda_{2,p}(\Gamma,V_D)\le\alpha(0)$. By Theorem~\ref{at:1} with $b=0$, $(\omega,\omega)\notin\Sigma_{p,\Gamma}$ for $\lambda_{1,p}(\Gamma,V_D)<\omega<\alpha(0)$, that is, there is no eigenvalue in $(\lambda_{1,p}(\Gamma,V_D),\alpha(0))$. Hence $\lambda_{2,p}(\Gamma,V_D)=\alpha(0)$, and the infimum in \eqref{a:104} is attained. Since $\tilde I_0(u)=\int_\Gamma|u'|^p\,dx$ on $M$,
\begin{align*}
\alpha(0) =\inf_{P \in \mathbb{P}}\, \max_{u \in P[-1,+1]} I_0(u) =\inf_{P \in \mathbb{P}}\, \max_{u \in P[-1,+1]}
\int_{\Gamma} |u'|^p\,dx.
\end{align*}
 \qed
\end{proof}

\section{Properties of the first nontrivial curve}{\label{as:6}}

In this section, we establish that the first nontrivial curve $\mathfrak{C}$ is Lipschitz continuous, strictly decreasing, and exhibits specific asymptotic behavior on the compact connected quantum graph.

\begin{prop}{\label{ap:5}}
The first nontrivial curve, described by the map
$\mathfrak{C} : \mathbb{R}^+ \to \mathbb{R}^2$ given by
$\mathfrak{C}(b) = (b + \alpha(b),\, \alpha(b))$, is Lipschitz continuous.
\end{prop}

\begin{proof}
Let $0\le b<b'$. On $M$ we have $0\le\tilde I_b-\tilde I_{b'}=(b'-b)\int_\Gamma(u^+)^p\,dx\le b'-b$. Hence, for every $P\in\mathbb P$, $\max_P\tilde I_{b'}\le\max_P\tilde I_b\le\max_P\tilde I_{b'}+b'-b$, and \eqref{a:28} gives $0\le\alpha(b)-\alpha(b')\le b'-b$. Thus $\alpha$ is nonincreasing and $1$-Lipschitz, and $\|\mathfrak C(b)-\mathfrak C(b')\|\le|b-b'|+2|\alpha(b)-\alpha(b')|\le3|b-b'|$. \qed
\end{proof}
For an open set $\mathcal U\subseteq\Gamma$ we put
\begin{equation}\label{a100:100}
X_{\mathcal U}=\{u\in X(\Gamma,V_D):u=0\text{ on }\Gamma\setminus\mathcal U\},\qquad
\lambda_{1,p}(\mathcal U)=\inf\Bigl\{\frac{\int_\Gamma|u'|^p\,dx}{\int_\Gamma|u|^p\,dx}:u\in X_{\mathcal U},\ u\not\equiv0\Bigr\}.
\end{equation}
The infimum is attained by Theorem~\ref{at:3} applied to the closed subspace $X_{\mathcal U}$ of $X(\Gamma,V_D)$, and a minimizer may be taken nonnegative.
\begin{lemma}{\label{al:7}}
Let $\mathcal{U}_1$, $\mathcal{U}_2$ be two open subsets of $\Gamma$ with $\mathcal{U}_1 \subset \mathcal{U}_2$,  $\mathcal{U}_2 \setminus V_D$ connected and $\ell(\mathcal U_2\setminus \mathcal U_1)>0$. Then $\lambda_{1,p}(\mathcal{U}_1) > \lambda_{1,p}(\mathcal{U}_2)$.
\end{lemma}

\begin{proof}
The inequality $\lambda_{1,p}(\mathcal{U}_1) \ge \lambda_{1,p}(\mathcal{U}_2)$ holds by the definition of the first eigenvalue $\lambda_{1,p}$. Suppose for contradiction that equality holds, that is, $\lambda_{1,p}(\mathcal{U}_1) = \lambda_{1,p}(\mathcal{U}_2)$. Let $\phi_{1,p}$ denote the extension by zero outside $\mathcal{U}_1$ of a nonnegative normalized minimizer for $\lambda_{1,p}(\mathcal{U}_1)$. Then
\begin{align*}
\int_{\mathcal{U}_2} |\phi_{1,p}'|^p\,dx
= \lambda_{1,p}(\mathcal{U}_1)\int_{\mathcal{U}_2}|\phi_{1,p}|^p\,dx
= \lambda_{1,p}(\mathcal{U}_2)\int_{\mathcal{U}_2}|\phi_{1,p}|^p\,dx,
\end{align*}
which means that $\phi_{1,p}$ attains $\lambda_{1,p}(\mathcal U_2)$ and satisfies $\int_\Gamma|\phi_{1,p}'|^{p-2}\phi_{1,p}'\phi'\,dx=\lambda_{1,p}(\mathcal U_2)\int_\Gamma\phi_{1,p}^{p-1}\phi\,dx$ for all $\phi\in X_{\mathcal U_2}$. The proof of Lemma~\ref{al:100}\textit{(i)}, applied on $\mathcal U_2\setminus V_D$ (the Kirchhoff condition holds at the vertices of $\Gamma$ lying in $\mathcal U_2\setminus V_D$), gives $\phi_{1,p}>0$ on $\mathcal U_2\setminus V_D$. This contradicts $\phi_{1,p}=0$ on $\mathcal U_2\setminus\mathcal U_1$, a set of positive length. Therefore, the inequality is strict. \qed
\end{proof}

\begin{lemma}{\label{al:9}}
Let $(\lambda^+, \lambda^-) \in \mathfrak{C}$ and let $\lambda^+(x)$ and $\lambda^-(x)$ be $L^\infty$ functions satisfying
\begin{equation}{\label{a:43}}
\lambda_{1,p}(\Gamma,V_D) \le \lambda^+(x) \le \lambda^+,
\qquad \lambda_{1,p}(\Gamma,V_D) \le \lambda^-(x) \le \lambda^- \quad \text{a.e.\ in } \Gamma.
\end{equation}
Assume further that
\begin{equation}{\label{a:44}}
\lambda_{1,p}(\Gamma,V_D) < \lambda^+(x)
\quad\text{and}\quad
\lambda_{1,p}(\Gamma,V_D) < \lambda^-(x)
\quad\text{on subsets of positive measure.}
\end{equation}
Then any nontrivial solution $u$ of
\begin{equation}{\label{a:45}}
\begin{cases}
	-\bigl(|u'|^{p-2}u'\bigr)'
	= \lambda^+(x)(u^+)^{p-1} - \lambda^-(x)(u^-)^{p-1}
	&\text{on the edges of } \Gamma, \\[6pt]
	u(v) = 0 & \forall\, v \in V_D, \\[8pt]
	\displaystyle\sum_{e \in E_v(\Gamma)} |u_e'(v)|^{p-2}u_e'(v)=0, &v\in V(\Gamma)\setminus V_D.
\end{cases}
\end{equation}
changes sign in $\Gamma$, and moreover
\begin{align*}
\lambda^+(x) = \lambda^+ \;\text{ a.e.\ on }
\{x \in \Gamma : u(x) > 0\}
\quad\text{and}\quad
\lambda^-(x) = \lambda^- \;\text{ a.e.\ on }
\{x \in \Gamma : u(x) < 0\}.
\end{align*}
\end{lemma}

\begin{proof}
Let $u$ be a nontrivial solution of \eqref{a:45}. Replacing $u$ by $-u$ if necessary, we may assume that the point
$(\lambda^+, \lambda^-) \in \mathfrak{C}$ satisfies $\lambda^+ \ge \lambda^-$. We first show that $u$ changes sign in $\Gamma$. Suppose, for contradiction, that it does not. There are two cases:
\begin{enumerate}[label=(\roman*)]
\item $u \ge 0$ a.e.\ on $\Gamma$,
\item $u \le 0$ a.e.\ on $\Gamma$.
\end{enumerate}
In case (i), $u$ solves
\begin{align*}
\begin{cases}
	-\bigl(|u'|^{p-2}u'\bigr)' = \lambda^+(x)|u|^{p-2}u
	&\text{on the edges of } \Gamma, \\[6pt]
	u(v) = 0 & \forall\, v \in V_D, \\[8pt]
	\displaystyle\sum_{e \in E_v(\Gamma)} |u_e'(v)|^{p-2}u_e'(v)=0, &v\in V(\Gamma)\setminus V_D.
\end{cases}
\end{align*}
By Lemma~\ref{al:100}\textit{(iii)} with $m=\lambda^+(x)$, $\int_\Gamma(\lambda^+(x)-\lambda_{1,p}(\Gamma,V_D))\phi_{1,p}^p\,dx\le0$. This is impossible: $\lambda^+(x)\ge\lambda_{1,p}(\Gamma,V_D)$ a.e., $\lambda^+(x)>\lambda_{1,p}(\Gamma,V_D)$ on a set of positive measure by \eqref{a:44}, and $\phi_{1,p}>0$ a.e. Case (ii) is similar. Therefore $u$ changes sign. Suppose now for contradiction that either
\begin{equation}{\label{a:48}}
\ell\{x \in \Gamma : \lambda^+(x) < \lambda^+,\; u(x) > 0\} > 0
\end{equation}
or
\begin{equation}{\label{a:49}}
\ell\{x \in \Gamma : \lambda^-(x) < \lambda^-,\; u(x) < 0\} > 0.
\end{equation}
Assume \eqref{a:48} holds (the argument for \eqref{a:49} is analogous). Set $b = \lambda^+ - \lambda^- \ge 0$, so that $\lambda^- = \alpha(b)$ as defined in \eqref{a:28}. We construct a path $P \in \mathbb{P}$ such that
\begin{equation}{\label{a:50}}
\max_{u \in P[-1,+1]} \tilde{I}_b(u) < \lambda^-,
\end{equation}
which contradicts the definition of $\alpha(b)$ as a minimax value. The path is built from a sign-changing function $w \in X(\Gamma, V_D)$ satisfying
\begin{equation}{\label{a:51}}
\frac{\int_\Gamma|(w^+)'|^p\,dx}{\int_\Gamma(w^+)^p\,dx} < \lambda^+
\quad\text{and}\quad
\frac{\int_\Gamma|(w^-)'|^p\,dx}{\int_\Gamma(w^-)^p\,dx} < \lambda^-.
\end{equation}
To find $w$, let $\mathcal{U}_1$ be a connected component of $\{x \in \Gamma : u(x) > 0\}$ with
\begin{align*}
\ell\{x \in \mathcal{U}_1 : \lambda^+(x) < \lambda^+\} > 0,
\end{align*}
which is possible by \eqref{a:48}. Let $\mathcal{U}_2$ be any connected component of $\{x \in \Gamma : u(x) < 0\}$, which is nonempty since $u$ changes sign.
 We claim that
\begin{equation}{\label{a:52}}
\lambda_{1,p}(\mathcal{U}_1) < \lambda^+
\quad\text{and}\quad
\lambda_{1,p}(\mathcal{U}_2) \le \lambda^-,
\end{equation}
where $\lambda_{1,p}(\mathcal{U}_i)$ is defined by~\eqref{a100:100}.
The function $u\chi_{\mathcal U_i}$ belongs to $X_{\mathcal U_i}$, as $u=0$ on $\partial\mathcal U_i$. Testing \eqref{a:45} with $u\chi_{\mathcal U_i}$ and using \eqref{a:43}, we obtain

\begin{align*}
\int_{\mathcal U_1}|u'|^p\,dx=\int_{\mathcal U_1}\lambda^+(x)u^p\,dx<\lambda^+\int_{\mathcal U_1}u^p\,dx ,
\end{align*}
the inequality being strict because $\lambda^+(x)<\lambda^+$ on a subset of $\mathcal U_1$ of positive measure, where $u>0$,
which gives $\lambda_{1,p}(\mathcal{U}_1) < \lambda^+$. The second inequality in \eqref{a:52} follows similarly, with $\le$ since only $\lambda^-(x)\le\lambda^-$ is known on $\mathcal U_2$. We now enlarge $\mathcal{U}_1$ and $\mathcal{U}_2$ slightly to two disjoint open sets $\tilde{\mathcal{U}}_1$ and
$\tilde{\mathcal{U}}_2$ satisfying
\begin{align*}
\lambda_{1,p}(\tilde{\mathcal{U}}_1) < \lambda^+
\quad\text{and}\quad
\lambda_{1,p}(\tilde{\mathcal{U}}_2) < \lambda^-.
\end{align*}
For $\xi \ge 0$, define $\mathcal{U}_1(\xi) = \{x \in \mathcal{U}_1 : d(x,\mathcal{U}_1^c) > \xi\}$. Then $\lambda_{1,p}(\mathcal{U}_1(\xi)) \ge \lambda_{1,p}(\mathcal{U}_1)$. Conversely, let $\psi\ge0$ be a minimizer for $\lambda_{1,p}(\mathcal U_1)$ and, for $\delta>0$, let $\psi_\delta=(\psi-\delta)^+$. Since $\psi$ is uniformly continuous on $\Gamma$ and vanishes on $\Gamma\setminus\mathcal U_1$, there is $\xi_\delta>0$ such that $\psi<\delta$ on $\{x:d(x,\mathcal U_1^c)\le\xi_\delta\}$, so $\psi_\delta\in X_{\mathcal U_1(\xi)}$ for $\xi\le\xi_\delta$, and $\psi_\delta\to\psi$ in $W^{1,p}(\Gamma)$ as $\delta\to0$. Hence $\lambda_{1,p}(\mathcal{U}_1(\xi)) \to \lambda_{1,p}(\mathcal{U}_1)$ as $\xi \to 0$, and there exists $\xi_0 > 0$ such that
\begin{equation}{\label{a:54}}
\lambda_{1,p}(\mathcal{U}_1(\xi)) < \lambda^+,
\quad \text{for all } 0 \le \xi \le \xi_0.
\end{equation}
We claim that there is a point $x_0\in\partial\mathcal{U}_2 \setminus V_D$ such that $B(x_0,r)\setminus\overline{\mathcal{U}_2}$ has positive length for every $r>0$. Otherwise every $x\in\partial\mathcal{U}_2 \setminus V_D$ has a ball $B(x,r_x)$ with $\ell(B(x,r_x)\setminus\overline{\mathcal{U}_2})=0$. The open set $B(x,r_x)\setminus\overline{\mathcal{U}_2}$ is then empty, so $\overline{\mathcal{U}_2}\setminus V_D$ is open as well as closed in $\Gamma\setminus V_D$. Then $\overline{\mathcal{U}_2} \setminus V_D=\Gamma\setminus V_D$ by \eqref{ah:100}, and $\overline{\mathcal U_2}=\Gamma$ since $V_D$ is finite. But $\mathcal{U}_1$ is open, disjoint from $\mathcal{U}_2$, and therefore disjoint from $\overline{\mathcal{U}_2}$, contradicting $\mathcal{U}_1\neq\emptyset$. Fix such an $x_0$ and choose $0<\xi<\xi_0$ so small that $B(x_0,\xi/2)\cap V_D=\emptyset$. Balls in $\Gamma$ are connected: every point of $B(x_0,\xi/2)$ is joined to $x_0$ by a shortest path. Also $x_0\notin\mathcal{U}_1$, because $\mathcal{U}_1$ is open and disjoint from $\mathcal{U}_2$, so $d(x_0,\mathcal{U}_1(\xi))\ge\xi$. Setting
\begin{align*}
\tilde{\mathcal{U}}_1 = \mathcal{U}_1(\xi),
\qquad
\tilde{\mathcal{U}}_2 = \mathcal{U}_2 \cup B(x_0, \xi/2),
\end{align*}
we have $\tilde{\mathcal{U}}_1 \cap \tilde{\mathcal{U}}_2 = \emptyset$ because $d(x_0,\mathcal{U}_1(\xi))\ge\xi>\xi/2$, and $\lambda_{1,p}(\tilde{\mathcal{U}}_1) < \lambda^+$ by \eqref{a:54}. The set $\tilde{\mathcal{U}}_2$ is open, connected, disjoint from $V_D$ (recall $u=0$ on $V_D$), and $\ell(\tilde{\mathcal{U}}_2\setminus\mathcal{U}_2)>0$ by the choice of $x_0$, so Lemma~\ref{al:7} together with \eqref{a:52} gives $\lambda_{1,p}(\tilde{\mathcal{U}}_2) < \lambda^-$.

Define $w = w_1 - w_2$, where $w_i$ is the extension by zero outside $\tilde{\mathcal{U}}_i$ of the eigenfunction associated with $\lambda_{1,p}(\tilde{\mathcal{U}}_i)$. Then $w$ satisfies \eqref{a:51}. For this $w$,
\begin{align*}
\tilde{I}_b\!\left(\frac{w}{\|w\|_p}\right)
&< (\lambda^+-b)
\frac{\int_\Gamma(w^+)^p\,dx}{\|w\|_p^p}
+ \lambda^-
\frac{\int_\Gamma(w^-)^p\,dx}{\|w\|_p^p}
= \lambda^-,
\end{align*}
and
\begin{align*}
\tilde{I}_b\!\left(\frac{w^+}{\|w^+\|_p}\right)
&= \frac{\int_\Gamma|(w^+)'|^p\,dx}
{\|w^+\|_p^p}
- b\frac{\int_\Gamma(w^+)^p\,dx}
{\|w^+\|_p^p}
< \lambda^+-b=\lambda^-.
\end{align*}
Similarly,
$\tilde{I}_b(w^-/\|w^-\|_p)<\lambda^--b$.
Put $\nu=\lambda^-=\lambda^+-b$. The functions
$w^+=w_1$ and $w^-=w_2$ satisfy
\[
\int_\Gamma |(w^+)'|^p\,dx
<\lambda^+\int_\Gamma(w^+)^p\,dx,
\qquad
\int_\Gamma |(w^-)'|^p\,dx
<\lambda^-\int_\Gamma(w^-)^p\,dx.
\]
Put $a_\pm=\int_\Gamma(w^\pm)^p\,dx$. For $t\in[0,1]$, these inequalities and $\lambda^+-b=\lambda^-=\nu$ give
\begin{align*}
\tilde I_b\Bigl(\frac{w^+-tw^-}{\|w^+-tw^-\|_p}\Bigr)&<\frac{(\lambda^+-b)a_++t^p\lambda^-a_-}{a_++t^pa_-}=\nu,\\
\tilde I_b\Bigl(\frac{tw^++(1-t)w^-}{\|tw^++(1-t)w^-\|_p}\Bigr)&<\frac{t^p(\lambda^+-b)a_++(1-t)^p(\lambda^--b)a_-}{t^pa_++(1-t)^pa_-}\le\nu,\\
\tilde I_b\Bigl(\frac{(1-t)w^+-w^-}{\|(1-t)w^+-w^-\|_p}\Bigr)&<\frac{(1-t)^p(\lambda^+-b)a_++\lambda^-a_-}{(1-t)^pa_++a_-}=\nu .
\end{align*}
These three paths join
$\frac{w}{\|w\|_p}$ to $\frac{w^+}{\|w^+\|_p}$,
$\frac{w^+}{\|w^+\|_p}$ to $\frac{w^-}{\|w^-\|_p}$, and
$\frac{-w^-}{\|w^-\|_p}$ to $\frac{w}{\|w\|_p}$, respectively,
and $\tilde I_b<\nu$ along each of them.
Let $\mathcal O=\{\tilde I_b<\nu-b\}$.
Since $\nu=\alpha(b)$, the argument in the
proof of Theorem~\ref{at:1} shows that the
only critical points of $\tilde I_b$ in
$\mathcal O$ are $\phi_{1,p}$ and possibly
$-\phi_{1,p}$, so by Lemma~\ref{al:6}
the component of $\mathcal O$ containing
$\frac{w^-}{\|w^-\|_p}$ contains
$\phi_{1,p}$ or $-\phi_{1,p}$.
Let $\gamma_4$ be a path in this component
from $\frac{w^-}{\|w^-\|_p}$ to that point.
Since $\tilde I_b(-v)\le\tilde I_b(v)+b$
on $M$, we have $\tilde I_b<\nu$
along $-\gamma_4$. If $\gamma_4$ ends at $\phi_{1,p}$, we go from $-\phi_{1,p}$ to
$\frac{-w^-}{\|w^-\|_p}$ along $-\gamma_4$,
then to $\frac{w}{\|w\|_p}$,
$\frac{w^+}{\|w^+\|_p}$,
$\frac{w^-}{\|w^-\|_p}$ along the three
paths above, and to $\phi_{1,p}$
along $\gamma_4$. If $\gamma_4$ ends at $-\phi_{1,p}$, we go from $-\phi_{1,p}$ to
$\frac{w^-}{\|w^-\|_p}$ along $\gamma_4$,
then to $\frac{w^+}{\|w^+\|_p}$,
$\frac{w}{\|w\|_p}$,
$\frac{-w^-}{\|w^-\|_p}$, and to
$\phi_{1,p}$ along $-\gamma_4$. In both cases we obtain $P\in\mathbb P$
along which $\tilde I_b<\nu$. As
$P[-1,1]$ is compact,
\eqref{a:50} holds, which is the
desired contradiction.
\qed
\end{proof}

\begin{corollary}{\label{ac:1}}
Let $(\lambda^+, \lambda^-) \in \mathfrak{C}$ and let $\lambda^+(x)$, $\lambda^-(x)$ be $L^\infty$ functions on $\Gamma$ satisfying \eqref{a:43} and \eqref{a:44}. If either $\lambda^+(x) < \lambda^+$ a.e.\ in $\Gamma$ or $\lambda^-(x) < \lambda^-$ a.e.\ in $\Gamma$, then \eqref{a:45} has only the trivial solution.
\end{corollary}

\begin{proof}
This follows from Lemma~\ref{al:9}.  \qed
\end{proof}

\begin{prop}{\label{ap:6}}
The first nontrivial curve described by the map
$\mathfrak{C} : \mathbb{R}^+ \to \mathbb{R}^2$ such that
\begin{align*}
\mathfrak{C}(b) = (b + \alpha(b),\, \alpha(b))
\end{align*}
is strictly decreasing, meaning $b < b'$ implies
$b + \alpha(b) < b' + \alpha(b')$ and $\alpha(b) > \alpha(b')$.
\end{prop}

\begin{proof}
Let $b<b'$. If $b+\alpha(b)\ge b'+\alpha(b')$, then $\alpha(b)>\alpha(b')$, and Corollary~\ref{ac:1}, applied at $(\lambda^+,\lambda^-)=(b+\alpha(b),\alpha(b))$ with the constants $\lambda^+(x)=b'+\alpha(b')$ and $\lambda^-(x)=\alpha(b')$, shows that \eqref{a1:1} has no nontrivial solution at $(b'+\alpha(b'),\alpha(b'))$. This contradicts Proposition~\ref{at4:5}. If $\alpha(b)\le\alpha(b')$, then $b+\alpha(b)<b'+\alpha(b')$, and the same argument at $(b'+\alpha(b'),\alpha(b'))$ gives a contradiction. \qed
\end{proof}

The function $\alpha(b)$ is continuous, strictly decreasing and $\alpha(b)>\lambda_{1,p}(\Gamma,V_D)$, so its limit as $b \to \infty$ exists. The following proposition determines this limit.

\begin{prop}{\label{ap:7}}
Let $V_D\neq\emptyset$. The limit of $\alpha(b)$ as $b \to \infty$ equals $\lambda_{1,p}(\Gamma,V_D)$.
\end{prop}

\begin{proof}
Suppose for contradiction that there exists $\epsilon > 0$ such that
\begin{align*}
\max_{u \in P[-1,+1]} \tilde{I}_b(u)
\ge \lambda_{1,p}(\Gamma,V_D) + \epsilon,
\quad \text{for all } P \in \mathbb{P} \text{ and all } b \ge 0.
\end{align*}
Choose $\phi \in X(\Gamma, V_D)$ such that there is no $r \in \mathbb{R}$ satisfying $\phi(x) \le r\phi_{1,p}$ a.e.\ in $\Gamma$. Such a $\phi$ exists because $V_D\neq\emptyset$: fix $v\in V_D$, an edge $e$ incident to $v$, a number $\theta\in(1-\tfrac1p,1)$, and let $\phi\ge0$ be a function in $X(\Gamma,V_D)$ which equals $d(x,v)^{\theta}$ on a neighbourhood of $v$ in $e$ and vanishes outside $e$. Then $\phi\in W^{1,p}(\Gamma)$ since $p(\theta-1)>-1$, while $\phi_{1,p}$ is $C^1$ on $e$ with $\phi_{1,p}(v)=0$, so $\phi_{1,p}(x)\le C\,d(x,v)$ near $v$ and $\phi/\phi_{1,p}\ge d(x,v)^{\theta-1}/C$ is unbounded near $v$. Consider the path $P \in \mathbb{P}$ defined by
\begin{align*}
P(t) = \frac{t\phi_{1,p} + (1-|t|)\phi}
{\|t\phi_{1,p} + (1-|t|)\phi\|_p},
\quad t \in [-1,+1].
\end{align*}
The set $P[-1,+1]$ is compact, so the maximum of $\tilde{I}_b$ on $P[-1,+1]$ is attained at some $t_b \in [-1,+1]$. Setting
$w_{t_b} = t_b\phi_{1,p} + (1-|t_b|)\phi$, we obtain
\begin{equation}{\label{a:75}}
\int_\Gamma|w_{t_b}'|^p\,dx
- b\int_\Gamma(w_{t_b}^+)^p\,dx
\ge (\lambda_{1,p}(\Gamma,V_D) + \epsilon)
\int_\Gamma|t_b\phi_{1,p} + (1-|t_b|)\phi|^p\,dx
\end{equation}
for all $b \ge 0$. As $b \to +\infty$, passing to a subsequence we may assume $t_b \to \bar{t} \in [-1,+1]$. Since $w_{t_b}$ is bounded in $X(\Gamma,V_D)$ independently of $t$, letting $b \to +\infty$ forces $\int_\Gamma(w_{t_b}^+)^p\,dx \to 0$. As $w_{t_b} \to w_{\bar{t}}$ in $X(\Gamma,V_D)$, we get
\begin{align*}
\int_\Gamma
\bigl((\bar{t}\phi_{1,p} + (1-|\bar{t}|)\phi)^+\bigr)^p\,dx = 0,
\end{align*}
which, by the choice of $\phi$, is only possible if $\bar{t} = -1$. Thus $t_b \to -1$, and passing to the limit in \eqref{a:75} yields
\begin{align*}
\int_\Gamma|\phi_{1,p}'|^p\,dx
&\ge (\lambda_{1,p}(\Gamma,V_D) + \epsilon)
\int_\Gamma|\phi_{1,p}|^p\,dx, \\[0.5em]
\lambda_{1,p}(\Gamma,V_D)\int_\Gamma|\phi_{1,p}|^p\,dx
&\ge (\lambda_{1,p}(\Gamma,V_D) + \epsilon)
\int_\Gamma|\phi_{1,p}|^p\,dx,
\end{align*}
a contradiction. \qed
\end{proof}

\section{Neumann-Kirchhoff vertex condition}{\label{as:8}}
In this section, we study the eigenvalue problem arising when the $p$-Laplacian $\bigl(|u'|^{p-2}u'\bigr)'$ serves as the differential operator on each edge of the graph, subject to purely Neumann-Kirchhoff vertex conditions at all vertices $v \in V(\Gamma)$ \cite{kostrykin1999kirchhoff}. Specifically, we consider the problem
\begin{equation}{\label{a:76}}
\begin{cases}
-\bigl(|u'|^{p-2}u'\bigr)' = \lambda|u|^{p-2}u
&\text{on the edges of } \Gamma, \\[8pt]
\displaystyle\sum_{e \in E_v(\Gamma)} |u_e'(v)|^{p-2}u_e'(v)=0,
&\text{for all vertices } v \in V(\Gamma).
\end{cases}
\end{equation}

Let $\Gamma$ be a compact connected quantum graph and $p \in (1,\infty)$. A real number $\lambda \in \mathbb{R}$ is called a Neumann eigenvalue of the $p$-Laplacian on $\Gamma$ if there exists a nontrivial function $u \in W^{1,p}(\Gamma)$ such that
\begin{align*}
\int_{\Gamma} |u'(x)|^{p-2}u'(x)\,\phi'(x)\,dx
= \lambda\int_{\Gamma} |u(x)|^{p-2}u(x)\,\phi(x)\,dx
\end{align*}
for all $\phi \in W^{1,p}(\Gamma)$. Then $u$ is called an eigenfunction corresponding to $\lambda$. The first eigenvalue $\lambda_{1,p}(\Gamma,\emptyset) = 0$ of \eqref{a:76} is simple, with eigenspace $\mathbb{R}$ (see \cite{del2016clustering}). Hence all eigenfunctions corresponding to $\lambda_{1,p}(\Gamma,\emptyset)$ have constant sign on $\Gamma$, a property that fails for higher eigenvalues (Lemma~\ref{al:100}\textit{(iii)}). Moreover, $\mathcal{F}_{p,\Gamma}$ contains the points $(0,0)$ and $(\lambda_{2,p}(\Gamma,\emptyset), \lambda_{2,p}(\Gamma,\emptyset))$ (where $\lambda_{2,p}(\Gamma,\emptyset)$ denotes the second eigenvalue of \eqref{a:76}), as well as the two lines $\{0\} \times \mathbb{R}$ and $\mathbb{R} \times \{0\}$. The nontrivial part of $\mathcal{F}_{p,\Gamma}$ is denoted by $\widetilde{\mathcal{F}}_{p,\Gamma}$ and defined as
\begin{align*}
\widetilde{\mathcal{F}}_{p,\Gamma}
= \mathcal{F}_{p,\Gamma}
\setminus \bigl((\{0\} \times \mathbb{R})
\cup (\mathbb{R} \times \{0\})\bigr),
\end{align*}
which is contained in $(\bar\lambda,\infty)\times(\bar\lambda,\infty)$ by Proposition~\ref{at:6}, where $\bar\lambda>0$ is defined in \eqref{a12:3}.

For $V_D=\emptyset$ we use the functional $\tilde I_b$, the manifold $M$, the class $\mathbb P$ and the minimax value $\alpha(b)$ of Section~\ref{as:4}, with $X(\Gamma,\emptyset)=W^{1,p}(\Gamma)$ and $\phi_{1,p}=\ell(\Gamma)^{-1/p}$, where $\ell(\Gamma)$ is the total length of $\Gamma$. By Lemma~\ref{all1:1}, the critical values of $\tilde I_b$ correspond to the points of $\mathcal F_{p,\Gamma}$ on the line $(b,0)+t(1,1)$. The proofs of Sections~\ref{as:4}--\ref{as:6} were written for $V_D\neq\emptyset$, but they use the Dirichlet set only in a few places, as the following result shows.

\begin{prop}\label{ap:100}
Let $V_D=\emptyset$, so that $X(\Gamma,\emptyset)=W^{1,p}(\Gamma)$, $\lambda_{1,p}(\Gamma,\emptyset)=0$ and $\phi_{1,p}=\ell(\Gamma)^{-1/p}$. Then Theorem~\ref{at:1}\textit{(i)}--\textit{(ii)} and all results of Sections~\ref{as:4}--\ref{as:6} except Proposition~\ref{ap:7} hold with $\Sigma_{p,\Gamma}$ replaced by $\mathcal F_{p,\Gamma}$.
\end{prop}
\begin{proof}
The proofs use the Dirichlet set only through the following three facts, each of which also holds for $V_D=\emptyset$.
\begin{enumerate}[label=(\alph*)]
	\item Simplicity of $\lambda_{1,p}$, positivity of $\phi_{1,p}$ on $\Gamma\setminus V_D$, and the fact that eigenfunctions of constant sign belong to $\lambda_{1,p}$. These are used in Propositions~\ref{ap:2} and~\ref{ap:4}, in the proof of Theorem~\ref{at:1} and in Lemmas~\ref{al:7} and~\ref{al:9}, always through Lemma~\ref{al:100}, which was proved for every $V_D$ satisfying \eqref{ah:100}, including $V_D=\emptyset$.
	\item Boundedness. In Lemmas~\ref{al:2}, \ref{al:3}, \ref{al:4}, \ref{al:101}, \ref{al:6} and Proposition~\ref{ap:4}, boundedness in $X(\Gamma,V_D)$ is deduced from a bound on $\int_\Gamma|u'|^p\,dx$ together with $\|u\|_p=1$, which gives a bound in $W^{1,p}(\Gamma)$, and no Poincar\'e inequality is used. The proof of Proposition~\ref{ap:4} uses the full $W^{1,p}$ norm for this reason.
	\item Positivity of $\lambda_{1,p}$. It is used only in Proposition~\ref{ap:7}, which is not claimed here and is replaced by Proposition~\ref{ap:8}.
\end{enumerate}
The remaining ingredients, namely the $C^1$ structure of $M$, the uniform convexity of $W^{1,p}(\Gamma)$ required in \cite[Theorem 2.1]{cuesta2003minimax}, the $(S_+)$ property, Ekeland's principle and the compact embeddings of Theorems~\ref{at:3} and~\ref{at2:3}, do not depend on $V_D$. In Lemma~\ref{al:9} the sets $X_{\mathcal U}$ are defined with $X(\Gamma,\emptyset)=W^{1,p}(\Gamma)$, and the requirement $x_0\notin V_D$ is void. \qed
\end{proof}

By Corollary~\ref{ar:2} and Proposition~\ref{ap:100}, the first nontrivial curve $\mathfrak{C}$ passes through the point $(\lambda_{2,p}(\Gamma,\emptyset), \lambda_{2,p}(\Gamma,\emptyset))$, where
$\lambda_{2,p}(\Gamma,\emptyset)$ is the least eigenvalue of \eqref{a:76} larger than $0$, and $\lambda_{2,p}(\Gamma,\emptyset)$ has the following variational characterization:
\begin{align*}
\lambda_{2,p}(\Gamma,\emptyset)
= \inf_{P \in \mathbb{P}}\;\max_{u \in P([-1,+1])}
\int_{\Gamma}|u'|^p\,dx,
\end{align*}
where $\mathbb{P}$ is as introduced above.

One of the principal differences between the Dirichlet-Kirchhoff Fu\v{c}\'{i}k spectrum $\Sigma_{p,\Gamma}$ and the Neumann Fu\v{c}\'{i}k spectrum $\mathcal{F}_{p,\Gamma}$ concerns the asymptotic behavior of the first nontrivial curve $\mathfrak{C}$, which is governed by the number $\bar\lambda$ defined in \eqref{a12:3}. By the compact embedding of $W^{1,p}(\Gamma)$ into $C(\Gamma)$ (Theorem~\ref{at2:3}), the set of admissible functions in \eqref{a12:3} is weakly closed, and the infimum is attained. Moreover $\bar\lambda>0$: if $\bar\lambda=0$, a minimizer $u$ satisfies $u'=0$, so $u$ is a constant vanishing somewhere, that is, $u\equiv0$, contradicting $\|u\|_p=1$.  For $x\in\Gamma$ we put $\lambda(x)=\lambda_{1,p}(\Gamma,\{x\})$, the infimum \eqref{a1:3} over $X(\Gamma,\{x\})=\{u\in W^{1,p}(\Gamma):u(x)=0\}$; if $x$ lies inside an edge, we first insert it as a vertex of degree two, which does not change the problem, since the Kirchhoff condition at a vertex of degree two is the continuity of $u'$, and then impose the Dirichlet condition there. The admissible set in \eqref{a12:3} is the union over $x\in\Gamma$ of the unit spheres of $X(\Gamma,\{x\})$, so that $\bar\lambda=\inf_{x\in\Gamma}\lambda(x)$; and if $u$ is a minimizer in \eqref{a12:3} vanishing at $x_0$, then $\lambda(x_0)\le\int_\Gamma|u'|^p\,dx=\bar\lambda$. Hence
\begin{equation}\label{n:lbarmin}
\bar\lambda=\min_{x\in\Gamma}\lambda(x),
\end{equation}
and the minimum is attained at the zero of every minimizer in \eqref{a12:3}.

The following example illustrates the variational definition of \(\bar{\lambda}\) and shows how it can be computed in concrete cases.

\begin{ex}\label{aex:2}
Let $\Gamma=[0,L]$ be a single edge with Kirchhoff, that is Neumann, conditions at both ends. A function vanishing at $x_0$ splits into its restrictions to $[0,x_0]$ and $[x_0,L]$, so the infimum in \eqref{a12:3} over such functions is $\min\{\mu(x_0),\mu(L-x_0)\}$, where $\mu(t)$ is the first eigenvalue on an interval of length $t$ with a Neumann condition at one end and a Dirichlet condition at the other. By even reflection about the Neumann end, $\mu(t)$ is the first Dirichlet eigenvalue on an interval of length $2t$, that is, $\mu(t)=\frac{p}{p'}\bigl(\frac{\pi_p}{2t}\bigr)^p$ in the normalisation of Remark~\ref{ar:3}. As $\mu$ is decreasing,
\begin{align*}
\bar\lambda=\mu(L)=\frac{p}{p'}\Bigl(\frac{\pi_p}{2L}\Bigr)^p,
\end{align*}
attained exactly when $x_0$ is an endpoint, for instance by $u(x)=\sin_p\bigl(\pi_p(L-x)/(2L)\bigr)$ when $x_0=L$. For $p=2$, $\bar\lambda=\pi^2/(4L^2)$. A path graph with two edges of lengths $L_1,L_2$ is the interval of length $L_1+L_2$, because the Kirchhoff condition at a vertex of degree two is a transmission condition \cite[Example~2.2]{berkolaiko2017elementary}.
\end{ex}

We turn to the asymptotic behavior of
$\alpha(b) $.

\begin{lemma}{\label{al:12}}
Let $p>1$ and let $u$ be a minimizer in \eqref{a12:3}. Then $u$ does not change sign in $\Gamma$ and $u$ vanishes at exactly one point in $\Gamma$.
\end{lemma}
\begin{proof}
Let $u$ be a minimizer of \eqref{a12:3}. We first show that if $u$ vanishes at some $x_0 \in \Gamma$, then
\begin{equation}\label{a:83}
\int_\Gamma |u'|^{p-2}u'\phi'\,dx = \bar{\lambda} \int_\Gamma |u|^{p-2}u\phi\, dx
\end{equation}
for all $\phi \in W_{x_0}= \{\phi \in W^{1,p}(\Gamma): \phi(x_0)=0\}$. Since $u$ vanishes at $x_0$ and is a minimizer of \eqref{a12:3},
\begin{align*}
\bar{\lambda} = \inf \left\{ \int_{\Gamma} |u'|^p \, dx: \, u \in W_{x_0} \text{ and } \|u\|_p = 1 \right\},
\end{align*}
and this infimum is attained at $u$. By the Lagrange multiplier rule there exists $\lambda \in \mathbb{R}$ such that
\begin{align*}
\int_{\Gamma} |u'|^{p-2} u' \phi' \, dx
= \lambda \int_{\Gamma} |u|^{p-2} u \phi \, dx, \quad \text{for all } \phi \in W_{x_0}.
\end{align*}
Taking $\phi = u$ shows $\lambda = \bar{\lambda}$, and \eqref{a:83} follows. We next show that $u$ vanishes at exactly one point in $\Gamma$. Suppose, for contradiction, that $u$ vanishes at two distinct points $x_1$ and $x_2$ in $\Gamma$. The function $w = |u|$ is also a minimizer in \eqref{a12:3} and vanishes at $x_1$ and $x_2$. Therefore $w$ satisfies \eqref{a:83}
, that is,
\begin{equation}{\label{a:84}}
\int_{\Gamma} |w'|^{p-2} w' \phi' \, dx
= \bar{\lambda} \int_{\Gamma} |w|^{p-2} w \phi \, dx,
\end{equation}
for all $\phi \in W_{x_1}$ and similarly for all $\phi \in W_{x_2}$. Let $\chi(x)=\min\{1,d(x,x_2)/d(x_1,x_2)\}$, so that $\chi\in W^{1,p}(\Gamma)$, $\chi(x_1)=1$ and $\chi(x_2)=0$. Every $\phi\in W^{1,p}(\Gamma)$ is the sum of $\phi-\phi(x_1)\chi\in W_{x_1}$ and $\phi(x_1)\chi\in W_{x_2}$, and therefore $w$ satisfies \eqref{a:84} for all $\phi\in W^{1,p}(\Gamma)$. With $\phi \equiv1$, \eqref{a:84} becomes
\begin{align*}
\bar{\lambda} \int_{\Gamma} |w|^{p-2} w \, dx =0.
\end{align*}
As $\bar{\lambda} > 0$ and $w\ge0$, this forces $w\equiv0$, contradicting $\|w\|_p=1$. It remains to verify that the minimizer $u$ has constant sign. Assume, for instance, that $u^+ \not\equiv 0$ and $u(x_0)=0$. Then $u^+(x) = \max (u(x),0)$ satisfies $u^+(x_0)=0$, and \eqref{a:83} with $\phi = u^+$ reads
\begin{align*}
\int_{\Gamma} |(u^+)'|^p \, dx
&= \bar{\lambda} \int_{\Gamma} (u^+)^p \, dx.
\end{align*}
Thus $u^+/\|u^+\|_p$ is also a minimizer for \eqref{a12:3}, and by the first part of the proof $u^+$ vanishes only at $x_0$, so that $u \ge 0$ in $\Gamma$.    \qed
\end{proof}

{
For $x\in\Gamma$ let $\lambda(x)$ be as in \eqref{n:lbarmin}. By \eqref{n:lbarmin} and Lemma~\ref{al:12}, the minimum is attained exactly at the zeros of the minimizers in \eqref{a12:3}, each minimizer having exactly one zero. Where this zero lies depends on the topology of $\Gamma$. For $p=2$, imposing a Dirichlet condition at one point is a rank-one perturbation of the Kirchhoff Laplacian, and the interlacing inequalities of \cite{berkolaiko2019surgery} give $0=\lambda_{1,2}(\Gamma,\emptyset)\le\lambda(x)\le\lambda_{2,2}(\Gamma,\emptyset)$ for every $x\in\Gamma$, in accordance with Proposition~\ref{at:6} below. On an interval, $\bar\lambda$ is the first eigenvalue of a mixed Neumann--Dirichlet problem (Example~\ref{aex:2}); the first curve of the Fu\v{c}\'{\i}k spectrum with weights for the $p$-Laplacian with mixed boundary conditions is studied in \cite{leadi2007first}. The following proposition locates the zero on the equilateral star; it is not the point of symmetry of $\Gamma$.

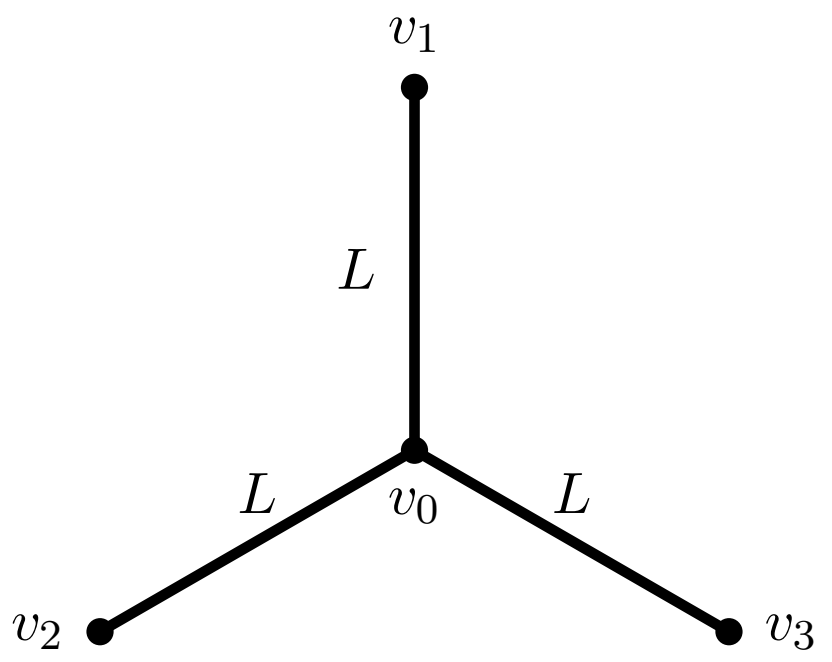
\begin{figure}[ht]
\begin{center}
\begin{tikzpicture}
	\draw[line width=2pt] (0,0) -- (0,2.4);
	\draw[line width=2pt] (0,0) -- (-2.08,-1.2);
	\draw[line width=2pt] (0,0) -- (2.08,-1.2);
	\fill (0,0) circle (0.09);
	\fill (0,2.4) circle (0.09);
	\fill (-2.08,-1.2) circle (0.09);
	\fill (2.08,-1.2) circle (0.09);
	\node[below] at (0,-0.12) {$v_0$};
	\node[above] at (0,2.5) {$v_1$};
	\node[left] at (-2.18,-1.2) {$v_2$};
	\node[right] at (2.18,-1.2) {$v_3$};
	\node[left] at (-0.12,1.2) {$L$};
	\node[above] at (-1.04,-0.55) {$L$};
	\node[above] at (1.04,-0.55) {$L$};
\end{tikzpicture}
\end{center}
\caption{The equilateral three-edge star of Proposition~\ref{aex:100}, with Kirchhoff conditions at all vertices.}
\label{n:figstar}
\end{figure}
\begin{prop}\label{aex:100}
Let $\Gamma$ be the star with central vertex $v_0$ and three pendant edges $e_1,e_2,e_3$ of length $L$ ending at $v_1,v_2,v_3$, with Kirchhoff conditions at all vertices (Figure~\ref{n:figstar}), and let $1<p<\infty$. Put $c_p=2^{1/(p-1)}$ and $\cot_p=\cos_p/\sin_p$, which is strictly decreasing on $(0,\pi_p)$ and vanishes at $\pi_p/2$ \cite{lang2011eigenvalues}. Let $z_p$ be the unique root in $(0,\pi_p/2)$ of
\begin{equation}\label{n:starp}
\cot_p(z)=c_p\cot_p\Bigl(\frac{\pi_p}{2}-z\Bigr).
\end{equation}
Then $z_p<\pi_p/4$,
\begin{align*}
\bar\lambda=\lambda(v_1)=\lambda(v_2)=\lambda(v_3)=\frac{p}{p'}\Bigl(\frac{z_p}{L}\Bigr)^p<2^{-p}\lambda(v_0),
\end{align*}
and $\lambda(x)>\bar\lambda$ for every other $x\in\Gamma$. Thus the minimum in \eqref{n:lbarmin} is attained exactly at the three pendant vertices, and not at the centre. For $p=2$, $\bar\lambda=\arctan^2(1/\sqrt2)/L^2\approx0.379/L^2$.
\end{prop}

\begin{proof}
If the zero is placed at the centre, each edge carries the mixed problem of Example~\ref{aex:2} on an interval of length $L$, so that $\lambda(v_0)=\frac{p}{p'}\bigl(\frac{\pi_p}{2L}\bigr)^p$. Now impose the Dirichlet condition at $v_1$ and write $\lambda=\frac{p}{p'}k^p$ with $0<kL<\pi_p/2$. Every solution of $-(\varphi_p(u'))'=\lambda\varphi_p(u)$ on an edge is a multiple of a translate of $\sin_p(k\,\cdot)$ \cite{lang2011eigenvalues}. We take $u=\sin_p(kt)$ on $e_1$, $t$ being the distance to $v_1$, and $u=B\sin_p(ks+\pi_p/2)$ on $e_2$ and $e_3$, $s$ being the distance to $v_2$, respectively $v_3$, so that $u(v_1)=0$ and $u'=0$ at $v_2,v_3$. Since $\sin_p(\pi_p-\theta)=\sin_p\theta$ and $\cos_p(\pi_p-\theta)=-\cos_p\theta$, continuity and the Kirchhoff condition at $v_0$ read
\begin{align*}
\sin_p(kL)=B\sin_p\Bigl(\frac{\pi_p}{2}-kL\Bigr),\qquad
\varphi_p\bigl(\cos_p(kL)\bigr)=2\,\varphi_p\Bigl(B\cos_p\Bigl(\frac{\pi_p}{2}-kL\Bigr)\Bigr),
\end{align*}
and eliminating $B$ gives \eqref{n:starp} with $z=kL$. On $(0,\pi_p/2)$ the left-hand side of \eqref{n:starp} decreases strictly from $+\infty$ to $0$ and the right-hand side increases strictly from $0$ to $+\infty$, so \eqref{n:starp} has exactly one root $z_p$ there. The corresponding $u$ is positive on $\Gamma\setminus\{v_1\}$, hence it is the first eigenfunction by Lemma~\ref{al:100}\textit{(iii)} with $V_D=\{v_1\}$, and $\lambda(v_1)=\frac{p}{p'}(z_p/L)^p$. At $z=\pi_p/4$ the left-hand side is smaller than the right-hand side, since $c_p>1$, so $z_p<\pi_p/4$ and $\lambda(v_1)<2^{-p}\lambda(v_0)$. By symmetry, $\lambda(v_2)=\lambda(v_3)=\lambda(v_1)$.

Let $x$ lie inside $e_i$. Then $\Gamma\setminus\{x\}$ has two components, the segment $S$ between $x$ and $v_i$, of length $t<L$, and the rest. For $u\in X(\Gamma,\{x\})$ the integrals of $|u'|^p$ and $|u|^p$ split over the two components, so $\lambda(x)\ge\min\{\mu(t),\lambda_B\}$. Here $\mu(t)>\mu(L)=\lambda(v_0)>\lambda(v_1)$ is the mixed eigenvalue of Example~\ref{aex:2}, and $\lambda_B$ is the infimum of the Rayleigh quotient over the nonzero $u\in W^{1,p}(\Gamma)$ that vanish on $S$. Such functions belong to $X(\Gamma,\{v_i\})$, so $\lambda_B\ge\lambda(v_i)=\lambda(v_1)$. 
The infimum $\lambda_B$ is attained by Theorem~\ref{at:3}, and if $\lambda_B=\lambda(v_i)$, a minimizer would be a first eigenfunction for $V_D=\{v_i\}$ vanishing on $S$, a set of positive length. This is impossible by Lemma~\ref{al:100}\textit{(ii)}, since $\Gamma\setminus\{v_i\}$ is connected. Hence $\lambda(x)>\lambda(v_1)$. Together with $\lambda(v_0)>\lambda(v_1)$, this proves that $\bar\lambda=\lambda(v_1)$ is attained exactly at $v_1,v_2,v_3$. For $p=2$ we have $c_2=2$ and $\cot(\pi/2-z)=\tan z$, so \eqref{n:starp} becomes $\tan^2z=\frac12$. \qed
\end{proof}}

\begin{prop}{\label{at:6}}
Let $p>1$. Then $\widetilde{\mathcal{F}}_{p,\Gamma}\subseteq(\bar\lambda,\infty)\times(\bar\lambda,\infty)$. In particular, $\alpha(b)>\bar\lambda$ for every $b\ge0$, and $\lambda_{2,p}(\Gamma,\emptyset)>\bar\lambda$.
\end{prop}
\begin{proof}
Let $(\lambda^+,\lambda^-)\in  \widetilde{\mathcal{F}}_{p,\Gamma}$ and let $u$ be a corresponding eigenfunction. The Neumann-Kirchhoff version of \eqref{a2:7} with the test function $\phi =1$ reads
\begin{align*}
\lambda^+\int_{\Gamma} (u^+)^{p-1} \, dx - \lambda^-\int_\Gamma (u^-)^{p-1} \, dx =0.
\end{align*}
If $u\ge0$, then $\int_\Gamma(u^-)^{p-1}\,dx=0$, and the identity gives $\lambda^+=0$ because $u\not\equiv0$, so that $(\lambda^+,\lambda^-)$ lies on a trivial line. The case $u\le0$ is the same. Hence $u$ changes sign in $\Gamma$. Taking $u^+$ and $u^-$ as test functions in \eqref{a2:7}, we get
\begin{equation}{\label{a:85}}
\lambda^+ = \frac{\int_{\Gamma} |(u^{+})'|^p \, dx}{\int_{\Gamma} (u^{+})^p \, dx} \ge \bar{\lambda}
\quad \text{and} \quad
\lambda^- = \frac{\int_{\Gamma} |(u^{-})'|^p \, dx}{\int_{\Gamma} (u^{-})^p \, dx} \ge \bar{\lambda}.
\end{equation}
Suppose that $\lambda^+ = \bar{\lambda}$. Then, by \eqref{a:85}, $u^+/\|u^+\|_p$ is a minimizer of \eqref{a12:3}. Since $u$ changes sign, $u^+$ vanishes at more than one point, which contradicts Lemma~\ref{al:12}, and the same argument applies to $\lambda^-$. For $b\ge0$ the point $(b+\alpha(b),\alpha(b))$ belongs to $\widetilde{\mathcal F}_{p,\Gamma}$ by Propositions~\ref{at4:5} and~\ref{ap:100}, since $\alpha(b)>0$. Hence $\alpha(b)>\bar\lambda$, and $b=0$ gives $\lambda_{2,p}(\Gamma,\emptyset)=\alpha(0)>\bar\lambda$. \qed
\end{proof}

\begin{prop}{\label{ap:8}}
Let $\alpha(b)$ be defined in \eqref{a:28} and let $V_D=\emptyset$. Then
\begin{align*}
\lim_{b\to+\infty}\alpha(b)=\bar{\lambda},
\end{align*}
where $\bar{\lambda}$ is given by \eqref{a12:3}.
\end{prop}

\begin{proof} By Proposition~\ref{at:6}, $\alpha(b)>\bar\lambda$ for every $b\ge0$, and therefore $\lim_{b\to\infty}\alpha(b)\ge\bar\lambda$. It remains to prove the reverse inequality. Part of the argument is taken from Proposition \ref{ap:7}. Suppose, for contradiction, that there exists $\epsilon>0$ such that
\begin{align*}
\max_{u \in P[-1,+1]}\tilde{I_b}(u) \ge \bar{\lambda} + \epsilon, \quad \text{ for all } P \in \mathbb{P} \text{ and all } b \ge 0.
\end{align*}
Choose a minimizer $u$ in \eqref{a12:3} with $u\le 0$ in $\Gamma$ and consider the path $P \in \mathbb{P}$ defined by
\begin{align*}
P(t) = \frac{t\phi_{1,p} + (1-|t|)u}{\|t\phi_{1,p} + (1-|t|)u\|_p} , \quad t \in [-1,+1].
\end{align*}
It is well defined because $u$ vanishes somewhere and $\phi_{1,p}$ does not vanish. By compactness of $P[-1,+1]$, the maximum of $\tilde{I_b}$ on $P[-1,+1]$ is attained at some $t_b \in [-1,+1]$. Writing $w_{t_b} = t_b \phi_{1,p}+ (1-|t_b|)u$, we have
\begin{equation}{\label{a:86}}
\int_\Gamma |w_{t_b}'|^p \, dx - b \int_\Gamma (w_{t_b}^+)^p \, dx \ge (\bar{\lambda}+ \epsilon) \int_\Gamma |t_b\phi_{1,p} + (1-|t_b|)u|^p \, dx
\end{equation}
for all $b \ge 0$. Letting $b \to +\infty$ along a subsequence, we may assume that $t_b \to \bar{t} \in [-1,+1]$. Since $w_{t_b}$ is bounded in $W^{1,p}(\Gamma)$ independently of $b$, dividing \eqref{a:86} by $b$ gives $\int_\Gamma (w_{t_b}^+)^p \, dx \to 0$ as $b \to +\infty$, and the convergence $w_{t_b}\to w_{\bar{t}}$ in $W^{1,p}(\Gamma)$ yields
\begin{align*}
\int_\Gamma \big((\bar{t}\,\phi_{1,p} + (1-|\bar{t}|)u)^+\big)^p \, dx = 0,
\end{align*}
which means $(w_{\bar{t}})^+ = 0$ in $\Gamma$, that is, $w_{\bar{t}}(x) \le 0 \text{ in } \Gamma$, or equivalently
\begin{align*}
\bar{t}\,\phi_{1,p} \le -(1-|\bar{t}|)u \quad \text{ in } \Gamma.
\end{align*}
Since $u$ vanishes somewhere and $\phi_{1,p} \equiv \text{ Constant} >0$, this inequality forces $\bar{t} \le 0$. From \eqref{a:86} we deduce that
\begin{equation}{\label{a:87}}
\int_\Gamma |w_{\bar{t}}'|^p \, dx \ge (\bar{\lambda}+ \epsilon) \int_\Gamma |\bar{t}\phi_{1,p} + (1-|\bar{t}|)u|^p \, dx
\end{equation}

Using $\phi'_{1,p} \equiv 0$, we have
$|w'_{\bar{t}}|^p = (1 - |\bar{t}|)^p |u'|^p$. Moreover, since $\bar{t}\phi_{1,p} \le 0$ and $(1 - |\bar{t}|)u \le 0$ in $\Gamma$,
the inequality $(a+b)^p \ge a^p + b^p$ for all $a,b \ge 0$ gives
\begin{align*}
|\bar{t}\phi_{1,p} + (1 - |\bar{t}|)u|^p
\ge |\bar{t}\phi_{1,p}|^p + |(1 - |\bar{t}|)u|^p.
\end{align*}
Inequality \eqref{a:87} therefore becomes
\begin{equation}{\label{a:88}}
\int_\Gamma (1 - |\bar{t}|)^p |u'|^p \, dx
\ge (\bar{\lambda} + \epsilon)
\int_\Gamma \left( |\bar{t}\phi_{1,p}|^p + (1 - |\bar{t}|)^p |u|^p \right) dx.
\end{equation}
If $\bar t=-1$, then \eqref{a:88} reads $0\ge(\bar\lambda+\epsilon)\int_\Gamma\phi_{1,p}^p\,dx$, a contradiction. If $\bar t\in(-1,0]$, then $0<(1-|\bar t|)^p\le1$, and, $u$ being a minimizer in \eqref{a12:3},
\begin{align*}
\int_\Gamma |u'|^p \, dx = \bar{\lambda} \int_\Gamma |u|^p \, dx.
\end{align*}
As $|\bar{t}|^p\int_\Gamma |\phi_{1,p}|^p \, dx \ge 0$, inequality \eqref{a:88} becomes
\begin{align*}
(1 - |\bar{t}|)^p\int_\Gamma |u'|^p \, dx
\ge (\bar{\lambda} + \epsilon) (1 - |\bar{t}|)^p\int_\Gamma |u|^p dx.
\end{align*}
Equivalently,
\begin{align*}
\bar{\lambda}(1 - |\bar{t}|)^p\int_\Gamma |u|^p \, dx
&\ge (\bar{\lambda} + \epsilon) (1 - |\bar{t}|)^p\int_\Gamma |u|^p dx,\\
\bar{\lambda} &\ge \bar{\lambda} + \epsilon,
\end{align*}
which is impossible. \qed \end{proof}

By Propositions~\ref{at:6} and~\ref{ap:8}, $\bar\lambda$ is the largest number $\mu$ such that $\widetilde{\mathcal F}_{p,\Gamma}\subseteq(\mu,\infty)\times(\mu,\infty)$, because $(b+\alpha(b),\alpha(b))\in\widetilde{\mathcal F}_{p,\Gamma}$ and $\alpha(b)\to\bar\lambda$ as $b\to\infty$.

\section{Differentiability of the first curve at the diagonal}\label{as:10}

Throughout this section $1<p<\infty$ and $\lambda_2=\lambda_{2,p}(\Gamma,V_D)$, with $V_D$ satisfying \eqref{ah:100}, possibly empty (Proposition~\ref{ap:100}).  By Corollary~\ref{ar:2} and \eqref{a:29}, $\mathfrak C$ passes through $(\lambda_2,\lambda_2)$ and is the union of the two arcs
\begin{align*}
\mathfrak C^{+}=\{(b+\alpha(b),\alpha(b)):b\ge0\},\qquad
\mathfrak C^{-}=\{(\alpha(b),b+\alpha(b)):b\ge0\},
\end{align*}
which are reflections of each other in the diagonal. With $\Lambda$ as in \eqref{n:Lambda}, $\alpha(b)\to\Lambda$ as $b\to\infty$ by Propositions~\ref{ap:7} and~\ref{ap:8}. By Proposition~\ref{ap:6} the map $b\mapsto b+\alpha(b)$ is a strictly increasing homeomorphism of $[0,\infty)$ onto $[\lambda_2,\infty)$ and $\alpha$ is a strictly decreasing homeomorphism of $[0,\infty)$ onto $(\Lambda,\lambda_2]$.
Hence $\mathfrak C^{+}\cap\mathfrak C^{-}=\{(\lambda_2,\lambda_2)\}$, the arc $\mathfrak C^{+}$ is the graph of the strictly decreasing homeomorphism
\begin{align*}
g:[\lambda_2,\infty)\to(\Lambda,\lambda_2],\qquad g(b+\alpha(b))=\alpha(b),
\end{align*}
and $\mathfrak C^{-}$ is the graph of $g^{-1}$. The curve $\mathfrak C$ is therefore the graph of the function
\begin{equation}\label{a10:1}
G:(\Lambda,\infty)\to(\Lambda,\infty),\qquad G(\lambda)=\begin{cases} g^{-1}(\lambda), & \Lambda<\lambda\le\lambda_2,\\[2pt] g(\lambda), & \lambda\ge\lambda_2,\end{cases}
\end{equation}
which is well defined because $g(\lambda_2)=\lambda_2$, is continuous and strictly decreasing, and is an involution, $G\circ G=\operatorname{id}$, with the single fixed point $\lambda_2$. We say that $\mathfrak C$ is differentiable at $(\lambda_2,\lambda_2)$ if $G$ is differentiable at $\lambda_2$.

\begin{remark}\label{ar10:1}
Two consequences of the involution property require no analysis. If the one-sided derivatives $G'(\lambda_2^{+})$ and $G'(\lambda_2^{-})$ exist, are finite and nonzero, then differentiating $G\circ G=\operatorname{id}$ from the right at $\lambda_2$, where $G$ maps $[\lambda_2,\lambda_2+h)$ onto $(G(\lambda_2+h),\lambda_2]$, gives
\begin{align*}
G'(\lambda_2^{-})\,G'(\lambda_2^{+})=1 .
\end{align*}
If $G$ is differentiable at $\lambda_2$, then $G'(\lambda_2)^2=1$, and since $G$ is decreasing, $G'(\lambda_2)=-1$. Consequently, whenever $\mathfrak C$ is smooth at the diagonal, it is orthogonal to it. The results below identify the two one-sided slopes.
\end{remark}
We use the notation of Section~\ref{as:4}, where $M$ and $\tilde I_b$ are defined, $\|\tilde I_b'(u)\|_*$ denotes the norm given in \eqref{a:17}, and $\mathbb{P}$ and $\alpha(b)$ are defined in \eqref{a:28}. The constant $c_e$ and the balls $B_\rho(v)\subseteq M$ are those of Lemma~\ref{n:proj}, and we write $\lambda_{1,p}=\lambda_{1,p}(\Gamma,V_D)$. We assume
\begin{equation}\label{a10:2}
\text{every eigenfunction of \eqref{a1:1} with eigenvalue }\lambda_2\text{ is a multiple of a fixed }\phi_{2,p}\in M .
\end{equation}
The critical points of $\tilde I_0$ on $M$ with critical value $\lambda_2$ are exactly the eigenfunctions in $M$ with eigenvalue $\lambda_2$. Thus \eqref{a10:2} says that the critical set of $\tilde I_0$ at the level $\lambda_2$ is $\{\phi_{2,p},-\phi_{2,p}\}$. For $p=2$, \eqref{a10:2} is the simplicity of $\lambda_2$. We put
\begin{equation}\label{a10:3}
A=\int_\Gamma(\phi_{2,p}^+)^p\,dx,\qquad B=\int_\Gamma(\phi_{2,p}^-)^p\,dx,\qquad A+B=1,\qquad r=\frac AB ,
\end{equation}
and note that $A,B>0$, since $\phi_{2,p}$ changes sign by Lemma~\ref{al:100} \textit{(iii)}, as $\lambda_2>\lambda_{1,p}$.

The quantitative deformation lemma and the non-critical interval property used below are Lemmas~\ref{al:deform} and~\ref{al:noncrit} of Section~\ref{as:4}.

\begin{lemma}\label{al10:1}
Let $0<\varepsilon<\lambda_2-\lambda_{1,p}$ and put
\begin{align*}
O_\varepsilon=\{u\in M:\tilde I_0(u)<\lambda_2-\varepsilon\}.
\end{align*}
Then $O_\varepsilon$ has exactly two path components, $C_\varepsilon^-\ni-\phi_{1,p}$ and $C_\varepsilon^+\ni\phi_{1,p}$, and both are open in $M$.
\end{lemma}

\begin{proof}
By Lemma~\ref{n:proj}, the path components of an open subset of $M$ are open, which is the last assertion, and any two points of $B_r(x)$, $x\in M$, $0<r\le1/(2c_e)$, are joined by a path in $B_{C_0r}(x)$ with $C_0=2(1+c_e\|x\|)$.

By Theorem \ref{at:1} with $b=0$ and Corollary \ref{ar:2} there is no eigenvalue in the interval $(\lambda_{1,p},\lambda_2)$, and $\lambda_{1,p}$ is simple (Lemma~\ref{al:100}), so the only critical points of $\tilde I_0$ with critical value below $\lambda_2$ are $\pm\phi_{1,p}$.

We claim that for every $r>0$ there is $\eta>0$ with
\begin{equation}\label{a10:5}
\{u\in M:\tilde I_0(u)\le\lambda_{1,p}+\eta\}\subseteq B_r(\phi_{1,p})\cup B_r(-\phi_{1,p}).
\end{equation}
Otherwise there are $u_n\in M$ with $\tilde I_0(u_n)\to\lambda_{1,p}$ and $\|u_n\mp\phi_{1,p}\|\ge r$. The sequence is bounded in $X(\Gamma,V_D)$, so along a subsequence $u_n\rightharpoonup u$ and, by Theorem \ref{at:3}, $u_n\to u$ in $L^p(\Gamma)$. Hence $u\in M$, and $\int_\Gamma|u'|^p\,dx\le\lambda_{1,p}$ by weak lower semicontinuity, so that $u=\pm\phi_{1,p}$ and $\int_\Gamma|u_n'|^p\,dx\to\int_\Gamma|u'|^p\,dx$. Weak convergence of $u_n'$ together with convergence of the norms implies $u_n'\to u'$ in $L^p(\Gamma)$ by the uniform convexity of $L^p(\Gamma)$, so $u_n\to\pm\phi_{1,p}$ in $X(\Gamma,V_D)$, contradicting $\|u_n\mp\phi_{1,p}\|\ge r$.

Since $\tilde I_0$ is continuous on $M$ and $\tilde I_0(\pm\phi_{1,p})=\lambda_{1,p}<\lambda_2-\varepsilon$, there is $\rho>0$ with $B_\rho(\pm\phi_{1,p})\subseteq O_\varepsilon$. Put
\begin{align*}
r=\min\Bigl\{\frac{\rho}{C_0},\ \frac1{2c_e}\Bigr\},\qquad C_0=2\bigl(1+c_e\|\phi_{1,p}\|\bigr),
\end{align*}
and let $\eta>0$ be given by \eqref{a10:5} for this $r$, reduced if necessary so that $\lambda_{1,p}+\eta<\lambda_2-\varepsilon$. The inclusion \eqref{a10:5} persists under such a reduction.

Let $u\in O_\varepsilon$. If $\tilde I_0(u)\le\lambda_{1,p}+\eta$, then
\eqref{a10:5} implies that
\[
u\in B_r(\phi_{1,p})\cup B_r(-\phi_{1,p}).
\]
If $\tilde I_0(u)>\lambda_{1,p}+\eta$, then
the interval
$[\lambda_{1,p}+\eta,\tilde I_0(u)]$
contains no critical value of $\tilde I_0$, and
hence, by Lemma~\ref{al:noncrit}, there exists
a path joining $u$ to some point $y$ with
$\tilde I_0(y)\le\lambda_{1,p}+\eta$, along
which $\tilde I_0$ is non-increasing.
In particular, this path is contained in
$\{\tilde I_0\le\tilde I_0(u)\} \subseteq O_\varepsilon$. Moreover, \eqref{a10:5} implies
\[
y\in B_r(\phi_{1,p})\cup B_r(-\phi_{1,p}).
\]
Thus, in either case, $u$ can be joined
within $O_\varepsilon$ to a point in
$B_r(\phi_{1,p})\cup B_r(-\phi_{1,p})$.
By \eqref{a10:4}, this point can
be joined to $\phi_{1,p}$ or $-\phi_{1,p}$
by a path contained in
\[
B_{C_0r}(\pm\phi_{1,p})
\subseteq B_\rho(\pm\phi_{1,p})
\subseteq O_\varepsilon.
\]
Every point of $O_\varepsilon$ therefore
belongs to the path component of either
$\phi_{1,p}$ or $-\phi_{1,p}$, and hence
$O_\varepsilon$ has at most
two path components.

It has at least two: a path in $O_\varepsilon$ from $-\phi_{1,p}$ to $\phi_{1,p}$ would belong to $\mathbb P$ and satisfy
\begin{align*}
\max_{P[-1,1]}\tilde I_0\le\lambda_2-\varepsilon<\lambda_2=\alpha(0),
\end{align*}
contradicting the definition \eqref{a:28} of $\alpha(0)$. \qed
\end{proof}

\begin{theorem}\label{at10:1}
Let $1<p<\infty$ and assume \eqref{a10:2}. Then $\alpha$ is differentiable from the right at $b=0$, with
\begin{equation}\label{a10:6}
\alpha'(0^+)=-\max\{A,B\}.
\end{equation}
Consequently, with $G$ as in \eqref{a10:1},
\begin{align*}
G'(\lambda_2^{+})=-\max\Bigl\{r,\frac1r\Bigr\},\qquad G'(\lambda_2^{-})=-\min\Bigl\{r,\frac1r\Bigr\},
\end{align*}
and $\mathfrak C$ is differentiable at $(\lambda_2,\lambda_2)$ if and only if $r=1$, that is, if and only if
\begin{align*}
\int_\Gamma(\phi_{2,p}^+)^p\,dx=\int_\Gamma(\phi_{2,p}^-)^p\,dx;
\end{align*}
in that case the slope is $-1$.
\end{theorem}

\begin{proof}
For $\rho>0$ put
\begin{align*}
N_\rho=B_\rho(\phi_{2,p})\cup B_\rho(-\phi_{2,p}),
\end{align*}
and fix $\rho\in(0,\min\{\|\phi_{2,p}\|,1/c_e\})$, so that the two balls are disjoint. For $u,v\in M$ one has $|u^+-v^+|\le|u-v|$ pointwise, hence $\|u^+-v^+\|_p\le\|u-v\|_p$. The mean value theorem for $t\mapsto t^p$ on $[0,1]$, together with $\|u^+\|_p,\|v^+\|_p\le1$, gives
\begin{equation*}
\Bigl|\int_\Gamma(u^+)^p\,dx-\int_\Gamma(v^+)^p\,dx\Bigr|\le p\,\|u-v\|_p\le p\,c_e\|u-v\| .
\end{equation*}
Since $(-\phi_{2,p})^+=\phi_{2,p}^-$, the quantity $\int_\Gamma(u^+)^p\,dx$ differs from $A$ by at most $p\,c_e\rho$ on $B_\rho(\phi_{2,p})$, and from $B$ by at most $p\,c_e\rho$ on $B_\rho(-\phi_{2,p})$. We write $\kappa=p\,c_e\rho$.

\emph{Step 1:} We claim that there are $\beta>0$ and $\varepsilon_0>0$, both depending on $\rho$, such that
\begin{equation}\label{a10:7}
\|\tilde I_0'(u)\|_*\ge\beta\qquad\text{for all }u\in M\setminus N_{\rho/8}\text{ with }|\tilde I_0(u)-\lambda_2|\le2\varepsilon_0 .
\end{equation}
Otherwise there are $u_n\in M\setminus N_{\rho/8}$ with $\tilde I_0(u_n)\to\lambda_2$ and $\|\tilde I_0'(u_n)\|_*\to0$, and Lemma \ref{al:3} provides a subsequence converging to a critical point at the level $\lambda_2$. The limit still satisfies $\|u\mp\phi_{2,p}\|\ge\rho/8$, so it is neither $\phi_{2,p}$ nor $-\phi_{2,p}$, contradicting \eqref{a10:2}.
Fix
\begin{align*}
\delta=\frac{\rho}{8C_*},\qquad S=M\setminus N_{\rho/2},
\end{align*}
where $C_*$ is the constant of Lemma~\ref{al:deform} for $b=0$, $c=\lambda_2$. Then $\delta\le\rho/16\le1/(16c_e)$. Reduce $\varepsilon_0$ so that
\begin{align*}
\frac{10\varepsilon_0}{\delta}\le\beta,\qquad {\varepsilon_0\le1},\qquad 2\varepsilon_0<\lambda_2-\lambda_{1,p}.
\end{align*}
We claim that the hypothesis \eqref{a:defhyp} of Lemma \ref{al:deform} holds with $c=\lambda_2$ for every $\varepsilon\in(0,\varepsilon_0]$. Let such an $\varepsilon$ be given and let $u\in M$ satisfy $\operatorname{dist}(u,S)\le {2C_*\delta=\rho/4}$ and $|\tilde I_0(u)-\lambda_2|\le2\varepsilon$. Choose $w\in S$ with $\|u-w\|< {3\rho/8}$. Then $w\notin N_{\rho/2}$, that is, $\|w\mp\phi_{2,p}\|\ge\rho/2$, and therefore
\begin{align*}
\|u\mp\phi_{2,p}\|\ge\|w\mp\phi_{2,p}\|-\|u-w\|>{\frac{\rho}{2}-\frac{3\rho}{8}=\frac{\rho}{8}},
\end{align*}
so that $u\notin N_{{\rho/8}}$. Moreover $|\tilde I_0(u)-\lambda_2|\le2\varepsilon\le2\varepsilon_0$, so \eqref{a10:7} applies to $u$ and gives
\begin{align*}
\|\tilde I_0'(u)\|_*\ge\beta\ge {\frac{10\varepsilon_0}{\delta}\ge\frac{10\varepsilon}{\delta}},
\end{align*}
which is \eqref{a:defhyp}. We write $\eta_\varepsilon$ for the corresponding deformation. The bound $\|\eta_\varepsilon(t,u)-u\|\le {C_*\delta=\rho/8}$ in \eqref{a:defconc} is uniform in $\varepsilon$, because $\delta$ and $C_*$ do not depend on $\varepsilon$. The choice $2\varepsilon_0<\lambda_2-\lambda_{1,p}$ ensures $|\tilde I_0(\pm\phi_{1,p})-\lambda_2|=\lambda_2-\lambda_{1,p}>2\varepsilon$, so $\eta_\varepsilon$ fixes $\pm\phi_{1,p}$ by \eqref{a:defconc} and $\eta_\varepsilon(1,\cdot)$ maps $\mathbb P$ into $\mathbb P$. Finally, $u\notin N_\rho$ implies $u\in S$, because $N_{\rho/2}\subseteq N_\rho$.

\emph{Step 2:} Let $0<b\le\varepsilon_0$ and $P\in\mathbb P$. On $M$ we have $\int_\Gamma(u^+)^p\,dx\le1$, hence $\tilde I_b\ge\tilde I_0-b$. If $\max_P\tilde I_0>\lambda_2+\varepsilon_0$, then
\begin{align*}
\max_P\tilde I_b>\lambda_2+\varepsilon_0-b\ge\lambda_2 .
\end{align*}
Otherwise $P\subseteq\{\tilde I_0\le\lambda_2+\varepsilon_0\}$. Put $\tilde P=\eta_{\varepsilon_0}(1,P)\in\mathbb P$. By \eqref{a:28} we have $\max_{\tilde P}\tilde I_0\ge\alpha(0)=\lambda_2$, and the maximum is attained at $\eta_{\varepsilon_0}(1,u)$ for some $u\in P$. If $u\notin N_\rho$, then $u\in S$ and \eqref{a:defconc} gives $\tilde I_0(\eta_{\varepsilon_0}(1,u))\le\lambda_2-\varepsilon_0$, which is impossible. Hence $u\in N_\rho$, and since $\tilde I_0$ does not increase along the deformation,
\begin{align*}
\tilde I_0(u)\ \ge\ \tilde I_0\bigl(\eta_{\varepsilon_0}(1,u)\bigr)\ \ge\ \lambda_2 .
\end{align*}
Consequently
\begin{align*}
\max_P\tilde I_b\ \ge\ \tilde I_b(u)=\tilde I_0(u)-b\int_\Gamma(u^+)^p\,dx\ \ge\ \lambda_2-b\bigl(\max\{A,B\}+\kappa\bigr).
\end{align*}
In both cases the right-hand side bounds $\max_P\tilde I_b$ from below, so taking the infimum over $P\in\mathbb P$,
\begin{equation}\label{a10:8}
\alpha(b)\ \ge\ \lambda_2-b\bigl(\max\{A,B\}+\kappa\bigr)\qquad(0<b\le\varepsilon_0).
\end{equation}

\emph{Step 3:} Let $\varepsilon'\in(0,\varepsilon_0]$ and $\varepsilon\in(0,\varepsilon']$. By \eqref{a:28} there is $P\in\mathbb P$ with $\max_P\tilde I_0\le\lambda_2+\varepsilon$. Put $\tilde P=\eta_{\varepsilon'}(1,P)\in\mathbb P$, so that $\max_{\tilde P}\tilde I_0\le\lambda_2+\varepsilon$. We claim that
\begin{equation}\label{a10:9}
u\in\tilde P\ \text{ and }\ \tilde I_0(u)>\lambda_2-\varepsilon'\quad\Longrightarrow\quad u\in N_\rho .
\end{equation}
If $u\notin N_\rho$, its preimage $u_0\in P$ satisfies $\|u_0-u\|\le {C_*\delta}=\rho/8$ by \eqref{a:defconc}, and therefore
\begin{align*}
\|u_0\mp\phi_{2,p}\|\ge\|u\mp\phi_{2,p}\|-\|u-u_0\|\ge\rho-\frac{\rho}{8}=\frac{7\rho}{8}>\frac{\rho}{2},
\end{align*}
so that $u_0\notin N_{\rho/2}$, that is, $u_0\in S$. Since also $\tilde I_0(u_0)\le\lambda_2+\varepsilon\le\lambda_2+\varepsilon'$, \eqref{a:defconc} gives $\tilde I_0(u)\le\lambda_2-\varepsilon'$.

Apply Lemma \ref{al10:1} with $\varepsilon'/2$ in place of $\varepsilon$, and write $O=O_{\varepsilon'/2}$ and $C^\pm=C^\pm_{\varepsilon'/2}$. Parametrise $\tilde P$ on $[-1,1]$ and put
\begin{align*}
t^*=\sup\{t:\tilde P(t)\in C^-\},\qquad t^{**}=\inf\{t>t^*:\tilde P(t)\in C^+\}.
\end{align*}
The sets $\tilde P^{-1}(C^-)$ and $\tilde P^{-1}(C^+)$ are disjoint and open in $[-1,1]$ and contain $-1$ and $1$ respectively. Then $\tilde P(t^*)\notin C^-$, because $\tilde P^{-1}(C^-)$ is open and $t^*$ is its supremum, and $\tilde P(t^*)\notin C^+$, because $t^*$ is a limit of points of $\tilde P^{-1}(C^-)$. Thus $\tilde P(t^*)\notin O$, that is, $\tilde I_0(\tilde P(t^*))\ge\lambda_2-\varepsilon'/2$, and $\tilde P(t^*)\in N_\rho$ by \eqref{a10:9}. Let $U$ be the ball among $B_\rho(\phi_{2,p})$ and $B_\rho(-\phi_{2,p})$ containing $\tilde P(t^*)$. The same argument applies at $t^{**}$. For $t\in(t^*,t^{**})$ the point $\tilde P(t)$ lies in neither $C^-$ nor $C^+$, so $\tilde P([t^*,t^{**}])$ is a connected subset of $M\setminus O$, therefore of $N_\rho$. Being connected and containing $\tilde P(t^*)$, it is contained in $U$. As $U$ is open and $\tilde P$ continuous, there are $a<t^*$ and $c>t^{**}$ with
\begin{align*}
\tilde P(a)\in C^-,\qquad \tilde P(c)\in C^+,\qquad \tilde P([a,c])\subseteq U .
\end{align*}
Let $Q\in\mathbb P$ be the concatenation of a path in $C^-$ from $-\phi_{1,p}$ to $\tilde P(a)$, of $\tilde P|_{[a,c]}$, and of a path in $C^+$ from $\tilde P(c)$ to $\phi_{1,p}$. Then
\begin{equation}\label{a10:10}
Q\subseteq O\cup U,\qquad \tilde I_0\le\lambda_2+\varepsilon\ \text{ on }Q,\qquad \tilde I_0<\lambda_2-\frac{\varepsilon'}{2}\ \text{ on }Q\setminus U .
\end{equation}
Replacing $Q$ by $t\mapsto-Q(-t)$, which again belongs to $\mathbb P$ and along which $\tilde I_0$ takes the same values because $I_0$ is even, we may prescribe which of the two balls is $U$.

\emph{Step 4: } Let $0<b\le\varepsilon_0/2$, put $\varepsilon'=2b\max\{A,B\}$, so that $\varepsilon'\le\varepsilon_0$ because $\max\{A,B\}\le1$, let $\varepsilon\in(0,\varepsilon']$, and let $Q$ be the path of Step 3 with
\begin{align*}
U=B_\rho(\phi_{2,p})\ \text{ if }A\ge B,\qquad U=B_\rho(-\phi_{2,p})\ \text{ if }A<B .
\end{align*}
On $Q\setminus U$ we have, by \eqref{a10:10},
\begin{align*}
\tilde I_b\le\tilde I_0<\lambda_2-\frac{\varepsilon'}{2}=\lambda_2-b\max\{A,B\},
\end{align*}
while on $Q\cap U$ we have $\int_\Gamma(u^+)^p\,dx\ge\max\{A,B\}-\kappa$ and $\tilde I_0\le\lambda_2+\varepsilon$, so that
\begin{align*}
\tilde I_b(u)\le\lambda_2+\varepsilon-b\bigl(\max\{A,B\}-\kappa\bigr).
\end{align*}
Hence $\alpha(b)\le\max_Q\tilde I_b\le\lambda_2-b(\max\{A,B\}-\kappa)+\varepsilon$, and since $\varepsilon\in(0,\varepsilon']$ was arbitrary,
\begin{equation}\label{a10:11}
\alpha(b)\ \le\ \lambda_2-b\bigl(\max\{A,B\}-\kappa\bigr)\qquad(0<b\le\varepsilon_0/2).
\end{equation}

\emph{Step 5:} By \eqref{a10:8} and \eqref{a10:11}, for every admissible $\rho$,
\begin{align*}
-\max\{A,B\}-\kappa\ \le\ \liminf_{b\to0^+}\frac{\alpha(b)-\lambda_2}{b}\ \le\ \limsup_{b\to0^+}\frac{\alpha(b)-\lambda_2}{b}\ \le\ -\max\{A,B\}+\kappa ,
\end{align*}
where $\kappa=p\,c_e\rho$. Letting $\rho\to0$ gives \eqref{a10:6}. The limit $b\to0^+$ is taken first: $\beta$ and $\varepsilon_0$ depend on $\rho$. Finally, for $b>0$ and $\lambda=b+\alpha(b)$,
\begin{align*}
\frac{G(\lambda)-G(\lambda_2)}{\lambda-\lambda_2}=\frac{\alpha(b)-\lambda_2}{b+\alpha(b)-\lambda_2}=\frac{\bigl(\alpha(b)-\lambda_2\bigr)/b}{1+\bigl(\alpha(b)-\lambda_2\bigr)/b}\ \longrightarrow\ \frac{-\max\{A,B\}}{1-\max\{A,B\}}=-\frac{\max\{A,B\}}{\min\{A,B\}}
\end{align*}
as $b\to0^+$, by \eqref{a10:3} and $\min\{A,B\}>0$. Since $b\mapsto b+\alpha(b)$ is a homeomorphism of $[0,\infty)$ onto $[\lambda_2,\infty)$, this limit is $G'(\lambda_2^{+})$. The value of $G'(\lambda_2^{-})$ follows from Remark \ref{ar10:1}, and differentiability of $G$ at $\lambda_2$ amounts to the equality of the two one-sided derivatives, that is, to $r=1$. \qed
\end{proof}

{
\begin{remark}\label{ar10:2}
For $p=2$ and $\lambda_2$ simple, Theorem~\ref{at10:1} also follows from known results. Near $(\lambda_2,\lambda_2)$ the nontrivial part of $\Sigma_{2,\Gamma}$ consists of two curves through $(\lambda_2,\lambda_2)$, given by the Lyapunov--Schmidt reduction of \cite{gallouet1981resultats,ruf1981nonlinear}, which applies on $\Gamma$ since $X(\Gamma,V_D)$ embeds compactly in $L^2(\Gamma)$; pairing the equation with $\pm\phi_{2,2}$ shows that their slopes are $-r$ and $-1/r$. A curve with slope $-s$ meets the line $\lambda^+=\lambda^-+b$ at $\lambda^-=\lambda_2-\frac{s}{1+s}b+o(b)$, which decreases in $s$, and $(b+\alpha(b),\alpha(b))$ is the first nontrivial point on this line by Theorem~\ref{at:1}; hence $\mathfrak C$ follows the steeper curve where $\lambda^+>\lambda^-$ and the other one where $\lambda^+<\lambda^-$, that is, it is the minimal curve through $(\lambda_2,\lambda_2)$ in the sense of \cite[Section~3]{perera2015note}, where the slopes are not computed. For $p=2$ the content of Theorem~\ref{at10:1} is thus a proof without the reduction, together with the case of a multiple eigenvalue (Theorem~\ref{n:p2}). For $p\neq2$ no such reduction is available, and we are not aware of a result of this type; the curves through the eigenvalues of the $p$-Laplacian constructed in \cite{micheletti2001fucik} come without information on their slopes. In finite dimension, Holubov\'a and Ne\v{c}esal \cite{holubova2025fucik} obtain, for a symmetric matrix and a simple eigenvalue with eigenvector $u_0$, the slopes $-\|u_0^+\|^2/\|u_0^-\|^2$ for $\pm u_0$, in agreement with Theorem~\ref{at10:1}; see Remark~\ref{n:hn} for multiple eigenvalues.
\end{remark}
}

\begin{corollary}\label{ac10:1}
Let $1<p<\infty$ and assume \eqref{a10:2}. Let $\sigma:\Gamma\to\Gamma$ be a bijection that maps vertices onto vertices and each closed edge isometrically onto a closed edge, with $\sigma(V_D)=V_D$. Then $\phi_{2,p}\circ\sigma=\pm\phi_{2,p}$, and if $\phi_{2,p}\circ\sigma=-\phi_{2,p}$, then $\mathfrak C$ is differentiable at $(\lambda_2,\lambda_2)$, with slope $-1$. This is the case on an interval $[0,L]$ whose endpoints either both belong to $V_D$ or both do not, $\sigma$ being the reflection $x\mapsto L-x$.
\end{corollary}

\begin{figure}[ht]
\begin{center}
\begin{tikzpicture}

	\draw[line width=2pt] (0,0) -- (0,2.8);
	\draw[line width=2pt] (0,0) -- (-2.17,-1.25);
	\draw[line width=2pt] (0,0) -- (1.39,-0.8);

	\fill (0,0) circle (0.09);

	\draw[thick,fill=white] (0,2.8) circle (0.09);
	\draw[thick,fill=white] (-2.17,-1.25) circle (0.09);
	\draw[thick,fill=white] (1.39,-0.8) circle (0.09);

	\node[below right] at (0.05,0.25) {$v_0$};
	\node[above] at (0,2.92) {$v_1$};
	\node[left] at (-2.27,-1.25) {$v_2$};
	\node[right] at (1.49,-0.8) {$v_3$};

	\node[left] at (-0.12,1.4) {$\ell_1$};
	\node[below] at (-1.08,-0.62) {$\ell_2$};
	\node[below] at (0.78,-0.42) {$\ell_3$};

\end{tikzpicture}
\end{center}
\caption{The three-edge star with $\ell_1>\ell_2>\ell_3$, Dirichlet conditions at the pendant vertices $v_1,v_2,v_3$ and the Kirchhoff condition at $v_0$.}
\label{afig:star}
\end{figure}
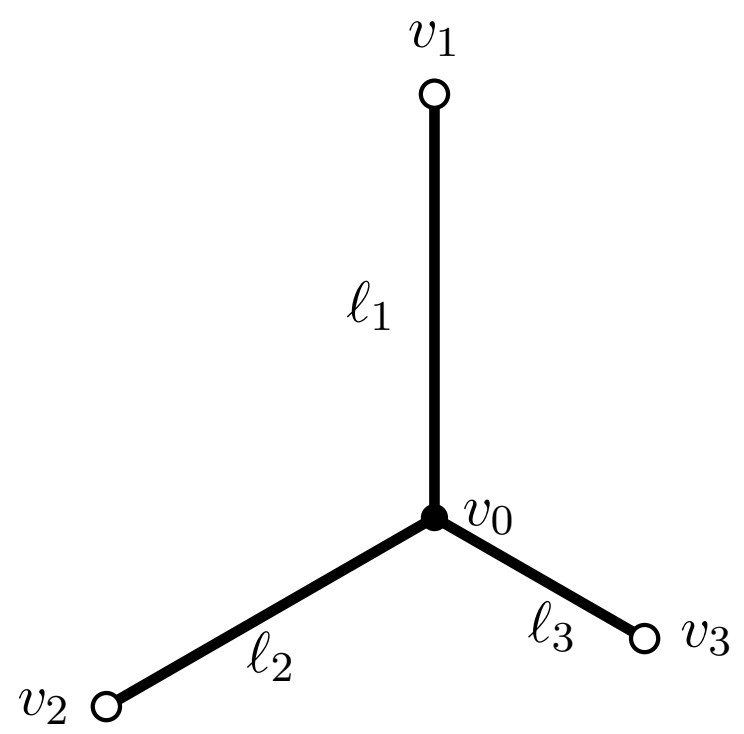

\begin{proof}
By hypothesis $\sigma$ maps $V_D$ onto $V_D$, hence $V(\Gamma)\setminus V_D$ onto itself, so it preserves the vertex conditions, and $u\mapsto u\circ\sigma$ is a linear bijection of $X(\Gamma,V_D)$ onto itself. On each edge $\sigma$ is an isometry, so that $(u\circ\sigma)'=\pm u'\circ\sigma$ with a sign depending only on the edge, and the change of variables $y=\sigma(x)$ gives, for $u,\phi\in X(\Gamma,V_D)$,
\begin{align*}
\int_\Gamma|(u\circ\sigma)'|^{p-2}(u\circ\sigma)'(\phi\circ\sigma)'\,dx &=\int_\Gamma|u'|^{p-2}u'\phi'\,dy,\\
\int_\Gamma|u\circ\sigma|^{p-2}(u\circ\sigma)(\phi\circ\sigma)\,dx&=\int_\Gamma|u|^{p-2}u\,\phi\,dy .
\end{align*}
As $\phi\circ\sigma$ runs over $X(\Gamma,V_D)$ when $\phi$ does, $u\circ\sigma$ is an eigenfunction with eigenvalue $\lambda$ whenever $u$ is, and $\|u\circ\sigma\|_p=\|u\|_p$. By \eqref{a10:2}, $\phi_{2,p}\circ\sigma=t\phi_{2,p}$ with $|t|=1$. If $t=-1$, then $\phi_{2,p}^-=(\phi_{2,p}\circ\sigma)^+=\phi_{2,p}^+\circ\sigma$, so that
\begin{align*}
B=\int_\Gamma(\phi_{2,p}^-)^p\,dx=\int_\Gamma(\phi_{2,p}^+\circ\sigma)^p\,dx=\int_\Gamma(\phi_{2,p}^+)^p\,dx=A,
\end{align*}
that is, $r=1$, and Theorem~\ref{at10:1} applies.

Now let $\Gamma=[0,L]$ with $V_D=\{0,L\}$ or $V_D=\emptyset$, so that the vertex condition is $u=0$, respectively $u'=0$, at both endpoints, and let $u$ be an eigenfunction with eigenvalue $(p-1)k^p$, $k>0$. We use the functions $\sin_p$ and $\cos_p=\sin_p'$ recalled in the paragraph preceding Figure~\ref{afig:star}. Here $\sin_p(\pi_p/2)=1$, and the zeros of $\cos_p$ are the points of $\frac{\pi_p}2+\pi_p\mathbb Z$. By uniqueness for the initial value problem \cite{lang2011eigenvalues}, $u(x)=c\sin_p(kx)$ if $V_D=\{0,L\}$ and $u(x)=c\sin_p(kx+\pi_p/2)$ if $V_D=\emptyset$, with $c\neq0$, and the condition at $x=L$ gives $kL\in\pi_p\mathbb N$. Hence every eigenvalue is simple, \eqref{a10:2} holds, and up to a constant factor,
\begin{align*}
\phi_{2,p}(x)=\sin_p\Bigl(\frac{2\pi_p x}{L}\Bigr)\quad\text{if }V_D=\{0,L\},\qquad \phi_{2,p}(x)=\sin_p\Bigl(\frac{\pi_p x}{L}+\frac{\pi_p}2\Bigr)\quad\text{if }V_D=\emptyset .
\end{align*}
Since $\sin_p$ is odd and $\sin_p(\theta+\pi_p)=-\sin_p\theta$, we have $\sin_p(2\pi_p-\theta)=-\sin_p\theta$ and $\sin_p(\pi_p-\theta)=\sin_p\theta$ for all $\theta$, and therefore $\phi_{2,p}(L-x)=-\phi_{2,p}(x)$ in both cases. \qed
\end{proof}

For $b\ge0$ let $K_b$ be the set of critical points of $\tilde I_b$ on $M$ at the level $\alpha(b)$, and put
\begin{align*}
m_b=\min_{u\in K_b}\int_\Gamma(u^+)^p\,dx,\qquad M_b=\max_{u\in K_b}\int_\Gamma(u^+)^p\,dx .
\end{align*}
Under \eqref{a10:2} we have $K_0=\{\phi_{2,p},-\phi_{2,p}\}$, $m_0=\min\{A,B\}$ and $M_0=\max\{A,B\}$. The next result holds for every $b\ge0$ and without \eqref{a10:2}.

\begin{prop}\label{n:bpos}
Let $1<p<\infty$ and $b\ge0$. Then $K_b$ is nonempty and compact, $0<m_b\le M_b<1$, and
\begin{equation}\label{n:dini}
-M_b\le\liminf_{h\to0}\frac{\alpha(b+h)-\alpha(b)}{h}\le\limsup_{h\to0}\frac{\alpha(b+h)-\alpha(b)}{h}\le-m_b,
\end{equation}
where $h\to0^+$ if $b=0$. In particular, if $u\mapsto\int_\Gamma(u^+)^p\,dx$ is constant on $K_b$, for instance if $K_b$ is a single point, then $\alpha$ is differentiable at $b$, from the right if $b=0$, and $\alpha'(b)=-m_b=-M_b$.
\end{prop}

\begin{proof}
By Proposition~\ref{ap:3}, $\alpha(b)$ is a critical value of $\tilde I_b$, so $K_b\neq\emptyset$. A sequence in $K_b$ is a Palais--Smale sequence at the level $\alpha(b)$, so by Lemma~\ref{al:3} it has a subsequence converging to a critical point at this level; hence $K_b$ is compact. Let $u\in K_b$. By the Lagrange multiplier rule, $u$ is a solution of \eqref{a1:1} with $(\lambda^+,\lambda^-)=(b+\alpha(b),\alpha(b))$. If $u$ did not change sign, it would be an eigenfunction of constant sign with eigenvalue $b+\alpha(b)$ or $\alpha(b)$, both larger than $\lambda_{1,p}$ by Proposition~\ref{ap:3}, which contradicts Lemma~\ref{al:100} \textit{(iii)}. Hence $0<\int_\Gamma(u^+)^p\,dx<1$, and $0<m_b\le M_b<1$ by compactness.

Fix $\rho\in(0,1/c_e)$, put $N_r(K_b)=\{u\in M:\operatorname{dist}(u,K_b)<r\}$ for $r>0$, and $\kappa=p\,c_e\rho$. By the first inequality in the proof of Theorem~\ref{at10:1},
\begin{equation}\label{n:mass}
m_b-\kappa\le\int_\Gamma(u^+)^p\,dx\le M_b+\kappa\qquad\text{for }u\in N_\rho(K_b).
\end{equation}

We proceed as in Step~1 of the proof of
Theorem~\ref{at10:1}, replacing
$\tilde I_0$, $\lambda_2$, and $N_r$ by
$\tilde I_b$, $\alpha(b)$, and $N_r(K_b)$,
respectively.

Recall that the set of critical points of
$\tilde I_b$ at the level $\alpha(b)$ is
precisely $K_b$. Therefore, by
Lemma~\ref{al:3}, there exist constants
$\beta>0$ and $\varepsilon_0>0$ such that
$\|\tilde I_b'(u)\|_*\geq\beta$
for every $u\in M\setminus N_{\rho/8}(K_b)$
with $|\tilde I_b(u)-\alpha(b)|
\leq 2\varepsilon_0$. Next, choose
$\delta=\frac{\rho}{8C_*}, \text{ and }
S=M\setminus N_{\rho/2}(K_b)$,
where $C_*$ is the constant appearing in
Lemma~\ref{al:deform}, corresponding to
$b$ and $c=\alpha(b)$.
By decreasing $\varepsilon_0$ if necessary,
we may assume that
\[
\frac{10\varepsilon_0}{\delta}\leq\beta,
\qquad
\varepsilon_0\leq 1,
\qquad
2\varepsilon_0<\alpha(b)-\lambda_{1,p}.
\]
We can now apply the same argument as in
Step~1 of Theorem~\ref{at10:1}.
Indeed, the distance function
$u\mapsto\operatorname{dist}(u,K_b)$
is $1$-Lipschitz, so the estimates used
there remain valid with $N_r(K_b)$
in place of $N_r$. Consequently, for every
$\varepsilon\in(0,\varepsilon_0]$,
Lemma~\ref{al:deform} provides a deformation
$\eta_\varepsilon$ satisfying
\[
\|\eta_\varepsilon(t,u)-u\|
\leq\frac{\rho}{8}
\]
for all $t\in[0,1]$ and $u\in M$. Moreover $\tilde I_b(-\phi_{1,p})=\lambda_{1,p}$ and $\tilde I_b(\phi_{1,p})=\lambda_{1,p}-b$ are smaller than $\alpha(b)-2\varepsilon$, and hence the map $\eta_\varepsilon(1,\cdot)$ sends $\mathbb P$ into $\mathbb P$.

\emph{Lower bounds.} 
Let $0<|h|\leq\varepsilon_0/2$ with $b+h\geq0$,
and fix $P\in\mathbb P$. Since
\[
\tilde I_{b+h}(u)
=\tilde I_b(u)-h\int_\Gamma(u^+)^p\,dx,
\qquad
0\leq\int_\Gamma(u^+)^p\,dx\leq1
\]
for all $u\in M$, we consider two cases.
If $\max_P\tilde I_b>\alpha(b)+\varepsilon_0$,
then
\[
\max_P\tilde I_{b+h}
>\alpha(b)+\varepsilon_0-|h|
\geq\alpha(b)+|h|.
\]
Otherwise, applying the argument of Step~2
to $\eta_{\varepsilon_0}$, we find
$u\in P\cap N_\rho(K_b)$ with
$\tilde I_b(u)\geq\alpha(b)$.
Thus, by \eqref{n:mass},
\[
\tilde I_{b+h}(u)\geq
\begin{cases}
\alpha(b)-h(M_b+\kappa), & h>0,\\[2mm]
\alpha(b)+|h|(m_b-\kappa), & h<0.
\end{cases}
\]
 Since $m_b-\kappa\le1$, in both cases
\begin{equation}\label{n:low}
\alpha(b+h)\ge\alpha(b)-h(M_b+\kappa)\quad h>0,\qquad \alpha(b+h)\ge\alpha(b)+|h|(m_b-\kappa)\quad h<0.
\end{equation}

\emph{Upper bounds.} 
Let $0<|h|\le\varepsilon_0$ with $b+h\ge0$, let $\varepsilon\in(0,|h|]$, let $P\in\mathbb P$ with $\max_P\tilde I_b\le\alpha(b)+\varepsilon$, and put $\tilde P=\eta_{|h|}(1,P)$. The argument proving \eqref{a10:9} shows that every $u\in\tilde P$ with $\tilde I_b(u)>\alpha(b)-|h|$ lies in $N_\rho(K_b)$. For $u\in\tilde P\setminus N_\rho(K_b)$ we therefore have
\begin{align*}
\tilde I_{b+h}(u)\le\tilde I_b(u)+\max\{-h,0\}\le\alpha(b)-|h|+\max\{-h,0\},
\end{align*}
which equals $\alpha(b)-h$ if $h>0$ and $\alpha(b)$ if $h<0$. For $u\in\tilde P\cap N_\rho(K_b)$, by \eqref{n:mass} and $\tilde I_b(u)\le\alpha(b)+\varepsilon$, we have
\[
\tilde I_{b+h}(u)\leq
\begin{cases}
\alpha(b)+\varepsilon-h(m_b-\kappa),
& h>0,\\[2mm]
\alpha(b)+\varepsilon+|h|(M_b+\kappa),
& h<0.
\end{cases}
\]
 Since $m_b-\kappa\le1$ and $\varepsilon$ is arbitrary,
\begin{equation}\label{n:up}
\alpha(b+h)\le\alpha(b)-h(m_b-\kappa)\quad h>0 ,\qquad \alpha(b+h)\le\alpha(b)+|h|(M_b+\kappa)\quad h<0.
\end{equation}
Dividing \eqref{n:low} and \eqref{n:up} by $h$ gives
\begin{align*}
-(M_b+\kappa)\le\frac{\alpha(b+h)-\alpha(b)}{h}\le-(m_b-\kappa)\qquad(0<|h|\le\varepsilon_0/2,\ b+h\ge0),
\end{align*}
where $\kappa=p\,c_e\rho$ and $\rho$ is arbitrary. Letting first $h\to0$ and then $\rho\to0$ gives \eqref{n:dini}. \qed
\end{proof}

\begin{theorem}\label{n:p2}
Let $p=2$, let $E_2$ be the eigenspace of $\lambda_2$, and recall that, in the notation of Proposition~\ref{n:bpos},
\begin{align*}
M_0=\max\Bigl\{\int_\Gamma(u^+)^2\,dx:\ u\in E_2,\ \|u\|_2=1\Bigr\}.
\end{align*}
Then $\frac12\le M_0<1$, $\alpha'(0^+)=-M_0$, and
\begin{align*}
G'(\lambda_2^{+})=-\frac{M_0}{1-M_0},\qquad G'(\lambda_2^{-})=-\frac{1-M_0}{M_0}.
\end{align*}
In particular, $\mathfrak C$ is differentiable at $(\lambda_2,\lambda_2)$ if and only if $\int_\Gamma(u^+)^2\,dx=\int_\Gamma(u^-)^2\,dx$ for every $u\in E_2$. If $\lambda_2$ is simple, then $M_0=\max\{A,B\}$, and these formulas are those of Theorem~\ref{at10:1}.
\end{theorem}

\begin{proof}
The set $K_0$ of Proposition~\ref{n:bpos} is $E_2\cap M$, and $K_0=-K_0$. Hence $M_0<1$ by Proposition~\ref{n:bpos}, and $M_0\ge\max\{\int_\Gamma(u^+)^2\,dx,\int_\Gamma(u^-)^2\,dx\}\ge\frac12$ for every $u\in K_0$. By Proposition~\ref{n:bpos}, $\liminf_{b\to0^+}(\alpha(b)-\lambda_2)/b\ge-M_0$.

For the upper bound, let $u_*\in K_0$ with $\int_\Gamma(u_*^+)^2\,dx=M_0$. Testing the equations of $\phi_{1,p}$ and $u_*$ with each other gives $\lambda_{1,p}\int_\Gamma\phi_{1,p}u_*\,dx=\int_\Gamma\phi_{1,p}'u_*'\,dx=\lambda_2\int_\Gamma\phi_{1,p}u_*\,dx$, and since $\lambda_{1,p}<\lambda_2$ both integrals vanish. For $\theta\in[-\pi/2,\pi/2]$ put $w(\theta)=\sin\theta\,\phi_{1,p}+\cos\theta\,u_*$. Then $\|w(\theta)\|_2=1$ and $w(\pm\pi/2)=\pm\phi_{1,p}$, so that $w$, reparametrised on $[-1,1]$, belongs to $\mathbb P$, and, since $|\sin\theta|\ge2|\theta|/\pi$,
\begin{align*}
\tilde I_0(w(\theta))=\lambda_{1,p}\sin^2\theta+\lambda_2\cos^2\theta\le\lambda_2-c\,\theta^2,\qquad c=\frac{4(\lambda_2-\lambda_{1,p})}{\pi^2}.
\end{align*}
Moreover $\|w(\theta)-u_*\|_2=2|\sin(\theta/2)|\le|\theta|$, so the first inequality in the proof of Theorem~\ref{at10:1}, with $p=2$ and $\|\cdot\|_2$ in place of $c_e\|\cdot\|$, gives $\int_\Gamma(w(\theta)^+)^2\,dx\ge M_0-2|\theta|$. Hence, for $b>0$,
\begin{align*}
\tilde I_b(w(\theta))\le\lambda_2-c\,\theta^2-b\,(M_0-2|\theta|)\le\lambda_2-bM_0+\max_{s\in\mathbb R}\bigl(2bs-cs^2\bigr)=\lambda_2-bM_0+\frac{b^2}{c},
\end{align*}
and \eqref{a:28} gives
\begin{equation}\label{n:p2up}
\alpha(b)\le\lambda_2-bM_0+\frac{\pi^2b^2}{4(\lambda_2-\lambda_{1,p})}\qquad(b>0).
\end{equation}
Thus $\alpha'(0^+)=-M_0$. The formulas for $G'(\lambda_2^\pm)$ follow as in Step~5 of the proof of Theorem~\ref{at10:1}, with $\max\{A,B\}$ replaced by $M_0$. Finally, $M_0/(1-M_0)=(1-M_0)/M_0$ if and only if $M_0=\frac12$, which, by the lower bound for $M_0$, happens if and only if $\int_\Gamma(u^+)^2\,dx=\int_\Gamma(u^-)^2\,dx$ for every $u\in K_0$. If $\lambda_2$ is simple, then $K_0=\{\pm\phi_{2,p}\}$ and $M_0=\max\{A,B\}$. \qed
\end{proof}

\begin{remark}\label{n:hn}
(i) The lower bound in Theorem~\ref{n:p2} holds for every $p$: by Proposition~\ref{n:bpos}, $\liminf_{b\to0^+}(\alpha(b)-\lambda_2)/b\ge-M_0$, where $M_0$ is the largest value of $\int_\Gamma(u^+)^p\,dx$ on the set $K_0$ of eigenfunctions in $M$ with eigenvalue $\lambda_2$. Equality $\alpha'(0^+)=-M_0$ holds under \eqref{a10:2}, where $K_0=\{\pm\phi_{2,p}\}$ and $M_0=\max\{A,B\}$ (Theorem~\ref{at10:1}), and for $p=2$ (Theorem~\ref{n:p2}). The upper bound in Theorem~\ref{n:p2} uses the plane spanned by $\phi_{1,p}$ and $u_*$, on which $\tilde I_0\le\lambda_2$; for $p\neq2$ the functional $\int_\Gamma|u'|^p\,dx$ is not quadratic, and we know no substitute for this plane.

(ii) Theorem~\ref{n:p2} is consistent with the finite-dimensional results of \cite{holubova2025fucik}. For a symmetric matrix with eigenvalue $\lambda$ and eigenspace $E$, their necessary condition for a Fu\v{c}\'{\i}k curve leaving $(\lambda,\lambda)$ in the direction of a unit vector $u_0\in E$ requires the orthogonal projection of $|u_0|$ onto $E$ to be a multiple of $u_0$. Since $|u_0|=2u_0^+-u_0$, this means that $u_0$ is a critical point of $u\mapsto\|u^+\|^2$ on the unit sphere of $E$, and the tangent line of the curve then has slope $-\|u_0^+\|^2/(1-\|u_0^+\|^2)$. Theorem~\ref{n:p2} identifies the directions followed by the first curve: the maximum $M_0$ of this function on the side $\lambda^+>\lambda^-$, and its minimum $1-M_0$ on the side $\lambda^+<\lambda^-$.
\end{remark}

 We now compute $A$ and $B$ on a three-edge star and show that $A<B$ occurs for every $p$. Let $\Gamma$ be the star with interior vertex $v_0$ and three pendant edges $e_1,e_2,e_3$ of lengths $\ell_1>\ell_2>\ell_3>0$, with Dirichlet conditions at the pendant vertices, that is, $V_D=\{v_1,v_2,v_3\}$, and the Kirchhoff condition at $v_0$; see Figure~\ref{afig:star}. We identify $e_i$ with $[0,\ell_i]$, with $x=0$ at $v_0$ and $x=\ell_i$ at $v_i$, and we let $1<p<\infty$. We use the function $\sin_p$ of Remark~\ref{ar:3}, normalised by $\sin_p(0)=0$, $\sin_p'(0)=1$ and $|\sin_p|^p+|\sin_p'|^p=1$, so that the first Dirichlet eigenvalue on an interval of length $L$ is $(p-1)(\pi_p/L)^p=\frac{p}{p'}(\pi_p/L)^p$, as in Remark~\ref{ar:3}. It is odd, satisfies $\sin_p(\theta+\pi_p)=-\sin_p\theta$, vanishes exactly on $\pi_p\mathbb Z$, and $\cos_p=\sin_p'$ satisfies $(\varphi_p(\cos_p))'=-(p-1)\varphi_p(\sin_p)$, where $\varphi_p(s)=|s|^{p-2}s$; see \cite{lang2011eigenvalues}. We use $\cot_p=\cos_p/\sin_p$ as in Proposition~\ref{aex:100}. For $p=2$ these are $\sin$, $\cos$, $\cot$, and $\pi_2=\pi$.

For $\lambda=(p-1)k^p$ with $k>0$, every solution of $-(\varphi_p(u'))'=\lambda\varphi_p(u)$ on $e_i$ vanishing at $v_i$ has the form
\begin{equation}\label{a10:12}
u_i(x)=c_i\sin_p\bigl(k(\ell_i-x)\bigr),\qquad x\in[0,\ell_i],
\end{equation}
by uniqueness for the initial value problem of this equation \cite{lang2011eigenvalues}, and we put
\begin{equation}\label{a10:13}
S(k)=\varphi_p\bigl(\cot_p(k\ell_1)\bigr)+\varphi_p\bigl(\cot_p(k\ell_2)\bigr)+\varphi_p\bigl(\cot_p(k\ell_3)\bigr).
\end{equation}
For $p=2$, $S(k)=\cot(k\ell_1)+\cot(k\ell_2)+\cot(k\ell_3)$.

\begin{lemma}\label{al10:2}
Let $1<p<\infty$ and $\ell_1>\ell_2>\ell_3>0$.
\begin{enumerate}[label=(\roman*)]
	\item If $u$ is an eigenfunction with eigenvalue $(p-1)k^p$ and $u(v_0)\neq0$, then $\sin_p(k\ell_i)\neq0$ for $i=1,2,3$, $u$ is given by \eqref{a10:12} with $c_i=u(v_0)/\sin_p(k\ell_i)$, and $S(k)=0$. Conversely, if $S(k)=0$ and $\sin_p(k\ell_i)\neq0$ for $i=1,2,3$, then \eqref{a10:12} with $c_i=1/\sin_p(k\ell_i)$ defines an eigenfunction with eigenvalue $(p-1)k^p$ and $u(v_0)=1$.
	\item If $u$ is an eigenfunction with eigenvalue $(p-1)k^p$ and $u(v_0)=0$, then $k\ge\pi_p/\ell_2$.
	\item The function $S$ has exactly one zero $k_1$ in $(0,\pi_p/\ell_1)$ and exactly one zero $k_2$ in
	\begin{align*}
		I_2=\Bigl(\frac{\pi_p}{\ell_1},\ \min\Bigl\{\frac{2\pi_p}{\ell_1},\frac{\pi_p}{\ell_2}\Bigr\}\Bigr),
	\end{align*}
	and $\lambda_{1,p}=(p-1)k_1^p$, $\lambda_2=(p-1)k_2^p$.
	\item The eigenvalue $\lambda_2$ is simple, and its eigenfunction $\phi_{2,p}$ vanishes at exactly one point of $\Gamma\setminus V_D$, namely the point of $e_1$ at distance $\pi_p/k_2$ from $v_1$.
\end{enumerate}
\end{lemma}

\begin{proof}
\textit{(i)} By \eqref{a10:12} and continuity at $v_0$, $c_i\sin_p(k\ell_i)=u(v_0)$ for $i=1,2,3$. The outward derivative of $u$ at $v_0$ along $e_i$ is $u_i'(0)=-c_ik\cos_p(k\ell_i)$, so the Kirchhoff condition reads $\sum_{i=1}^3\varphi_p(c_i)\varphi_p(\cos_p(k\ell_i))=0$. If $u(v_0)\neq0$, the first relation gives $\sin_p(k\ell_i)\neq0$ and $c_i=u(v_0)/\sin_p(k\ell_i)$, and the second becomes $\varphi_p(u(v_0))\,S(k)=0$, that is, $S(k)=0$. Read backwards, the same computation gives the converse: with $c_i=1/\sin_p(k\ell_i)$ the function \eqref{a10:12} is continuous at $v_0$, with $u(v_0)=1$, and it satisfies the Kirchhoff condition there because $S(k)=0$.

\textit{(ii)} If $u(v_0)=0$, then for each $i$ either $c_i=0$ or $\sin_p(k\ell_i)=0$. An edge with $c_i\neq0$ has $\sin_p(k\ell_i)=0$, hence $|\cos_p(k\ell_i)|=1$, and contributes a nonzero term to the Kirchhoff condition. If exactly one edge had $c_i\neq0$, the sum could not vanish. Together with $u\not\equiv0$, this shows that there are two edges $e_i$ and $e_j$, $i\neq j$, with $\sin_p(k\ell_i)=\sin_p(k\ell_j)=0$. Then $k\ell_i,k\ell_j\in\pi_p\mathbb N$, and therefore
\begin{align*}
	k\ \ge\ \frac{\pi_p}{\min\{\ell_i,\ell_j\}}\ \ge\ \frac{\pi_p}{\ell_2},
\end{align*}
because the smaller of any two of the three lengths is at most $\ell_2$.

\textit{(iii)} For $\theta\notin\pi_p\mathbb Z$, the identities $(\varphi_p(\cos_p))'=-(p-1)\varphi_p(\sin_p)$ and $|\sin_p|^p+|\cos_p|^p=1$ give
\begin{align*}
	\frac{d}{d\theta}\,\varphi_p(\cot_p\theta)=\frac{d}{d\theta}\,\frac{\varphi_p(\cos_p\theta)}{\varphi_p(\sin_p\theta)}=-(p-1)\,\frac{|\sin_p\theta|^p+|\cos_p\theta|^p}{|\sin_p\theta|^p}=-\frac{p-1}{|\sin_p\theta|^p}.
\end{align*}
Wherever its three terms are finite, $S$ is therefore strictly decreasing, with
\begin{align*}
	S'(k)=-(p-1)\sum_{i=1}^3\frac{\ell_i}{|\sin_p(k\ell_i)|^p}<0 .
\end{align*}
Moreover, $\varphi_p(\cot_p\theta)\to+\infty$ as $\theta$ approaches $m\pi_p$ from the right, and $\varphi_p(\cot_p\theta)\to-\infty$ as $\theta$ approaches $m\pi_p$ from the left, for every $m\in\mathbb Z$.
On $(0,\pi_p/\ell_1)$ all $k\ell_i$ lie in $(0,\pi_p)$. In particular, $S(k)\to+\infty$ as $k\to0^+$, and $S(k)\to-\infty$ as $k\to(\pi_p/\ell_1)^-$, because the first term tends to $-\infty$ while the other two remain finite. The function $S$ therefore has exactly one zero $k_1$ in this interval. On $I_2$ we have $k\ell_1\in(\pi_p,2\pi_p)$ and $k\ell_2,k\ell_3\in(0,\pi_p)$, so $S$ is finite and strictly decreasing there. It tends to $+\infty$ at the left end of $I_2$, where the first term tends to $+\infty$, and to $-\infty$ at the right end, where the first or the second term tends to $-\infty$. Hence $S$ has exactly one zero $k_2$ in $I_2$. 

We have $k_1\ell_i<\pi_p$ for every $i$, so the eigenfunction associated with $k_1$ by \textit{(i)} is positive, and hence $\lambda_{1,p}=(p-1)k_1^p$ by Lemma~\ref{al:100}\textit{(iii)}. It remains to show that there is no eigenvalue $(p-1)k^p$ with $k_1<k<k_2$. If such an eigenvalue had an eigenfunction with $u(v_0)\neq0$, then $S(k)=0$ by \textit{(i)}, but $S$ has no zero in $(k_1,\pi_p/\ell_1)$ or in $(\pi_p/\ell_1,k_2)$, and $k=\pi_p/\ell_1$ is excluded because $\sin_p(k\ell_1)=0$ there. If it had an eigenfunction with $u(v_0)=0$, then $k\ge\pi_p/\ell_2>k_2$ by \textit{(ii)} and $k_2\in I_2$. Hence $\lambda_2=(p-1)k_2^p$ by \eqref{a:104}.

\textit{(iv)} By \textit{(ii)} and $k_2<\pi_p/\ell_2$, no eigenfunction for $\lambda_2$ vanishes at $v_0$. By \textit{(i)}, every eigenfunction for $\lambda_2$ is then determined by its value at $v_0$, and thus $\lambda_2$ is simple. With $c_1=1$,
\begin{align*}
	\phi_{2,p}(x)=\sin_p\bigl(k_2(\ell_1-x)\bigr)\ \text{ on }e_1,\qquad \phi_{2,p}(x)=\frac{\sin_p(k_2\ell_1)}{\sin_p(k_2\ell_i)}\,\sin_p\bigl(k_2(\ell_i-x)\bigr)\ \text{ on }e_i,\ i=2,3,
\end{align*}
up to normalisation. Since $k_2\ell_1\in(\pi_p,2\pi_p)$, the function $\sin_p(k_2(\ell_1-x))$ vanishes on $(0,\ell_1)$ exactly at $x=\ell_1-\pi_p/k_2$. For $i=2,3$ we have $k_2\ell_i\in(0,\pi_p)$, so the functions $\sin_p(k_2(\ell_i-x))$ do not vanish on $[0,\ell_i)$, and $\phi_{2,p}(v_0)=\sin_p(k_2\ell_1)\neq0$. \qed
\end{proof}

With the normalisation $c_1=1$, the function $\phi_{2,p}$ is positive on the segment of $e_1$ of length $\pi_p/k_2$ adjacent to $v_1$ and negative on the rest of $\Gamma\setminus V_D$, since
\begin{align*}
\phi_{2,p}(v_0)=\sin_p(k_2\ell_1)<0,\qquad c_i=\frac{\sin_p(k_2\ell_1)}{\sin_p(k_2\ell_i)}<0\quad(i=2,3).
\end{align*}
The ratio $r$ in \eqref{a10:3} does not depend on the normalisation of $\phi_{2,p}$, so both weights may be computed with $c_1=1$. Writing $k=k_2$, $\sigma_p=\int_0^{\pi_p}|\sin_p\theta|^p\,d\theta$, and substituting $\theta=k(\ell_i-x)$, we obtain
\begin{align}
\int_\Gamma(\phi_{2,p}^+)^p\,dx&=\frac{\sigma_p}{k},\label{a10:14}\\
\int_\Gamma(\phi_{2,p}^-)^p\,dx&=\frac1k\int_{\pi_p}^{k\ell_1}|\sin_p\theta|^p\,d\theta+\sum_{i=2}^{3}\Bigl|\frac{\sin_p(k\ell_1)}{\sin_p(k\ell_i)}\Bigr|^p\,\frac1k\int_0^{k\ell_i}|\sin_p\theta|^p\,d\theta.\label{a10:15}
\end{align}
Here \eqref{a10:14} is the integral over the segment on which $\phi_{2,p}>0$. In \eqref{a10:15} the first term comes from the rest of $e_1$, and the sum from $e_2$ and $e_3$. Each of the three contributions to \eqref{a10:15} is nonnegative. For $p=2$ one has $\sigma_2=\pi/2$, and \eqref{a10:14}--\eqref{a10:15} read
\begin{equation}\label{a:105}
\int_\Gamma(\phi_{2,p}^+)^2\,dx=\frac{\pi}{2k},\qquad
\int_\Gamma(\phi_{2,p}^-)^2\,dx=\frac{\ell_1-\pi/k}{2}-\frac{\sin(2k\ell_1)}{4k}+\sum_{i=2}^{3}\frac{\sin^2(k\ell_1)}{\sin^2(k\ell_i)}\Bigl(\frac{\ell_i}{2}-\frac{\sin(2k\ell_i)}{4k}\Bigr).
\end{equation}

\begin{prop}\label{ap10:1}
Let $1<p<\infty$ and $\Omega=\{(\ell_1,\ell_2,\ell_3):\ell_1>\ell_2>\ell_3>0\}$. For $\ell\in\Omega$ let $D(\ell)$ be the right-hand side of \eqref{a10:15} minus the right-hand side of \eqref{a10:14}, with $k=k_2$.
\begin{enumerate}[label=(\roman*)]
\item $\mathfrak C$ is differentiable at $(\lambda_2,\lambda_2)$ if and only if $D(\ell)=0$, that is, if and only if
\begin{equation}\label{a10:16}
\int_\Gamma(\phi_{2,p}^+)^p\,dx=\int_\Gamma(\phi_{2,p}^-)^p\,dx .
\end{equation}
Moreover $r<1$ if $D(\ell)>0$, and $r>1$ if $D(\ell)<0$.
\item $D$ is continuous on $\Omega$. Hence $\Omega_+=\{D>0\}$ and $\Omega_-=\{D<0\}$ are open, and $\mathfrak C$ has a corner at $(\lambda_2,\lambda_2)$ for every $\ell\in\Omega_+\cup\Omega_-$.
\item If $\ell_2<\ell_1<2\ell_2$, then $(\ell_1,\ell_2,\ell_3)\in\Omega_+$ for every sufficiently small $\ell_3>0$. In particular, for every $p$ there are stars on which $\mathfrak C$ has a corner at $(\lambda_2,\lambda_2)$.
\item For $p=2$, the set $\Omega_0=\{D=0\}$ is a proper real-analytic subvariety of $\Omega$, hence of measure zero, and $\mathfrak C$ has a corner at $(\lambda_2,\lambda_2)$ for almost every $\ell\in\Omega$.

\end{enumerate}
\end{prop}

\begin{proof}
\textit{(i)} By Lemma~\ref{al10:2}, $\lambda_2$ is simple, so \eqref{a10:2} holds and Theorem~\ref{at10:1} applies. The normalisation $c_1=1$ used in \eqref{a10:14}--\eqref{a10:15} multiplies $A$ and $B$ by the same positive factor, so $D$ has the sign of $B-A$, and $D=0$ if and only if $r=1$.

\textit{(ii)} The set $\{(k,\ell):\ell\in\Omega,\ k\in I_2(\ell)\}$ is open, because the endpoints of $I_2$ depend continuously on $\ell$, and on it $S$ is continuous and strictly decreasing in $k$. Given $\ell^0\in\Omega$ and a small $\eta>0$, we have $S(k_2(\ell^0)-\eta)>0>S(k_2(\ell^0)+\eta)$ at $\ell=\ell^0$, and these inequalities persist for $\ell$ near $\ell^0$, so that $|k_2(\ell)-k_2(\ell^0)|<\eta$. Thus $k_2$ is continuous on $\Omega$, and, as $\sin_p(k_2\ell_2)$ and $\sin_p(k_2\ell_3)$ do not vanish, so is $D$. The last assertion follows from \textit{(i)}.

\textit{(iii)} Fix $\ell_1,\ell_2$ with $\ell_2<\ell_1<2\ell_2$ and let $\ell_3\to0^+$. Then $I_2=(\pi_p/\ell_1,\pi_p/\ell_2)$, so $k_2$ stays in the compact interval $[\pi_p/\ell_1,\pi_p/\ell_2]$, $k_2\ell_3\to0^+$, and $\varphi_p(\cot_p(k_2\ell_3))\to+\infty$. Since $S(k_2)=0$,
\begin{align*}
	\varphi_p\bigl(\cot_p(k_2\ell_1)\bigr)+\varphi_p\bigl(\cot_p(k_2\ell_2)\bigr)=-\varphi_p\bigl(\cot_p(k_2\ell_3)\bigr)\to-\infty .
\end{align*}
The first term is bounded below: $k_2\ell_1\in(\pi_p,\pi_p\ell_1/\ell_2)$ with $\pi_p\ell_1/\ell_2<2\pi_p$, and $\varphi_p\circ\cot_p$ is decreasing on $(\pi_p,2\pi_p)$, so that $\varphi_p(\cot_p(k_2\ell_1))>\varphi_p(\cot_p(\pi_p\ell_1/\ell_2))$. Hence $\varphi_p(\cot_p(k_2\ell_2))\to-\infty$ with $k_2\ell_2\in(0,\pi_p)$, that is, $k_2\ell_2\to\pi_p^-$, so that
\begin{align*}
	\sin_p(k_2\ell_2)\to0,\qquad \sin_p(k_2\ell_1)\to\sin_p\Bigl(\frac{\pi_p\ell_1}{\ell_2}\Bigr)\neq0,\qquad \frac1{k_2}\int_0^{k_2\ell_2}|\sin_p\theta|^p\,d\theta\to\frac{\sigma_p\,\ell_2}{\pi_p}>0,
\end{align*}
the middle limit being nonzero because $\pi_p<\pi_p\ell_1/\ell_2<2\pi_p$. The term $i=2$ in \eqref{a10:15} therefore tends to $+\infty$, the remaining terms of \eqref{a10:15} are nonnegative, and \eqref{a10:14} satisfies $\sigma_p/k_2\le\sigma_p\ell_1/\pi_p$. Hence $D\to+\infty$ as $\ell_3\to0^+$.

\textit{(iv)} Let $p=2$. Since $S'<0$ on $I_2$, the zero $k_2$ is simple, and by the implicit function theorem $k_2$ depends real-analytically on $\ell\in\Omega$. As $\sin(k_2\ell_2)$ and $\sin(k_2\ell_3)$ do not vanish, the function $D$, given by \eqref{a:105}, is real-analytic on the connected open set $\Omega$, and it is not identically zero by \textit{(iii)}. Therefore $\Omega_0$ is a proper real-analytic subvariety of $\Omega$ and has measure zero. For $p\neq2$, the function $\sin_p$ is in general not real-analytic at $\pi_p/2$, which the numbers $k_2\ell_2$ and $k_2\ell_3$ may cross, and this argument does not apply. \qed

\end{proof}

\begin{ex}\label{n:eqstar}
Let $\Gamma$ be the star of Figure~\ref{afig:star} with $\ell_1=\ell_2=\ell_3=\ell$, Dirichlet conditions at $v_1,v_2,v_3$ and the Kirchhoff condition at $v_0$, and let $p=2$. The proofs of Lemma~\ref{al10:2}(i) and (ii) do not use the ordering of the lengths. By (i), an eigenvalue $k^2$ with an eigenfunction that does not vanish at $v_0$ satisfies $3\cot(k\ell)=0$, that is, $k\ell\in\frac{\pi}2+\pi\mathbb Z$. By the proof of (ii), an eigenfunction $u$ with $u(v_0)=0$ requires $k\ell\in\pi\mathbb N$. Then $|\cos(k\ell)|=1$, and the Kirchhoff condition at $v_0$ reads $c_1+c_2+c_3=0$. Hence $\lambda_{1,2}=(\pi/(2\ell))^2$, $\lambda_2=(\pi/\ell)^2$, and
\begin{align*}
E_2=\Bigl\{u=c_i\sin\Bigl(\frac{\pi x}{\ell}\Bigr)\ \text{on }e_i\ (i=1,2,3):\ c_1+c_2+c_3=0\Bigr\}
\end{align*}
has dimension two, so \eqref{a10:2} fails. As $\sin(\pi x/\ell)>0$ on $(0,\ell)$, for $u\in E_2\setminus\{0\}$ the quotient $\int_\Gamma(u^+)^2\,dx/\int_\Gamma u^2\,dx$ equals $\sum_{c_i>0}c_i^2/\sum_ic_i^2$. If exactly one $c_i$, say $c_1$, is positive, put $s=|c_2|$ and $t=|c_3|$. Then $c_1=s+t$, and the quotient equals
\begin{align*}
\frac{(s+t)^2}{(s+t)^2+s^2+t^2}\le\frac23,
\end{align*}
with equality if and only if $s=t$. If two of the $c_i$ are positive, the quotient is at most $\frac12$, because $(s+t)^2\ge s^2+t^2$. Theorem~\ref{n:p2} therefore applies with $M_0=\frac23$, attained at $c=(2,-1,-1)$, and
\begin{align*}
G'(\lambda_2^{+})=-2,\qquad G'(\lambda_2^{-})=-\frac12 .
\end{align*}
Thus $\mathfrak C$ has a corner at $(\lambda_2,\lambda_2)$, although $\Gamma$ is invariant under every permutation of its edges. The space $E_2$ contains functions with $\int_\Gamma(u^+)^2\,dx=\int_\Gamma(u^-)^2\,dx$, for instance $c=(1,-1,0)$, which is odd under the exchange of $e_1$ and $e_2$, but the slopes are determined by the maximizer $c=(2,-1,-1)$. Corollary~\ref{ac10:1} does not apply because \eqref{a10:2} fails. The maximizer vanishes at $v_0$, so its negative set is the union of the interiors of $e_2$ and $e_3$, and it has three nodal domains. Thus, on a graph, a solution associated with a point of $\mathfrak C$ may have more than two nodal domains, in contrast with the nodal domain property proved for domains in \cite{cuesta2000nodal}.
\end{ex}

\begin{remark}\label{ar10:3}
The proofs of Theorem~\ref{at10:1} and Corollary~\ref{ac10:1} use the graph only through Theorem~\ref{at:3}, Lemma~\ref{al:100}, Lemma~\ref{al:3}, Theorem~\ref{at:1} and Corollary~\ref{ar:2}, and they apply without change to the Dirichlet $p$-Laplacian on a bounded domain $\Omega\subset\mathbb R^N$ \cite{cuesta1999beginning}. On a domain, however, \eqref{a10:2} is hard to verify, and it fails on the most symmetric ones. On a ball $B\subset\mathbb R^N$, $N\ge2$, for every $p\in(1,\infty)$ the second eigenfunctions are not radial \cite{anoop2016structure}. If $u$ is one of them and $u\circ R=c(R)\,u$ for every $R$ in the rotation group $SO(N)$, then $c$ is continuous on $SO(N)$ with values in $\{\pm1\}$ and $c(\mathrm{id})=1$. As $SO(N)$ is connected, $c\equiv1$, and $u$ is radial because $SO(N)$ acts transitively on spheres. Hence some rotation of $u$ is an eigenfunction that is not a multiple of $u$, and \eqref{a10:2} fails. On a graph \eqref{a10:2} can be checked and $r$ computed, and both alternatives of Theorem~\ref{at10:1} occur for every $p$ (Corollary~\ref{ac10:1} and Proposition~\ref{ap10:1}). Without \eqref{a10:2}, a symmetric graph may still give a corner (Example~\ref{n:eqstar}).
\end{remark}

We close this section with three questions. First, for $p\neq2$ we do not know whether $\alpha'(0^+)=-M_0$ when \eqref{a10:2} fails. Second, for $p=2$ and $\lambda_2$ simple the corner of $\mathfrak C$ is the crossing of two smooth curves of $\Sigma_{2,\Gamma}$ (Remark~\ref{ar10:2}). For $p\neq2$ we do not know whether $\Sigma_{p,\Gamma}$ contains a second curve through $(\lambda_2,\lambda_2)$, with slope $-\min\{r,1/r\}$ on the side $\lambda^+>\lambda^-$. Third, for $b>0$ Proposition~\ref{n:bpos} gives only the bounds \eqref{n:dini}, and we do not know whether the one-sided derivatives of $\alpha$ exist. The proof of Theorem~\ref{at10:1} routes a near-optimal path through a neighborhood of the critical point of larger mass (Steps~3 and~4), and this uses Lemma~\ref{al10:1} and the evenness of $\tilde I_0$. For $b>0$ the functional $\tilde I_b$ is not even, and Lemma~\ref{al10:1} has no analogue.

\section{Nonresonance between \texorpdfstring{$(\lambda_{1,p}(\Gamma,V_D), \lambda_{1,p}(\Gamma,V_D))$}{(lambda1,p(Gamma,VD), lambda1,p(Gamma,VD))} and \texorpdfstring{$\mathfrak{C}$}{C}}\label{as:nonres}

This section deals with the problem
\begin{equation}{\label{a:57}}
\begin{cases}
-\bigl(|u'|^{p-2}u'\bigr)'=f(x,u) \quad &\text{on the edges of $\Gamma$, } \\[6pt]
u(v) = 0 & \forall v \in V_D, \\[8pt]
\displaystyle\sum_{e \in E_v(\Gamma)} |u_e'(v)|^{p-2}u_e'(v)=0, &v\in V(\Gamma)\setminus V_D.
\end{cases}
\end{equation}

Throughout this section $V_D\subseteq V(\Gamma)$ satisfies \eqref{ah:100} and may be empty. In the latter case $\lambda_{1,p}(\Gamma,V_D)=0$, the first eigenfunction $\phi_{1,p}$ is constant, and the results of Sections \ref{as:4}--\ref{as:6} quoted below hold by Proposition~\ref{ap:100}. Proposition~\ref{ap:7}, the only one that needs $V_D\neq\emptyset$, is not used.

We assume that $f(x,u)/ |u|^{p-2}u$ lies asymptotically between the point $(\lambda_{1,p},\lambda_{1,p})$, where $\lambda_{1,p}=\lambda_{1,p}(\Gamma,V_D)$, and a point $(\lambda^+,\lambda^-)\in\mathfrak{C}$.  Let $f : \Gamma \times \mathbb{R}\to \mathbb{R}$ satisfy the $L^{\infty}$-Carath\'eodory conditions. Given a point $(\lambda^+, \lambda^-) \in \mathfrak{C}$, we assume that, uniformly with respect to $x$,
\begin{equation}{\label{a:58}}
\pi_{\pm}(x) \le \liminf_{s \to \pm \infty} \frac{f(x,s)}{|s|^{p-2}s} \le \limsup_{s \to \pm \infty} \frac{f(x,s)}{|s|^{p-2}s} \le \Pi_{\pm}(x)
\end{equation}
where $\pi_{\pm}(x)$ and $\Pi_{\pm}(x)$ are $L^{\infty}$ functions such that
\begin{equation}{\label{a:59}}
\begin{cases}
\lambda_{1,p}(\Gamma, V_D) \le \pi_{+}(x) \le \Pi_{+}(x) \le \lambda^+ \quad \text{ a.e. in } \Gamma\\[6pt]
\lambda_{1,p}(\Gamma, V_D) \le \pi_{-}(x) \le \Pi_{-}(x) \le \lambda^- \quad \text{ a.e. in } \Gamma.
\end{cases}
\end{equation}
Put $F(x,s) = \int_0^s f(x,t)\,dt$. We also assume that, uniformly with respect to $x$,
\begin{equation}{\label{a:60}}
\theta_{\pm}(x) \le \liminf_{s \to \pm \infty} \frac{pF(x,s)}{|s|^{p}} \le \limsup_{s \to \pm \infty} \frac{pF(x,s)}{|s|^{p}} \le \Theta_{\pm}(x)
\end{equation}

where $\theta_{\pm}(x)$ and $\Theta_{\pm}(x)$ are $L^\infty$ functions such that
\begin{equation}{\label{a:61}}
\begin{cases}
\lambda_{1,p}(\Gamma, V_D) \le \theta_{+}(x) \le \Theta_{+}(x) \le \lambda^+ \quad \text{ a.e. in } \Gamma\\[6pt]
\lambda_{1,p}(\Gamma, V_D) \le \theta_{-}(x) \le \Theta_{-}(x) \le \lambda^- \quad \text{ a.e. in } \Gamma \\[6pt]
\theta_{+} (x) > \lambda_{1,p}(\Gamma, V_D) \quad \text{ and } \quad \theta_{-}(x) > \lambda_{1,p}(\Gamma, V_D) \text{ on subsets of positive measure, } \\[6pt]
\text{either } \Theta_{+}(x) <\lambda^+ \text{ a.e. in } \Gamma  \text{ or }  \Theta_{-}(x) <\lambda^- \text{ a.e. in } \Gamma .
\end{cases}
\end{equation}

\begin{theorem}{\label{at:5}}
Let \eqref{a:58}, \eqref{a:59}, \eqref{a:60}, and \eqref{a:61} hold and let $(\lambda^+,\lambda^-) \in \mathfrak{C}$. Then problem \eqref{a:57} admits at least one solution $u \in X(\Gamma, V_D)$.
\end{theorem}

For the proof of Theorem \ref{at:5}, define the energy functional $\mathcal{E} :X(\Gamma, V_D) \to \mathbb{R}$ by
\begin{equation}\label{a:62}
\mathcal{E}(u) = \frac{1}{p} \int_\Gamma |u'|^p \, dx - \int_\Gamma F(x,u) \, dx.
\end{equation}
Then $\mathcal{E}\in C^1(X(\Gamma,V_D),\mathbb{R})$, and its critical points are the weak solutions of \eqref{a:57}.

By \eqref{a:60}, which holds uniformly with respect to $x$, and the $L^\infty$-Carath\'eodory conditions on $f$, for every $\delta>0$ there is $a_\delta\in L^1(\Gamma)$ such that, for a.e.\ $x\in\Gamma$,
\begin{equation}{\label{a:71}}
\begin{cases}
	(\theta_+(x)-\delta)\frac{s^p}{p} - a_\delta(x) \le F(x,s) \le (\Theta_+(x)+\delta)\frac{s^p}{p} + a_\delta(x)
	\quad \forall s>0,\\[6pt]
	(\theta_-(x)-\delta)\frac{|s|^p}{p} - a_\delta(x) \le F(x,s) \le (\Theta_-(x)+\delta)\frac{|s|^p}{p} + a_\delta(x) \quad \forall s<0.
\end{cases}
\end{equation}

\begin{lemma}{\label{al:10}}
$\mathcal{E}$ satisfies the $(P.S.)$ condition on $X(\Gamma,V_D)$.
\end{lemma}

\begin{proof}
Let $(u_n)$ be a $(P.S.)$ sequence, that is,
\begin{align}
|\mathcal{E}(u_n)| &\le K, \label{a:64} \\
|\langle\mathcal{E}'(u_n), \phi\rangle| &\le \epsilon_n \|\phi\|\qquad \forall \phi \in X(\Gamma, V_D), \label{a:65}
\end{align}
with a constant $K$ and $\epsilon_n\to0$. By \eqref{a:58} and the $L^\infty$-Carath\'eodory conditions, there are $K_1>0$ and $a\in L^\infty(\Gamma)$ with
\begin{equation}\label{n:growth}
|f(x,s)|\le a(x)+K_1|s|^{p-1}\qquad\text{for a.e. }x\in\Gamma\text{ and all }s\in\mathbb R .
\end{equation}
It suffices to prove that $(u_n)$ is bounded in $X(\Gamma,V_D)$. If it is, then along a subsequence $u_n\rightharpoonup u$ in $X(\Gamma,V_D)$ and $u_n\to u$ in $L^p(\Gamma)$ by Theorem~\ref{at:3}. By \eqref{n:growth}, $(f(x,u_n))$ is bounded in $L^{p'}(\Gamma)$, so that $\int_\Gamma f(x,u_n)(u_n-u)\,dx\to0$, and \eqref{a:65} with $\phi=u_n-u$ gives $\int_\Gamma|u_n'|^{p-2}u_n'(u_n'-u')\,dx\to0$. Hence $u_n\to u$ in $X(\Gamma,V_D)$ by Lemma~\ref{n:splus}.

Suppose, for contradiction, that $\|u_n\|\to\infty$ along a subsequence, still denoted $(u_n)$, and put $w_n=u_n/\|u_n\|$. Passing to a further subsequence, $w_n\rightharpoonup w_0$ in $X(\Gamma,V_D)$, $w_n\to w_0$ in $L^p(\Gamma)$ and $w_n\to w_0$ a.e.\ in $\Gamma$. Put $g_n=f(x,u_n)/\|u_n\|^{p-1}$. By \eqref{n:growth},
\begin{equation}\label{n:gn}
|g_n|\le\frac{a}{\|u_n\|^{p-1}}+K_1|w_n|^{p-1},
\end{equation}
so $(g_n)$ is bounded in $L^{p'}(\Gamma)$ and, along a subsequence, $g_n\rightharpoonup f_0$ in $L^{p'}(\Gamma)$. Testing \eqref{a:65} with $w_n-w_0$ and dividing by $\|u_n\|^{p-1}$, we obtain $\int_\Gamma|w_n'|^{p-2}w_n'(w_n'-w_0')\,dx\to0$, because $\int_\Gamma g_n(w_n-w_0)\,dx\to0$ by H\"older's inequality. Hence $w_n\to w_0$ in $X(\Gamma,V_D)$ by Lemma~\ref{n:splus}, in particular $\|w_0\|=1$, and dividing \eqref{a:65} by $\|u_n\|^{p-1}$ and letting $n\to\infty$ gives
\begin{equation}\label{n:w0eq}
\int_\Gamma|w_0'|^{p-2}w_0'\phi'\,dx=\int_\Gamma f_0\,\phi\,dx\qquad\forall\phi\in X(\Gamma,V_D).
\end{equation}

\emph{Representation of $f_0$.} We claim that $f_0=0$ a.e.\ on $\{w_0=0\}$, and that
\begin{equation}\label{n:f0}
\pi_+w_0^{p-1}\le f_0\le\Pi_+w_0^{p-1}\ \text{ a.e. on }\{w_0>0\},\qquad \pi_-|w_0|^{p-1}\le-f_0\le\Pi_-|w_0|^{p-1}\ \text{ a.e. on }\{w_0<0\}.
\end{equation}
As $|w_n|^{p-1}\to|w_0|^{p-1}$ in $L^{p'}(\Gamma)$, \eqref{n:gn} gives $g_n\chi_{\{w_0=0\}}\to0$ in $L^{p'}(\Gamma)$, so $f_0=0$ a.e.\ on $\{w_0=0\}$. Let $\delta>0$. By \eqref{a:58} there is $S_\delta>0$ with $f(x,s)\le(\Pi_+(x)+\delta)s^{p-1}$ for a.e.\ $x\in\Gamma$ and all $s\ge S_\delta$. For a.e.\ $x$ with $w_0(x)>0$ we have $u_n(x)\to+\infty$, and therefore
\begin{align*}
h_n=\bigl(g_n-(\Pi_++\delta)(w_n^+)^{p-1}\bigr)^+\chi_{\{w_0>0\}}\to0\qquad\text{a.e. in }\Gamma .
\end{align*}
By \eqref{n:gn}, $0\le h_n\le G_n:=a\|u_n\|^{1-p}+K_1|w_n|^{p-1}$, and $G_n$ converges in $L^{p'}(\Gamma)$, so $h_n\to0$ in $L^{p'}(\Gamma)$ by the generalised dominated convergence theorem. The functions $\bigl(g_n-(\Pi_++\delta)(w_n^+)^{p-1}\bigr)\chi_{\{w_0>0\}}$ are bounded above by $h_n$ and converge weakly in $L^{p'}(\Gamma)$ to $\bigl(f_0-(\Pi_++\delta)w_0^{p-1}\bigr)\chi_{\{w_0>0\}}$, so $f_0\le(\Pi_++\delta)w_0^{p-1}$ a.e.\ on $\{w_0>0\}$. The other three bounds follow in the same way, with the dominating functions $G_n+(\|\pi_\pm\|_\infty+\delta)|w_n|^{p-1}$, which also converge in $L^{p'}(\Gamma)$. Letting $\delta\to0$ we obtain \eqref{n:f0}.

Put $\lambda^+(x)=f_0(x)/w_0(x)^{p-1}$ on $\{w_0>0\}$ and $\lambda^+(x)=\lambda^+$ elsewhere, and $\lambda^-(x)=-f_0(x)/|w_0(x)|^{p-1}$ on $\{w_0<0\}$ and $\lambda^-(x)=\lambda^-$ elsewhere. Then
\begin{align*}
f_0=\lambda^+(x)(w_0^+)^{p-1}-\lambda^-(x)(w_0^-)^{p-1},
\end{align*}
the functions $\lambda^\pm(x)$ satisfy \eqref{a:43} by \eqref{n:f0} and \eqref{a:59}, and since $\lambda^\pm>\lambda_{1,p}(\Gamma,V_D)$ by Proposition~\ref{ap:3},
\begin{equation}{\label{a:68}}
\lambda^+(x)> \lambda_{1,p}(\Gamma,V_D) \text{ on } \{w_0\le 0\}, \qquad \lambda^-(x)> \lambda_{1,p}(\Gamma,V_D) \text{ on } \{w_0\ge 0\}.
\end{equation}
By \eqref{n:w0eq}, $w_0$ is a nontrivial solution of \eqref{a:45}. If \eqref{a:44} holds, Lemma~\ref{al:9} shows that $w_0$ changes sign and that $\lambda^\pm(x)=\lambda^\pm$ a.e.\ on $\{\pm w_0>0\}$. Hence one of the following alternatives holds:
\begin{enumerate}[label=(\roman*)]
\item $\lambda^+(x) = \lambda_{1,p}(\Gamma,V_D)$ a.e.\ in $\Gamma$;
\item $\lambda^-(x) = \lambda_{1,p}(\Gamma,V_D)$ a.e.\ in $\Gamma$;
\item $w_0$ changes sign and is a solution of \eqref{a1:1} with the point $(\lambda^+,\lambda^-)\in \mathfrak{C}$.
\end{enumerate}

We show that each of them leads to a contradiction. Dividing \eqref{a:64} by $\|u_n\|^p$ and using $w_n\to w_0$ in $X(\Gamma,V_D)$, we get
\begin{equation}{\label{a7:70}}
\lim_{n\to\infty} \int_\Gamma \frac{pF(x,u_n)}{\|u_n\|^p} \, dx= \int_\Gamma |w_0'|^p \, dx.
\end{equation}
On the other hand, by \eqref{a:71}, for every $\delta>0$, a.e.\ $x\in\Gamma$ and all $s\in\mathbb R$,
\begin{align*}
(\theta_+-\delta)(s^+)^p+(\theta_--\delta)(s^-)^p-p\,a_\delta\le pF(x,s)\le(\Theta_++\delta)(s^+)^p+(\Theta_-+\delta)(s^-)^p+p\,a_\delta .
\end{align*}
We take $s=u_n(x)$, divide by $\|u_n\|^p$ and integrate. Since $w_n^\pm\to w_0^\pm$ in $L^p(\Gamma)$, letting $n\to\infty$ and then $\delta\to0$ we obtain from \eqref{a7:70}
\begin{equation}\label{n:lims}
\int_\Gamma\bigl(\theta_+(w_0^+)^p+\theta_-(w_0^-)^p\bigr)\,dx\le\int_\Gamma|w_0'|^p\,dx\le\int_\Gamma\bigl(\Theta_+(w_0^+)^p+\Theta_-(w_0^-)^p\bigr)\,dx .
\end{equation}
Suppose (i) holds. By \eqref{a:68}, $w_0>0$ a.e.\ in $\Gamma$, and \eqref{n:w0eq} with $\phi=w_0$ gives
\begin{equation}{\label{a7:69}}
\int_\Gamma |w_0'|^p \,dx = \lambda_{1,p}(\Gamma,V_D)\int_\Gamma w_0^p \, dx.
\end{equation}
The left inequality in \eqref{n:lims} then gives $\int_\Gamma(\theta_+-\lambda_{1,p}(\Gamma,V_D))w_0^p\,dx\le0$. This is impossible, because $\theta_+\ge\lambda_{1,p}(\Gamma,V_D)$ a.e.\ and $\theta_+>\lambda_{1,p}(\Gamma,V_D)$ on a set of positive measure by \eqref{a:61}, while $w_0>0$ a.e. Case (ii) leads to a contradiction in the same way, with $w_0<0$ a.e.\ and $\theta_-$. Suppose (iii) holds. Then \eqref{n:w0eq} with $\phi=w_0$ and the right inequality in \eqref{n:lims} give
\begin{align*}
\int_\Gamma\bigl((\lambda^+-\Theta_+)(w_0^+)^p+(\lambda^--\Theta_-)(w_0^-)^p\bigr)\,dx\le0 .
\end{align*}
By \eqref{a:61} both integrands are nonnegative, and either $\Theta_+<\lambda^+$ a.e.\ or $\Theta_-<\lambda^-$ a.e. Both $w_0^+$ and $w_0^-$ are positive on sets of positive measure, so the integral is positive, a contradiction. Hence $(u_n)$ is bounded. \qed
\end{proof}

\begin{lemma}{\label{al:11}}
There exist $\delta>0$ and $R>0$ such that
\begin{equation}{\label{a:70}}
\max \{\mathcal{E}(R\phi_{1,p}), \mathcal{E}(-R\phi_{1,p})\} < -\|a_\delta\|_1\le\inf_{P\in\mathbb P_1}\max_{u \in P[-1,+1]} \mathcal{E}(u),
\end{equation}
where $a_\delta$ is given by \eqref{a:71} and $\mathbb{P}_1 = \{P \in C ([-1,+1], X(\Gamma,V_D)) : P(\pm1) = \pm R\phi_{1,p}\}$.
\end{lemma}

\begin{proof}
We use the functions $a_\delta$ of \eqref{a:71}. The functional associated with the functions $\Theta_\pm(x)$ from \eqref{a:71} is
\begin{equation}{\label{a:72}}
\mathcal{L}(u) = \int_\Gamma |u'|^p \, dx - \int_\Gamma \Theta_+(x) (u^+)^p \, dx - \int_\Gamma \Theta_-(x) (u^-)^p \, dx,
\end{equation}
and we claim that
\begin{equation}{\label{a:73}}
\beta  = \inf_{P \in \mathbb{P}} \max_{u \in P[-1,+1]} \mathcal{L}(u) >0,
\end{equation}
where $\mathbb{P}$ is, as before, the set of all continuous paths from $-\phi_{1,p}$ to $\phi_{1,p}$ in $M$. Set $b = \lambda^+ - \lambda^- \ge 0$. The case $\lambda^+<\lambda^-$ can be reduced to this one by replacing $u$ with $-u$ and $f(x,s)$ with $-f(x,-s)$. This transformation interchanges the subscripts $+$ and $-$ in \eqref{a:58}--\eqref{a:61}, while the symmetry of $\mathfrak{C}$ yields $(\lambda^-,\lambda^+)\in \mathfrak{C}$.
 On $M$, \eqref{a:61} gives $\mathcal L(u)\ge\int_\Gamma|u'|^p\,dx-\lambda^+\int_\Gamma(u^+)^p\,dx-\lambda^-\int_\Gamma(u^-)^p\,dx=\tilde I_b(u)-\lambda^-$, and $\max_P\tilde I_b\ge\alpha(b)=\lambda^-$ for every $P\in\mathbb P$, because $(\lambda^+,\lambda^-)\in\mathfrak C$. Hence $\beta\ge0$. From $\theta_{\pm}(x) \le \Theta_{\pm}(x)$, the definition \eqref{a:72} and \eqref{a:61} we deduce that
\begin{align*}
\mathcal{L}(\pm \phi_{1,p})
&\le \int_\Gamma (\lambda_{1,p}(\Gamma,V_D) - \theta_\pm (x)) \phi_{1,p}^p \, dx <0.
\end{align*}
The restriction $\tilde{\mathcal{L}}$ of $\mathcal{L}$ to $M$ therefore has the mountain pass geometry
\begin{align*}
\max\{\tilde{\mathcal{L}}(+\phi_{1,p}), \tilde{\mathcal{L}}(-\phi_{1,p})\} < 0 \le  \max_{u \in P[-1,+1]} \tilde{\mathcal{L}}(u),
\end{align*}
for every path $P \in \mathbb{P}$. The $(P.S.)$ condition for $\tilde{\mathcal L}$ on $M$ is proved exactly as in Lemma \ref{al:3}: the weights $\Theta_\pm\in L^\infty(\Gamma)$ enter only through terms that converge by the compact embedding of Theorem~\ref{at:3}. Hence $\beta $ is a critical value of $\tilde{\mathcal{L}}$, that is, there exist $u \in M$ and $t \in \mathbb{R}$ with
\begin{equation}{\label{a7:73}}
\begin{cases}
	\mathcal{L}(u) = \beta, \\[5pt]
	(\mathcal{L'}(u), \phi) = t(J'(u), \phi), \qquad \forall \phi \in X(\Gamma,V_D).
\end{cases}
\end{equation}
Suppose, for contradiction, that $\beta = 0$. Taking $\phi = u$ in \eqref{a7:73} shows $t=0$, so $u$ is a nontrivial solution of \eqref{a:45} with $\lambda^\pm(x)=\Theta_\pm(x)$. By \eqref{a:61}, this contradicts Lemma~\ref{al:9}. We turn to \eqref{a:70}. By \eqref{a:61} and Lemma~\ref{al:100}, $\int_\Gamma(\theta_\pm-\lambda_{1,p}(\Gamma,V_D))\phi_{1,p}^p\,dx>0$, and we fix $\delta>0$ smaller than $\beta$ and than both integrals. The left-hand inequalities in \eqref{a:71} give
\begin{align*}
\mathcal{E}(\pm R \phi_{1,p}) \le \frac{R^p}{p} \int_\Gamma \bigl(\lambda_{1,p}(\Gamma,V_D) - \theta_\pm +\delta\bigr)\phi_{1,p}^p \, dx + \|a_\delta\|_1\to-\infty\qquad(R\to\infty),
\end{align*}
and we choose $R$ with
\begin{equation}{\label{a:74}}
\mathcal{E}(\pm R \phi_{1,p}) < -\|a_\delta\|_1.
\end{equation} Let $P \in \mathbb{P}_1$. If $0 \in P[-1,+1]$, then by \eqref{a:74}
\begin{align*}
\mathcal{E}(\pm R \phi_{1,p}) < -\|a_\delta\|_1 \le 0 = \mathcal{E}(0) \le \max_{u \in P[-1,+1]} \mathcal{E}(u).
\end{align*}
In particular $\max_{u\in P[-1,+1]}\mathcal E(u)\ge-\|a_\delta\|_1$. If $0 \notin P[-1,+1]$, the normalized path $\tilde{P}(t) = P(t)/ \|P(t)\|_p$ belongs to $\mathbb{P}$. By \eqref{a:71},
\begin{align*}
\mathcal{E}(u) \ge \frac{\mathcal{L}(u)- \delta\|u\|_p^p}{p} - \|a_\delta\|_1.
\end{align*}
Therefore, by \eqref{a:73},
\begin{align*}
\max_{u\in P[-1,+1]}\frac{p\mathcal{E}(u)+ \delta\|u\|_p^p + p\|a_\delta\|_1}{\|u\|_p^p} \ge \max_{u\in \tilde{P}[-1,+1]} \mathcal{L}(u) \ge \beta,
\end{align*}
and consequently, since $\delta<\beta$,
\begin{align*}
\max_{u\in P[-1,+1]}\frac{p\mathcal{E}(u)+  p\|a_\delta\|_1}{\|u\|_p^p} \ge \beta - \delta >0,
\end{align*}
which implies
\begin{align*}
\max_{u \in P[-1,+1]} \mathcal{E}(u) > -\|a_\delta\|_1 .
\end{align*}
Thus $\max_P\mathcal E\ge-\|a_\delta\|_1$ for every $P\in\mathbb P_1$, and together with \eqref{a:74} this is \eqref{a:70}. \qed
\end{proof}

\textbf{Proof of Theorem \ref{at:5}.} Let $R$ and $\mathbb P_1$ be as in Lemma~\ref{al:11}. By Lemmas~\ref{al:10} and~\ref{al:11}, the minimax principle \cite[Theorem~2.9]{willem1996minimax}, applied to the class $\mathbb P_1$ of paths joining $-R\phi_{1,p}$ to $R\phi_{1,p}$, shows that $c=\inf_{P\in\mathbb P_1}\max_{u\in P[-1,+1]}\mathcal E(u)$ is a critical value of $\mathcal E$. A critical point at this level solves \eqref{a:57}. \qed

When $V_D=\emptyset$ the hypotheses \eqref{a:59} and \eqref{a:61} take a special form. There $\lambda_{1,p}(\Gamma,V_D)=0$, and by Proposition \ref{ap:8} and Proposition \ref{at:6} the curve $\mathfrak{C}$ satisfies $\alpha(b)>\bar{\lambda}$ for every $b\ge0$, with $\alpha(b)\to\bar{\lambda}$ as $b\to+\infty$, while $b+\alpha(b)$ increases to $+\infty$. The region of admissible pairs is therefore a half-strip, and no bound at all is imposed on the growth of $f$ at $+\infty$.

\begin{corollary}\label{ac:strip}
Let $V_D=\emptyset$ and let $f$ satisfy the $L^\infty$-Carath\'eodory conditions together with \eqref{a:58} and \eqref{a:60}, where $\pi_{\pm},\Pi_{\pm},\theta_{\pm},\Theta_{\pm}\in L^\infty(\Gamma)$ satisfy
\begin{align*}
0\le\pi_{+}\le\Pi_{+},\qquad 0\le\pi_{-}\le\Pi_{-}\le\bar{\lambda},\qquad 0\le\theta_{+}\le\Theta_{+},\qquad 0\le\theta_{-}\le\Theta_{-}\le\bar{\lambda}\quad\text{a.e. in }\Gamma,
\end{align*}
with $\theta_{+}>0$ and $\theta_{-}>0$ on subsets of positive measure. Then \eqref{a:57} admits at least one solution, whatever the size of $\Pi_{+}$ and $\Theta_{+}$.
\end{corollary}

\begin{proof}
By Propositions \ref{ap:5} and \ref{ap:6} (valid for $V_D=\emptyset$ by Proposition~\ref{ap:100}) the map $b\mapsto b+\alpha(b)$ is continuous and strictly increasing, and since $\alpha(0)=\lambda_{2,p}(\Gamma,\emptyset)$ (Corollary~\ref{ar:2} and Proposition~\ref{ap:100}) and $b+\alpha(b)\ge b$, it maps $[0,\infty)$ onto $[\lambda_{2,p}(\Gamma,\emptyset),\infty)$. Choose $b\ge0$ with
\begin{align*}
b+\alpha(b)\ \ge\ \max\bigl\{\|\Pi_{+}\|_{\infty},\ \|\Theta_{+}\|_{\infty}\bigr\},
\end{align*}
and set $(\lambda^+,\lambda^-)=(b+\alpha(b),\alpha(b))\in\mathfrak{C}$. Then $\Pi_{+}\le\lambda^+$ and $\Theta_{+}\le\lambda^+$ a.e., while Proposition \ref{at:6} gives $\Pi_{-}\le\bar{\lambda}<\alpha(b)=\lambda^-$ and $\Theta_{-}\le\bar{\lambda}<\lambda^-$ a.e. Hence \eqref{a:59} and \eqref{a:61} hold, the last requirement of \eqref{a:61} being satisfied on the negative side, and Theorem \ref{at:5} applies. \qed
\end{proof}

Corollary~\ref{ac:strip} is effective as soon as $\bar\lambda$ is known. Example~\ref{aex:2} gives it on an interval, where the half-strip also follows from \cite{defigueiredo1994first,arias2000antimaximum}, and Proposition~\ref{aex:100} gives it on the equilateral star.

\section*{Declarations}

\paragraph{Conflict of interest.} The authors declare that they have no conflict of interest.

\paragraph{Data availability.} Data sharing is not applicable to this article, as no datasets were generated  during the current study.
\paragraph{Funding.} The authors received no specific funding for this work.

\bibliographystyle{siam}
\bibliography{ref}

\end{document}